%% file: Dmods.tex
\documentclass[11pt]{article}
\usepackage{amsmath,amssymb,xypic,theorem,longtable}
\xyoption{curve}

\theoremstyle{change}
\theorembodyfont{\itshape}
\newtheorem{thm}[subsubsection]{Theorem}
\newtheorem{prop}[subsubsection]{Proposition}
\newtheorem{lemma}[subsubsection]{Lemma}
\newtheorem{cor}[subsubsection]{Corollary}
\newtheorem{defn}[subsubsection]{Definition}
{\theorembodyfont{\rmfamily}
\newtheorem{example}[subsubsection]{Example}
\newtheorem{remark}[subsubsection]{Remark}}

\makeatletter
\renewcommand{\subsection}{\@startsection%
{subsection}{2}{0mm}{\baselineskip}{-1em}%
{\normalfont\normalsize\bfseries}}
\renewcommand{\subsubsection}{\@startsection%
{subsubsection}{3}{0mm}{\baselineskip}{-1em}%
{\normalfont\normalsize\textit}}
\makeatother

\numberwithin{equation}{subsubsection}

\newcommand{\demobox}{\vrule height6pt width6pt depth0pt}

\newcommand{\nodemo}{\unskip\nobreak\hfil\qquad
\demobox\parfillskip=0pt\par}

\newenvironment{demo}{\noindent{\it Proof.}}
{{\unskip\nobreak\hfil\qquad
\demobox\parfillskip=0pt\par}
\medskip}

\include{mathdefs} 
\newcommand\hP{\hat P}
\newcommand\hI{\hat I}

\newcommand\fq{{\mathfrak{q}}}

\newcommand{\tensu}[1]{\underset{{#1}}{\otimes}}
\newcommand{\ctensu}[1]{\underset{{#1}}{\hat\otimes}}

\title{Arithmetic $\D$-modules on Locally Noetherian Formal Schemes}
\author{Richard Crew}
\date{}

\begin{document}
\maketitle

\begin{abstract}
  We extend Berthelot's theory of arithmetic $\D$-modules to a class
  of morphisms that are not necessarily of finite type, or even
  adic. This is the class of \emph{quasi-smooth} morphisms, which are
  formally smooth and \emph{universally noetherian} (any base change
  by a noetherian formal scheme is a morphism of noetherian formal
  schemes). We also work out the correspondence between modules over
  the various arithmetic differential operators of finite level, and
  the corresponding notions of stratification. Finally we prove a
  general theorem on Frobenius descent, proven by Berthelot in the
  smooth case. 
\end{abstract}

\input Dmods-intro

\noindent\textit{Notation and Conventions.} 
Rings are assumed to have an identity. Terminology and notation
regarding commutative algebra and formal schemes generally follows EGA
\cite{EGA}. For example if $A$ is a commutative ring and $M$ is an
$A$-module an $f\in A$, $M_f$ is the (algebraic) localization of $M$,
and if $M$ is a topological $A$-module $M_{\{f\}}$ is the formal
localization, i.e. the completion of $M_f$. A topological ring $R$ is
\textit{preadic} if it has the $J$-adic topology for some ideal
$J\sset R$, and \textit{adic} if it is preadic, separated and
complete.

Starting with section \ref{sec:smooth-morphisms} we will assume,
without explicit statement to the contrary that formal schemes are
\textit{locally noetherian}, i.e. locally of the form $\Spf{A}$ with
$A$ an adic noetherian ring. Even so, this condition will frequently
be stated explicitly for emphasis.

In any category $\cC$ with fibered products the notation $X_S(r)$
denotes the fibered product of $r+1$ copies of $X$ over $S$. If
$K\subset\bN$ is a finite subset we write $X_S(K)=X_S(r)$ where
$r=\#K$ is the cardinality of $K$, and the factors in the product
$X_S(K)$ are indexed by the elements of $K$ in increasing
order. Recall the usual notation for the projection maps relating the
$X_S(K)$: for $L\sset K$ the projection $p_{LK}:X_S(K)\to X_S(L)$ is
the product of the projections $p_i:X_S(K)\to X$ for all $i\in
L$. When $K=[0,1,\ldots,r]$ and $L=K\setminus\{i\}$ we will write
these projections as $p_L:X_S(r)\to X_S(r-1)$. For example, if $\cC$
is the category of schemes and $T$ is any scheme,
$p_{02}:X_S(2)\to X_S(1)$ is the morphism inducing
$(x,y,z)\mapsto (x,z)$ on $T$-valued points.

We use the same notation for tensor products of rings, or completed
tensor products of topological rings; thus if the category $\cC$ in
the last paragraph is the category of schemes and $X=\Sp{A}$,
$S=\Sp{R}$ then $X_S(r)=\Sp{A_R(r)}$. The tensor product of an
abelian group $M$ with $\bQ$ will usually be written $M_\bQ$, as in
\cite{berthelot:1996}.

When dealing with rings or geometric constructions in an affine
setting, completions will usually be denoted by a ``hat'' which is
dropped in purely geometric situations. For example if $A$ is an
$R$-algebra, $\diff_{A/R}$ is the usual module of 1-forms,
$\cdiff_{A/R}$ is its completion in the natural topology, but the
sheafification of $\cdiff_{A/R}$ for a morphism $\cX\to\cS$ is
$\diff_{\cX/\cS}$. Exceptions to this convention occur when a formal
construction is subject to further completion, as for example the ring
$\niv{\hD}{m}_{\cX/\cS}$, which is a completion of
$\niv{\D}{m}_{\cX/\cS}$ (in any case this notation is completely
entrenched in the literature). Another example: $\hat\cX_\cS(r)$ is
the formal completion of the product $\cX_\cS(r)$ along the diagonal.

In addition to the standard notations for multi-indices, we use the
following: for $K=(k_1,\ldots,k_n)\in\bZ^n$ and $a\in\bZ$ we write
$K<a$ (resp. $K\le a$) to mean $k_i<a$ (resp. $k_i\le a$) for all $i$.

\tableofcontents
\bigskip

\setlongtables
\begin{longtable}{ll}
  \caption*{\textbf{Index of Terminology.}}
  \endhead
  formally smooth, formally \'etale, formally unramified
  & defn.\  \ref{defn:formally-smooth}\\

  universally locally noetherian
  & defns.\  \ref{defn:universally-noetherian},
  \ref{defn:universally-noetherian-for-schemes}\\

  universally noetherian
  & defns.\  \ref{defn:universally-noetherian},
  \ref{defn:universally-noetherian-for-schemes}\\

  smooth, \'etale, unramified
  & defn.\  \ref{defn:smooth}\\

  quasi-smooth, quasi-\'etale, quasi-unramified
  & defn.\  \ref{defn:quasi-smooth}\\

  bilateralising ideal
  &defn.\  \ref{defn:bilateralising-ideal}\\

  $m$-bilateralising ideal
  &defn.\  \ref{defn:m-bilateralising}\\

  pro-quasicoherent $\O_\cX$-module
  &defn.\  \ref{defn:pro-qcoh}\\

  split-regular ideal
  &defn.\  \ref{defn:split-regular-ideal}\\

  open ideal adapted to split-regular ideal
  &defn.\  \ref{defn:adapted-to-split-regular-ideal}\\

  split-regular ideal sheaf
  &defn.\  \ref{defn:split-regular-sheaf}\\

  toplogically quasi-nilpotent
  &defns.\ \ref{defn:topologically-quasi-nilpotent},
  \ref{defn:topologically-quasi-nilpotent-isogeny}\\

\end{longtable}

\input Dmods1

\input Dmods2

\bibliographystyle{hplain}
\bibliography{Dmods}

\parindent=0pt

Email: \texttt{rcrew@ufl.edu}
\bigskip

Department of Mathematics\\
358 Little Hall\\
Gainesville FL 32611\\
USA

\end{document}

%% file: mathdefs.tex
\newcommand{\Ker}{{\operatorname{Ker}}}

\newcommand{\Sp}[1]{{\operatorname{Spec}({#1})}}

\newcommand{\Spalg}[2]{{\operatorname{Spec}_{#1}({#2})}}
\newcommand{\Spf}[1]{{\operatorname{Spf}({#1})}}
\newcommand{\Spfalg}[2]{{\mathit{Spf}_{#1}({#2})}}

\newcommand{\Hom}{{\operatorname{Hom}}}

\newcommand{\Tor}{{\operatorname{Tor}}}

\newcommand{\Diff}{{\operatorname{Diff}}}

\newcommand{\Der}{{\operatorname{Der}}}

\newcommand{\Exalcotop}{{\operatorname{Exalcotop}}}

\newcommand{\limdir}{\varinjlim}
\newcommand{\liminv}{\varprojlim}

\newcommand{\diff}{\Omega^1}
\newcommand{\cdiff}{\hat\Omega^1}

\newcommand{\tens}{\otimes}

\newcommand{\ctens}{\mathop{\hat\otimes}}

\newcommand{\fa}{{\mathfrak a}}
\newcommand{\fb}{{\mathfrak b}}

\newcommand{\fm}{{\mathfrak m}}
\newcommand{\fn}{{\mathfrak n}}
\newcommand{\fp}{{\mathfrak p}}
\newcommand{\bF}{{\mathbb{F}}}
\newcommand{\bA}{{\mathbb{A}}}

\newcommand{\bN}{{\mathbb N}}

\newcommand{\bQ}{{\mathbb Q}}

\newcommand{\bZ}{{\mathbb Z}}

\newcommand{\Sym}{{\textrm{Sym}}}

\newcommand{\fdim}[1]{\mathop{\textrm{fdim}({#1})}}

\renewcommand{\O}{{\cal O}}

\newcommand{\cB}{{\cal B}}
\newcommand{\cC}{{\cal C}}

\newcommand{\cI}{{\cal I}}
\newcommand{\cJ}{{\cal J}}
\newcommand{\cK}{{\cal K}}
\newcommand{\cM}{{\cal M}}

\newcommand{\cP}{{\cal P}}
\newcommand{\cQ}{{\cal Q}}
\newcommand{\cS}{{\cal S}}
\newcommand{\cT}{{\cal T}}
\newcommand{\cU}{{\cal U}}
\newcommand{\cV}{{\cal V}}

\newcommand{\cX}{{\cal X}}
\newcommand{\cY}{{\cal Y}}

\newcommand{\Xto}{\xrightarrow}

\newcommand{\inj}{\hookrightarrow}

\newcommand{\isom}{\Xto{\sim}}
\newcommand{\sset}{\subseteq}
\newcommand{\<}{\langle}
\renewcommand{\>}{\rangle}
\newcommand{\md}{\mathrm{d}}
\renewcommand{\d}{{\partial}}
\newcommand{\dpe}[2]{{#1}^{[{#2}]}}
\newcommand{\dpniv}[3]{{#1}^{\{{#2}\}_{#3}}}
\newcommand{\dpbrniv}[3]{{#1}^{\{{#2}\}_{(#3)}}}
\newcommand{\dpbrshort}[2]{{#1}^{\{{#2}\}}}
\newcommand{\dpabniv}[3]{{#1}^{\<{#2}\>_{(#3)}}}

\newcommand{\niv}[2]{{#1}^{({#2})}}
\newcommand{\bbinom}[3]{\genfrac{\{}{\}}{0pt}{}{#1}{#2}_{(#3)}}
\newcommand{\abinom}[3]{\genfrac{\<}{\>}{0pt}{}{#1}{#2}_{(#3)}}

\newcommand{\D}{{\cal D}}
\newcommand{\hD}{\hat{\cal D}}
\newcommand{\Ddag}{{\cal D}^\dagger}

\newcommand{\canj}{\circ}
\newcommand{\cani}{\bullet}



%% file: Dmods-intro.tex
\section*{Introduction}
\label{sec:intro}

The subject of this paper is an extension of Pierre Berthelot's theory
of arithmetic differential operators \cite{berthelot:1990},
\cite{berthelot:1996}, \cite{berthelot:2000},
\cite{berthelot:2002}. The original setting of the theory was that of
a formal scheme $\cX$ smooth over a noetherian $p$-adic formal scheme
$\cS$. The main constructions of the theory are of a sheaf of rings
$\Ddag_{\cX/\cS\bQ}$ on $\cX$ and the category of holonomic
$F\Ddag$-complexes, together with a set of cohomological operations
which are conjectured to satisfy Grothendieck's six-functor
formalism. These conjectures are still open, but Caro and Tsuzuki have
shown \cite{caro-tsuzuki:2012} that a full subcategory of the category
of holonomic $F\Ddag$-complexes, Caro's category of overholonomic
complexes is stable under the six operations and contains anything of
geometric origin. Nonetheless it can be difficult to show that any
particular object is overholonomic if it does not evidently come from
geometry, so there is some interest in Berthelot's original project
and in particular in studying the effect of cohomological operations
relative to closed immersion.

In this paper we introduce and study the class of \emph{universally
  noetherian} morphisms of formal schemes. This class includes
morphisms of finite type, but also morphisms that are not even adic;
for example $\Spf{\bZ_p[[X]]}\to\Spf{\bZ_p}$ is universally noetherian
when both rings have the topology defined by the maximal ideal. The
construction of ordinary and arithmetic differential operator rings
can be carried out for any morphism that is \emph{quasi-smooth}, i.e.\
universally noetherian and formally smooth.

There are two motivations for doing this. The first is that an
arithmetic differential operator ring relative to
$\Spf{\bZ_p[[X]]}\to\Spf{\bZ_p}$ appeared in the author's study
\cite{crew:2006} of cohomological operations relative to the inclusion
of a point in a smooth curve, and the study of local monodromy of
holonomic $F\Ddag$-modules in \cite{crew:2012}. The results of this
paper are meant to prepare the way for work in higher dimensions and
codimensions. The second is that such morphisms are the natural
setting for the categories of convergent and overconvergent
isocrystals on a scheme of finite type over a field of positive
characteristic. The methods of this paper allow one to construct such
categories for any quasi-smooth morphism $\cX\to\cS$ of noetherian
formal schemes; this will be the subject of a later publication. We
should note that Bernard Le Stum, in some recent preprints has also
constructed categories of this sort, though the relation of his
constructions with arithmetic differential operator rings is not clear
at the moment.

After some brief remarks on flat and formally smooth morphisms of
formal schemes, \S\ref{sec:formal-schemes} we introduce our basic
finiteness condition: a morphism $\cX\to\cS$ of noetherian formal
schemes is \emph{universally noetherian} if for every morphism
$\cT\to\cS$ with $\cT$ a noetherian formal scheme, $\cT\times_\cS\cX$
is noetherian. One can make a similar definition by replacing
``noetherian'' by ``locally noetherian'' but we will not use this
notion in this paper. This category has nice stablility properties: it
is stable under composition, pullbacks and fiber products. Morphisms
of finite type are universally noetherian. The inclusion
$\hat\cX_Y\to\cX$ of the completion of a noetherian formal scheme
along a closed subscheme is universally noetherian. Localizations are
universally noetherian. We can thus generate a large class of examples
of these morphisms.

If $\cX\to\cS$ is universally noetherian the ideal of the diagonal
$\cX\to\cX\times_\cS\cX$ is locally finitely generated. Thus if
$\cX/\cS$ is \emph{quasi-smooth}, i.e.\ formally of finite type and
universally noetherian, the construction of the usual differential
invariants of $\cX/\cS$ proceeds without essential changes; this we do
in \S\ref{sec:smooth-morphisms}. The usual characterizations of
smooth, \'etale and unramified morphisms has an analogue in the more
general case, as does first-order deformation theory. However when in
\S\ref{sec:arith-diff} we turn to the construction of Berthelot's
construction of arithmetic differential operators we run into the
following problem. In \cite{berthelot:1996} the ring
$\niv\hD m_{\cX/\cS}$ is the $p$-adic completion of
$\niv\D m_{\cX/\cS}$. In the general case the $p$-adic completion is
not what we want since $\niv\D m_{\cX/\cS}/p^n\niv\D m_{\cX/\cS}$ is
not a quasi-coherent sheaf on a scheme (this fact is used in
\cite{berthelot:1996}, for example to prove Cartan's theorem A and B
for $\niv\hD m_{\cX/\cS}$). What we want is to take a completion with
respect to an ideal of definition $J\subset\O_\cX$, but
$\niv\D m_{\cX/\cS}/J^n\niv\D m_{\cX/\cS}$ is not a ring unless the
$J$ is carefully chosen: we want the left or right ideal generated by
$J$ in $\niv\D m_{\cX/\cS}$ to be a bilateral ideal. Fortunately there
is a large supply of such ideals, which we call
\emph{$m$-bilateralsing} ideals. A similar problem turns up in section
\S\ref{sec:stratifications} when we want to construct $m$-PD-envelopes
of ideals in a formal setting.

In \S\ref{sec:stratifications} we develop the theory of
$m$-PD-stratifications and $m$-HPD-stra\-ti\-fications in this
context. Again this works as expected provided we place an additional
technical condition on the ideal whose $m$-PD-envelopes we wish to
construct (that of being split-regular, definition
\ref{defn:split-regular-ideal}). In particular, the category of
quasi-nilpotent $\niv\hD m_{\cX/\cS}$-modules is equivalent to the
category of $\O_\cX$-modules with an $m$-HPD-stratification. Finally
if $\cB$ is a commutative $\O_\cX$-algebra with a compatible left
$\niv\D m_{\cX/\cS}$-module structure (the multiplication
$\cB\tens_{\O_\cX}\cB\to\cB$ is linear), the above theory extends to
modules over the ring
$\niv\D m_{\cB/\cS}=\cB\tens_{\O_\cX}\niv\D m_{\cX/\cS}$. This
construction is the basis of our treatment of (over)convergent
isocrystals in the sequel. We end with Berthelot's Frobenius descent
theorem in this setting.

\bigskip

\noindent\textit{Acknowledgements.} 
I am indebted to many people for helpful conversations and a large
number of pointed questions that served to improve the originial
manuscript. I am particularly grateful for help from Pierre Berthelot
and Bernard Le Stum in Rennes, Francesco Baldassari and Bruno
Chiarellotto in Padua, Chris Lazda in Exeter and Tomoyuki Abe in
Tokyo. Most of this article was worked out during visits to IRMAR at
the University of Rennes and to the mathematics department of the
University of Padua. Much of the final draft was written during a
visit to the Kavli IPMU in Tokyo in 2024.  I would like to thank all
of these institutions for their support. Finally, words can not
adequately express what Pierre Berthelot and his work meant for all
those who are working in this subject. I offer this paper as a small
tribute to his memory.

\bigskip

%% file: Dmods1.tex
\section{Formal Geometry}
\label{sec:formal-schemes}

\subsection{Flatness and formal smoothness.}
\label{sec:formal-flatness}

We will use the same definition of formal smoothness for a morphism of
formal schemes as for ordinary schemes: 

\begin{defn}\label{defn:formally-smooth}
  a morphism $\cX\to\cS$ of formal schemes is \emph{formally smooth}
  (resp. \emph{formally unramified}, \emph{formally \'etale}) if for
  any commutative square
  \begin{displaymath}
    \xymatrix{
      Z_0\ar[r]\ar[d]&\cX\ar[d]\\
      Z\ar[r]\ar@{.>}[ur]&\cS
    }  
  \end{displaymath}
  in which $Z$ is an \emph{affine scheme} and $Z_0\to Z$ is a
  nilpotent closed immersion, there exists a morphism $Z\to\cX$
  (resp. there is at most one morphism, there exists a unique morphism)
  making the extended diagram commutative.   
\end{defn}

When $\cX=\Spf{B}$ and $\cS=\Spf{A}$ are affine $\cX\to\cS$ is
formally smooth if and only if $B$ is a formally smooth $A$-algebra in
the sense of \cite[IV, 14.3.1]{EGA}. The reader may check that most of
the elementary properties of formally smooth, formally unramified and
formally \'etale morphisms of ordinary schemes (e.g. \cite[IV
\S17]{EGA} propositions 17.1.3 and 17.1.4) are also valid in the
present context of locally noetherian formal schemes. There is
one important exception: with this definition, formally smoothness is
\textit{not} a local property, either on the base or the source. We
will return to this question in section \ref{sec:smooth-formal-case}.
For now we observe that if $f:\cX\to\cS$ is formally smooth
(resp. unramified, \'etale) and $U\sset\cX$, $V\sset\cS$ are open
formal subschemes such that $f(U)\sset V$ , the induced morphisms
$U\to V$ is formally smooth (resp. unramified, \'etale); this follows
from the definitions and the fact that $Z$ and $Z_0$ have the same
underlying topological space.

\begin{lemma}\label{lemma:completions-of-local-rings}
  Let $A$ be an adic noetherian ring with ideal of definition $J$,
  $x\in\Spf{A}$, $\fm\subset A$ the open prime ideal corresponding to
  $x$. If $M$ is a coherent $A$-module let $\cM$ be the sheaf on
  $\Spf{A}$ corresponding to $M$ and $\cM_x$ the stalk of $\cM$ at
  $x$. Then the natural morphism $M_\fm\to\cM_x$ induces an
  isomorphism $\hat M_\fm\isom\hat\cM_x$ of the $J$-adic completions.
\end{lemma}
\begin{demo}
  We have $M_\fm=\limdir_fM_f$ and $\cM_x=\limdir_fM_{\{f\}}$ where
  $f$ runs through $A\setminus\fm$, and the natural map
  $M_\fm\to\cM_x$ is induced by $M_f\to M_{\{f\}}$. It suffices to
  show that
  \begin{displaymath}
    (\limdir_fM_f)\tens_AA/J^n\to(\limdir_fM_{\{f\}})\tens_AA/J^n
  \end{displaymath}
  is an isomorphism for all $n$. Since inductive limits commute with
  tensor products, this is the same as
  \begin{displaymath}
    \limdir_f(M_f\tens_AA/J^n)\to\limdir_f(M_{\{f\}}\tens_AA/J^n)
  \end{displaymath}
  and the assertion is clear, since
  $M_f\tens_AA/J^n\to M_{\{f\}}\tens_AA/J^n$ is an isomorphism.
\end{demo}

The morphism $M_\fm\to\cM_x$ is functorial in $M$, and also in $A$ in
the sense that if $A\to B$ is a continuous homomorphism of adic rings
yielding $f:\Spf{B}\to\Spf{A}$, and $\fn\subset B$ is an open prime
ideal corresponding to $y\in\Spf{B}$ such that $f(y)=x$, the diagram
\begin{displaymath}
  \xymatrix{
    (B\tens_AM)_\fn\ar[r]&f^*\cM_y\\
    M_\fm\ar[r]\ar[u]&\cM_x\ar[u]
  }
\end{displaymath}
is commutative. It follows that the isomorphisms of completions is
functorial in the same sense.

Recall that a morphism $f:\cX\to\cS$ of locally ringed spaces, and in
particular of locally noetherian formal schemes is \textit{flat} at a
point $x\in\cX$ if the morphism $\O_{f(x)}\to\O_x$ of local rings is
flat, and $f$ is flat if it is flat at every point of $x$.

\begin{lemma}\label{lemma:formal-flatness}
  For any morphism $f:\cX\to\cS$ of locally noetherian formal
  schemes, the following are equivalent:
  \begin{enumerate}
  \item $f$ is flat;
  \item $f$ is flat at every closed point of $\cX$;
  \item for every pair $\Spf{B}\sset\cX$, $\Spf{A}\sset\cS$ of open
    affines such that $f(\Spf{B})\sset\Spf{A}$, $B$ is a
    flat $A$-algebra.
  \end{enumerate}
\end{lemma}
\begin{demo}
  The implications (iii)$\implies$(i) and
  \ref{lemma:formal-flatness}.1$\implies$(ii) are clear. Suppose now
  that (ii) holds; we can also assume that $A$ and $B$ in (iii) are
  noetherian. It suffices to show that for every maximal ideal
  $\fn\subset B$ and $\fm$ equal to the inverse image of $\fn$ in $A$,
  $B_\fn$ is a flat $A_\fm$-algebra \cite[Ch. II \S3 no. 4
  Prop. 15]{bourbaki-AC}. By the faithful flatness of completions
  (e.g. \cite[$0_I$ 7.6.18]{EGA}) it suffices to show that
  $\hat B_\fn$ is a flat $\hat A_\fm$-algebra, where the completions
  are taken with respect to the adic topologies of $A$ and $B$.

  Since $B$ is adic and noetherian it is a Zariski ring, and every
  maximal ideal is open. Therefore $\fn\subset B$ corresponds to a
  closed point $x\in\cX$ and $y=f(x)$ corresponds to $\fm$. Now
  (ii) asserts $\O_x$ is a flat
  $\O_y$-algebra, and thus $\hat\O_x$ is a flat $\hat\O_y$-algebra. By
  lemma \ref{lemma:completions-of-local-rings} there are isomorphisms
  $\hat\O_x\simeq\hat B_\fn$ and $\hat\O_y\simeq\hat A_\fm$, and by
  functoriality the map $\hat A_\fm\to\hat B_\fn$ corresponds via
  these identifications to $\hat\O_y\to\hat\O_x$.
\end{demo}

\begin{prop}\label{prop:formally-smooth-implies-flat}
  A formally smooth morphism $f:\cX\to\cS$ of locally noetherian
  formal schemes is flat.
\end{prop}
\begin{demo}
  It suffices to show that for all $x\in\cX$ and $y=f(x)$ that $\O_x$
  is a flat $\O_y$-algebra. Pick open affine neighborhoods $x\in U$,
  $y\in V$ such that $f(U)\sset V$; by the remark just before the
  proposition, the morphism $U\to V$ induced by $f$ is formally
  smooth. If $U=\Spf{B}$ and $V=\Spf{A}$, the topological $A$-algebra
  $B$ is formally smooth. Then $\O_x$ is a formally smooth
  $\O_y$-algebra for the preadic topologies induced by $B$ and $A$,
  and therefore formally smooth for the preadic topologies defined by
  the maximal ideals of $\O_x$ and $\O_y$. The assertion then follows
  from theorem 19.7.1 of \cite[$0_{IV}$]{EGA}.
\end{demo}

\begin{cor}\label{cor:formally-smooth-implies-flat}
  If $A$ is an adic noetherian ring and $B$ is an adic noetherian ring
  and a formally smooth $A$-algebra, $B$ is flat over $A$.
\end{cor}

\subsection{Universally noetherian morphisms.}
\label{sec:formal-schemes-finiteness}

The basic finiteness condition of this paper is the following:

\begin{defn}\label{defn:universally-noetherian}
  A morphism $f:\cX\to\cS$ of locally noetherian formal schemes is
  \emph{universally noetherian} (resp. \emph{universally locally
    noetherian} if for every noetherian formal scheme $\cS'$ and every
  morphism $g:\cS'\to\cS$ of formal schemes, the fiber product
  $\cX\times_\cS\cS'$ is noetherian (resp. locally noetherian).
\end{defn}

We will mainly be concerned with universally noetherian morphisms, and
leave it to the reader to formulate in what follows the corresponding
results for unviersally locally noetherian morphisms. In the situation
of definition \ref{defn:universally-noetherian} we will also say that
$\cX$ is a universally noetherian formal $\cS$-scheme.  If $A$ and $B$
are adic noetherian rings and $A\to B$ is a continuous homomorphism,
we say that $B$ is a universally noetherian $A$-algebra if $\Spf{B}$
is a universally noetherian formal $\Spf{A}$-scheme. To check that
$\cX\to\cS$ is universally noetherian, it suffices to check the
condition of \ref{defn:universally-noetherian} for morphisms
$\cS'\to\cS$ with $\cS'$ formally affine and noetherian. In particular
if $A$ is a noetherian ring, an $A$-algebra $B$ is universally
noetherian if and only if $A\ctens_BC$ is a noetherian ring for any
adic noetherian $B$-algebra $C$.


It is immediate from the definition that a universally noetherian
morphism is quasicompact, and that if $f:\cX\to\cS$ is universally
noetherian and $\cS$ is noetherian, then so is $\cX$. 

\begin{lemma}
  If $f:\cX\to\cS$ is universally noetherian and $\cS'\to\cS$ is morphism with
  $\cS'$ locally noetherian, the fiber product $\cX\times_\cS\cS'$ is
  locally noetherian.
\end{lemma}
\begin{demo}
  Let $\{U_\alpha\}_{\alpha\in I}$ be a cover of $\cS'$ by noetherian
  formal schemes. By definition the
  $f^{-1}(U_\alpha)=\cX\times_\cS U_\alpha$ are noetherian, and since
  they cover $\cX$, $\cX$ is locally noetherian.
\end{demo}

Thus fibered products with a universally noetherian morphism do not
force us to leave the category of adic locally noetherian schemes.
The usual \textit{sorites} hold for the class of universally
noetherian morphisms:

\begin{prop}\label{prop:noeth-sorites}
  \begin{enumerate}
  \item An immersion is universally noetherian.  
  \item Let $f:\cX\to\cS$ be a morphism of locally noetherian schemes.
    If $f:\cX\to\cS$ is universally noetherian and $\cS'$ is any
    locally noetherian formal $\cS$-scheme, the base-change
    $\cX\times_\cS\cS'\to\cS'$ is universally noetherian.
  \item If $f:\cX\to\cS$ and $g:\cY\to\cX$ then $f\circ g:\cY\to\cS$
    is universally noetherian.
  \item If $f:\cX\to\cS$ and $g:\cY\to\cS$ are universally noetherian
    morphisms then so is $f\times g:\cX\times_\cS\cY\to\cS$.
  \end{enumerate}
\end{prop}
\begin{demo}
  For (i) it suffices to treat the case of open and closed
  immersions. Since any base-change of an open (resp. closed)
  immersion is open (resp. closed), the assertion is clear in the case
  of closed immersions, and for open immersions it suffices to add
  that an open immersion of locally noetherian formal schemes is
  quasi-compact.  Assertions (ii) and (iii) follow from the definition
  and the transitivity of fibered products, while (iv) follows from
  (i) and (ii).
\end{demo}

We will show that (iii) in the proposition has a partial converse
(proposition \ref{prop:properties-of-loc-noeth2}).  The property of
being universally noetherian is local on the base and, with a suitable
restriction, on the source:

\begin{prop}\label{prop:properties-of-loc-noeth}
  Let $f:\cX\to\cS$ be a morphism of locally noetherian schemes.
  \begin{enumerate}
  \item If $f$ is quasi-compact and $\{U_\alpha\}_{\alpha\in I}$ is an
    open cover of of $\cX$, then $f$ is universally noetherian if and
    only if each of the induced morphisms $U_\alpha\to\cS$ is
    universally noetherian.
  \item If $\{V_\alpha\}_{\alpha\in I}$ is an open cover of $\cS$,
    then $f$ is universally noetherian if and only if the morphisms
    $f^{-1}(V_\alpha)\to V_\alpha$ are universally noetherian.
  \end{enumerate}
\end{prop}
\begin{demo}
  Necessity in assertions (i) and (ii) follows from the sorites
  \ref{prop:noeth-sorites}. To prove the condition is sufficient in
  (i) we first observe that for any $\cS'\to\cS$ with $\cS'$
  noetherian, each of the $U_\alpha\times_\cS\cS'$ are
  noetherian. Since the $U_\alpha\times_\cS\cS'$ cover
  $\cX\times_\cS\cS'$, the latter is locally noetherian. On the other
  hand since $\cX\to\cS$ is quasicompact, so is
  $\cX\times_\cS\cS'\to\cS'$, from which it follows that
  $\cX\times_\cS\cS'$ is quasicompact, hence noetherian. This proves
  (i), and sufficiency in (ii) follows from (i), as one sees by taking
  $U_\alpha=f^{-1}(V_\alpha)$ and observing, first, that $f$ is
  necessarily quasi-compact since all of the
  $f^{-1}(V_\alpha)\to V_\alpha$ are, and second that the composite
  morphisms $U_\alpha\to V_\alpha\to\cS$ are universally noetherian.
\end{demo}

\begin{prop}\label{prop:properties-of-loc-noeth2}
  If $\cX\to\cS$ is universally noetherian and $\cY$ is a locally
  noetherian formal $\cS$-scheme, any morphism $\cX\to\cY$ is
  universally noetherian.
\end{prop}
\begin{demo}
  Suppose $\cT$ is a noetherian formal scheme and $\cT\to\cS$ is a
  morphism. By (ii) of the last proposition we may assume that $\cY$
  is affine, in which case $\cY\to\cS$ is separated and
  $\cX\times_\cY\cT\to\cX\times_\cS\cT$ is a closed immersion.  Since
  by hypothesis $\cX\times_\cS\cT$ is noetherian, it follows that
  $\cX\times_\cY\cT$ is noetherian as well, as required.
\end{demo}

\begin{defn}\label{defn:universally-noetherian-for-schemes}
  A morphism $f:X\to S$ of locally noetherian \emph{schemes} is
  universally noetherian if $X\times_SS'$ is noetherian for any
  morphism of \emph{schemes} $S'\to S$ with $S'$ noetherian.
\end{defn}

We say that an $A$-algebra $B$ is universally noetherian if
$\Sp{B}\to\Sp{A}$ is universally noetherian in the above sense. Such
$A$-algebras are sometimes called \textit{strongly noetherian} but I
avoid this terminology since it conflicts with a different notion
from the theory of adic spaces.


The condition of definition
\ref{defn:universally-noetherian-for-schemes} is \textit{a priori}
weaker than the condition that $f$ be universally noetherian in the
sense of definition \ref{defn:universally-noetherian} when $X$ and $S$
are regarded as (discrete) formal schemes. In fact the latter
condition requires that $X\times_S\cS'$ be noetherian for any morphism
$\cS'\to S$ with $\cS'$ a noetherian \textit{formal} scheme. In fact
\ref{defn:universally-noetherian} and
\ref{defn:universally-noetherian-for-schemes} are equivalent for a
morphism of schemes. To see this we will need a general result of
topological algebra, which combines \cite[Ch. III \S2 no. 12
Cor. 2]{bourbaki-AC} and \cite[Ch. III \S2 no. 10
Cor. 5]{bourbaki-AC}; see also the general discussion of completions
in \cite[III \S2 no.\ 12]{bourbaki-AC}. In the next proposition and
its corollary the ``natural'' topology of $\hat R$ is the topology it
has as the completion of the preadic topological ring $(R,I)$: the
ideals $\widehat{(I^n)}$ are a basis of the neighborhoods of $0$.  The
completion of an ideal $M\sset R$ will be denoted by $\hat M$, and
identified with an ideal of $\hat R$.

\begin{prop}\label{prop:general-noetherian-lemma}
  Suppose $R$ is a commutative ring and $I\subset R$ is a finitely
  generated ideal. Then for the $I$-preadic topology of $R$,
  \begin{equation}\label{eq:general-noetherian-lemma}
    \widehat{(I^n)}=(\hat I)^n=I^n\hat R.
  \end{equation}
  Furthermore the natural topology of $\hat R$ is the $\hat I$-adic
  topology and the natural homomorphism $R/I^n\to\hat R/\hat I^n$ is
  an isomorphism. Finally if $R/I$ is noetherian then so is $\hat R$,
  and in particular $\hat R$ is a Zariski ring.
\end{prop}

\begin{cor}\label{cor:completions-of-powers}
  Let $(R,J)$ be a preadic ring with $J$ finitely generated. If $R/J$
  is noetherian then for any ideal $M\subset R$,
  $\hat M^n=\widehat{M^n}$.
\end{cor}
\begin{demo}
  The proposition says that $\hat R$ is noetherian and therefore the
  ideal $\hat M^n\subset R$ is finitely generated, and thus closed
  since $\hat R$ is a Zariski ring. If $i:R\to\hat R$ is the canonical
  homomorphism, $i(M^n)\sset\hat M^n\sset\widehat{M^n}$ and $i(M^n)$
  is dense in $\widehat{M^n}$. Since both $\hat M^n$ and
  $\widehat{M^n}$ are closed, $\hat M^n=\widehat{M^n}$.
\end{demo}

\begin{lemma}\label{lemma:univ-noeth-ord-schemes}
  A universally noetherian morphism $X\to S$ of locally noetherian
  schemes is also universally noetherian when $X$ and $S$ are
  considered as formal schemes.
\end{lemma}
\begin{demo}
  It suffices to treat the affine case $X=\Sp{B}$, $S=\Sp{A}$. Let
  $(C,J)$ be an adic noetherian ring. Then $B\tens_AJ$ is an ideal of
  finite type in $B\tens_AC$ and an ideal of definition of
  $B\tens_AC$, and by hypothesis the ring
  $(B\tens_AC)/(B\tens_AJ)\simeq B\tens_A(C/J)$ is noetherian. The
  last proposition then shows that $B\ctens_AC$ is noetherian, as
  required.
\end{demo}

\begin{prop}\label{prop:noeth-if-red-is noeth}
  Suppose $f:\cX\to\cS$ is a morphism of locally noetherian formal
  schemes and $f_0:X\to S$ is the corresponding morphism of reduced
  closed subschemes. Then $f$ is universally noetherian if and only if
  $f_0$ is universally noetherian
\end{prop}
\begin{demo}
  Necessity: if $f$ is universally noetherian then so is
  $\cX\times_\cS S\to S$, and the closed immersion $X\to\cX\times_\cS
  S$ is universally noetherian as well.

  Sufficiency: it is enough to check the case where $\cX=\Spf{B}$ and
  $\cS=\Spf{A}$ are formally affine. Then $A$ and $B$ are noetherian,
  and if $I\subset A$, $J\subset B$ are maximal ideals of definition,
  $X=\Sp{B/J}$ and $S=\Sp{A/I}$. Suppose $B'$ is a noetherian
  $A$-algebra with ideal of definition $J'$ such that $IB'\sset J'$.
  By hypothesis the ring $(B/J)\tens_{A/I}(B'/J')$ is noetherian.  Set
  $B''=B\tens_AB'$ and $J''=J\tens B'+B\tens J'\subset B''$. By
  definition $B\ctens_AB'$ is the completion $\hat B''$ for the
  $J''$-adic topology. Since $J\subset B$ and $J'\subset B'$ are
  finitely generated, so is $J''$ and it follows from proposition
  \ref{prop:general-noetherian-lemma} that $\hat B''$,
  i.e. $B\ctens_AB'$ is noetherian, as required.
\end{demo}

Recall that $f$ is \textit{formally of finite type} if, in the
notation of the last proposition, $f_0$ is of finite type.

\begin{cor}\label{cor:fft-is-noeth}
  A morphism of formal schemes that is formally of finite type is
  universally noetherian.
\end{cor}
\begin{demo}
  Since a morphism of finite type is universally noetherian this
  follows from proposition \ref{prop:noeth-if-red-is noeth}.
\end{demo}

If $\cX$ is an locally noetherian formal scheme and $Y\subset\cX$ is a
closed subscheme, the completion $\hat\cX_Y$ of $\cX$ along $Y$ is
defined in the same way as in the case of ordinary schemes,
c.f. \cite[I Ch. 10]{EGA}.

\begin{cor}\label{cor:completion-is-noeth}
  If $f:\cX$ is a locally noetherian formal scheme and $Y\subset\cX$
  is a closed subscheme, the canonical morphism $i_Y:\hat\cX_Y\to\cX$
  is universally noetherian.
\end{cor}
\begin{demo}
  The morphism of reduced schemes induced by $i_Y$ is a closed
  immersion.
\end{demo}

From the corollary and the sorites we conclude that if $\cX\to\cS$ is
of finite type and $Y\subset\cX$ is closed, $\hat\cX_Y\to\cS$ is
universally noetherian; this is the main case of interest for the
construction of the categories of convergent and overconvergent
isocrystals. 

\begin{cor}\label{cor:localization-is noeth}
  Suppose $(A,I)$ is an adic noetherian ring, $S\subset A$ is a
  multiplicative system and $B$ is the completion of $S^{-1}A$ with
  respect to the ideal $S^{-1}I$. Then $B$ is a universally noetherian
  $A$-algebra. 
\end{cor}
\begin{demo}
  We may assume that $I$ is a maximal ideal of definition of $A$.
  Suppose $C$ is an adic noetherian ring with maximal ideal of
  definition $K$. If $A\to C$ is a continuous homomorphism then
  $CI\sset K$. The completion $J$ of $S^{-1}I$ is a maximal ideal of
  definition of $B$, $B/J\simeq S^{-1}A/S^{-1}I$ and the isomorphism
  $(B/J)\tens_{(A/I)}(C/K)\simeq S^{-1}C/S^{-1}K$ shows that
  $(B/J)\tens_{(A/I)}(C/K)$ is noetherian. Thus $B/J$ is a universally
  noetherian $A/I$-algebra, and it follows from proposition
  \ref{prop:noeth-if-red-is noeth} that $A\to B$ is universally
  noetherian.
\end{demo}

Applying the corollary when $S$ is the complement of an open prime
ideal, we get:

\begin{cor}\label{cor:immersion-of-complete-local-ring-is-noeth}
  Suppose $\cX$ is a locally noetherian formal scheme and $x$ is
  a point of $\cX$. Denote by $\hat\O_x$ the completion of the local
  ring of $x$ with respect to an ideal of definition of $\cX$. Then
  $\Spf{\hat\O_x}\to\cX$ is an adic universally noetherian morphism.
  \nodemo
\end{cor}

From corollary \ref{cor:completion-is-noeth} we see that
$\Spf{\hat\O_x}\to\cX$ is also universally noetherian if $\hat\O_x$ is
given the adic topology defined by the maximal ideal, although the
morphism $\Spf{\hat\O_x}\to\cX$ is not adic in this case.

The following proposition is an easy consequence of the fact that a
formal scheme has the same underlying topological space as its reduced
closed subscheme. Its equivalent properties define the notion of a
\textit{radicial} morphism of adic locally noetherian schemes.

\begin{prop}\label{prop:radicial-morphisms}
  For any universally noetherian morphism $f:\cY\to\cX$, the following
  are equivalent:
  \begin{enumerate}
  \item $f$ is universally injective, i.e. for any morphism
    $\cX'\to\cX$ with $\cX'$ locally noetherian,
    $\cY\times_\cX\cX'\to\cX'$ is injective.
  \item The morphism induced by $f$ on the reduced closed subschemes
    of $\cY$ and $\cX$ is radicial.
  \end{enumerate}\nodemo
\end{prop}

The following result of V\'amos \cite[Theorem 11]{vamos:1978} is
useful in constructing (counter)examples:

\begin{prop}\label{prop:vamos-thm}
  A field extension extension $L/K$ is finitely generated if and only
  if $L\tens_KL$ is noetherian.
\end{prop}

Since \cite[Theorem 11]{vamos:1978} the comes at the end of a long
discussion of more general topics we give a direct proof for the
reader's convenience.  The direct implication is clear since by
hypothesis $L$ is a localization of a finitely generated
$K$-algebra. Suppose conversely that $L\tens_KL$ is noetherian; the
reverse implication is proven in a sequence of cases:

\begin{enumerate}
\item \textit{$L/K$ inseparable algebraic}. In this case we show that
  $L/K$ is finite. Observe first that if $E/F$ is an extension of
  fields and $E\tens_FE$ is noetherian then $\diff_{E/F}$ is a
  finite-dimensional $F$-vector space. Let $p>0$ be the
  characteristic of $K$.
  \begin{enumerate}
  \item We first show that there is an $n\ge0$ such that
    $L^{p^n}\sset K$. If not, there is an infinite sequence
    $\{u_n\}_{n>0}$ of elements of $L$ such that $u_n^{p^n}\in K$ and
    $u_n^{p^{n-1}}\not\in K$. We claim that the $u_n$ are
    $p$-independent over $K$, i.e. that the monomials
    $\prod_{n\in S}u_n^{k(n)}$ are linearly independent over $K$,
    where $S$ runs through finite sets of positive integers and
    $0\le k(n)<p$.  Given this, theorem 21.4.5 of \cite[$0_{IV}$]{EGA}
    shows that the elements $du_n\in\diff_{L/K}$ for all $n>0$ are
    linearly independent over $L$. By the preceding observation this
    implies that $L\tens_KL$ is not noetherian, so $L^{p^n}\sset K$
    for $n\gg0$.

    We show that $u_1,\ldots,u_n$ are $p$-independent by induction on
    $n$. When $n=1$, the assertion is that
    $u_1,u_1^2,\ldots,u_1^{p-1}$ are linearly independent over $K$. If
    not, $u_1$ is separable over $K$, and since $u_1^p\in K$ we find
    that $u_1\in K$, a contradiction. Suppose now that
    $u_1,\ldots,u_{n-1}$ are $p$-independent and let
    $x=u_n^{p^{n-1}}$. A nontrivial $p$-dependence relation among
    $u_1,\ldots,u_n$ can be written $\sum_{i<p}a_iu_n^i=0$ with
    $a_i\in K(u_1,\ldots,u_{n-1})$ and not all of the $a_i$ are equal
    to zero since $u_1,\ldots,u_{n-1}$ are $p$-independent. Taking the
    $p^{n-1}$ power of this relation yields $\sum_{i<p}b_kx^i=0$ with
    $b_i=a_i^{p^{n-1}}\in K$ and not all of the $b_i$ are zero. Again
    this implies that $x$ is separable over $K$, and since
    $x^p=u_n^{p^n}\in K$ we again conclude that $x\in K$,
    i.e. $u_n^{p^{n-1}}\in K$, a contradiction.
  \item We now set $L_n=KL^{p^n}$, and show that $L_n/L_{n+1}$ is
    finite for all $n$. Then $L_0=L$ and by part (a) we know $L_n=K$
    for $n\gg0$ and, so this will show that $L/K$ is finite.

    Since $L\tens_KL$ is noetherian so is $L\tens_{L_1}L$, and by our
    earlier observation the $L$-space $\diff_{L/L_1}$ has finite
    dimension. By theorem 21.4.5 of \cite[$0_{IV}$]{EGA} $L/L_1$ has a
    finite $p$-base, so $L/L_1$ is finite since a $p$-base of $L/L_1$
    generates it as an extension. Since Frobenius is an injective
    homomorphism, $L^{p^n}$ is a finite extension of
    $L_1^{p^n}=K^{p^n}L^{p^{n+1}}$. We conclude that $L_n=KL^{p^n}$ is
    a finite extension of $L_{n+1}=KL^{p^{n+1}}$ for all $n\ge0$.
  \end{enumerate}
\item \textit{$L/K$ algebraic}. We again show that $L/K$ is finite.
  Let $K^s$ be the separable closure of $K$ in $L$, so that $L/K^s$ is
  purely inseparable. Since $L\tens_KL$ is noetherian so is
  $L\tens_{K^s}L$ and thus $L/K^s$ is finite by the last case.  If
  $K^s/K$ were not finite, $K^s\tens_KK^s$ would be an infinite direct
  sum of fields and the underlying topological space of
  $\Sp{K^s\tens_KK^s}$ would not be noetherian. Since
  $\Sp{L\tens_KL}\to\Sp{K^s\tens_KK^s}$ is faithfully flat the
  underlying space of $\Sp{L\tens_KL}$ would not be noetherian
  either. Therefore $K^s/K$ is finite, and so is $L/K$.
\item \textit{$L/K$ purely transcendental}. In this case we show that
  if $L/K$ is not finitely generated then the ring $L\tens_KL$ has an
  infinite strictly ascending chain of prime ideals.

  Suppose that $x_1$, $x_2$, $x_3\ldots$ is a sequence of
  algebraically independent elements and let
  $K_n=K[x_1,\ldots,x_n]$. The projection $f_n:L\tens_KL\to L\tens_{K_n}L$
  is surjective and $L\tens_{K_n}L$ is a domain, being the
  localization of a polynomial ring. Therefore $\fp_n=\Ker(f_n)$ is a
  prime ideal. Clearly $\fp_n\sset\fp_{n+1}$, and the inclusion is
  strict since $x_{n+1}\tens1-1\tens x_{n+1}$ is in $\fp_{n+1}$ but
  not in $\fp_n$.
\item \textit{General case}. Let $K\sset L_1\sset L$ be an
  intermediate field with $L_1/K$ purely transcendental and $L/L_1$
  algebraic. As before, $L\tens_{L_1}L$ is noetherian, so $L/L_1$ is
  finite and it suffices to show that $L_1/K$ is finitely
  generated. Suppose that it is not and let
  $\fp_0\subset \fp_1\subset \fp_2\subset\cdots$ be the chain of
  ideals of $L_1\tens_KL_1$ constructed in the last step. Since
  $L/L_1$ is finite, $L\tens_KL$ is a finite $L_1\tens_KL_1$-algebra
  and for all $n$ there is a prime $\fq_n\subset L\tens_KL$ above
  $\fp_n$. Then $\fq_1\subset\fq_2\subset\fq_3\cdots$ is an infinite
  strictly ascending chain of prime ideals of $L\tens_KL$, a
  contradiction.\nodemo
\end{enumerate}

\begin{cor}\label{cor:univ-noeth-field-extension}
  Let $L/K$ be an extension of fields. Then $L$ is a universally
  noetherian $K$-algebra if and only if $L$ is a finitely generated
  extension of $K$.
\end{cor}
\begin{demo}
  If $L/K$ is finitely generated, $L$ is a localization of a
  $K$-algebra of finite type, and the conclusion follows from
  corollaries \ref{cor:fft-is-noeth} and \ref{prop:noeth-sorites}
  (iii). Conversely if $L/K$ is universally noetherian, $L\tens_KL$ is
  a noetherian ring and the proposition shows that $L/K$ is finitely
  generated. 
\end{demo}

\begin{cor}\label{cor:field-extension-for-univ-noetherian}
  Suppose $f:\cX\to\cS$ is universally noetherian, $x\in\cX$ and
  $y=f(x)$. The field extension $\kappa(x)/\kappa(y)$ is finitely
  generated.
\end{cor}
\begin{demo}
  The morphism $\cX\times_\cS\Spf{\kappa(y)}\to\Spf{\kappa(y)}$ is
  also universally noetherian, so we may assume
  $\cS=\Spf{\kappa(y)}$. The morphism $\Spf{\hat\O_x}\to\cX$ is
  universally noetherian, where $\hat\O_x$ is the completion of the
  local ring at $x$ for the adic topology of $\O_x$. The closed
  immersion $\Sp{\kappa(x)}\to\Spf{\hat\O_\cX}$ is also universally
  noetherian. It follows that $\kappa(y)\to\kappa(x)$ is universally
  noetherian, and the assertion follows from the theorem.
\end{demo}

\begin{example}
  Suppose $k$ is a field. If $k[[x]]$ is given the $x$-adic topology,
  $k[x]\to k[[x]]$ is universally noetherian by corollary
  \ref{cor:completion-is-noeth}. On the other hand if $k[[x]]$ is
  given the discrete topology $k[x]\to k[[x]]$ is not universally
  noetherian; if it were $k(x)\to k((x))$ would be universally
  noetherian as well by corollary \ref{cor:localization-is
    noeth}. Since $k((x))/k(x)$ is not finitely generated this
  contradicts corollary \ref{cor:univ-noeth-field-extension}.

  There are more direct arguments for this example. Let $A$ be the
  integral closure of $k[x]$ in $k[[x]]$, i.e. the henselisation of
  $k[x]$ at $(x)$ and $L$ the fraction field of $A$. As in the proof
  of the proposition $\Sp{L\tens_{k(x)}L}$ is not a noetherian space,
  and since $\Sp{L\tens_{k(x)}L}$ is open in $\Sp{A\tens_{k[x]}A}$,
  $\Sp{A\tens_{k[x]}A}$ is not a noetherian space either. Finally
  $k[[x]]\tens_{k[x]}k[[x]]$ is a faithfully flat
  $A\tens_{k[x]}A$-algebra, so $\Sp{k[[x]]\tens_{k[x]}k[[x]]}$ is not
  a noetherian space.
\end{example}

\section{Differentials and Smoothness}
\label{sec:smooth-morphisms}

From now on all formal schemes are assumed to be locally
noetherian. For emphasis we will frequently restate this assumption
anyway.

\subsection{Differential invariants.}
\label{sec:differential-invariants}

As usual we start with the affine case, and then globalize. For the
whole of this section we fix a topological ring $R$ and a topological
$R$-algebra $A$, preadic with ideal of definition $J\subset A$.

\subsubsection{}
\label{sec:topologies}

We begin with a review of the topological aspects of the module of
relative 1-forms.  We denote by $I$ the diagonal ideal
$I=\Ker(A\tens_RA\to A)$, so that the ring of principal parts of order
$r$ and the module of relative 1-forms are
\begin{displaymath}
  P^r_{A/R}=(A\tens_RA)/I^{r+1},\qquad\diff_{A/R}=I/I^2.
\end{displaymath}
We denote by $d_0$, $d_1:A\to P^n_{A/R}$ the morphisms
$d_0(b)=b\tens1$, resp. $d_1(b)=1\tens b$.

Ideals of $A\ctens_RA$ will always have the induced topology.  We
topologize $\diff_{A/R}=I/I^2$ as a subquotient of $A\tens_RA$
(i.e. as a quotient of $I$ in the induced topology, or as a subobject
of $P^2_{A/R}$; these are the same). For $K=A\tens J+J\tens A$ this
coincides with the $K$-adic topology; this is evident if $A\tens_RA$
is noetherian (Artin-Rees), but in general it follows from the fact
that for any ideal $M\sset A$
\begin{displaymath}
  I\cap(A\tens M^2+M^2\tens A)\sset MI+I^2
\end{displaymath}
by \cite[$0_{IV}$ 20.4.5.1]{EGA}

(c.f.\ \cite[$0_{IV}$ Prop. 20.4.5]{EGA}). For the $A$-module structure
of $\diff_{A/R}$ defined by $d_0$ or $d_1$, we have
\begin{displaymath}
  J^n\diff_{A/R}=K^n\diff_{A/R}
\end{displaymath}
and the topology of $\diff_{A/R}$ is also the $J$-adic topology by
\cite[$0_{IV}$ Prop. 20.4.5]{EGA} (this is without any assumption that
$\diff_{A/R}$ is finitely generated). As in \cite[$0_{IV}$
\S20.7]{EGA} we denote by $\cdiff_{A/R}$ the completion of
$\diff_{A/R}$ with respect to its subquotient topology.  By the
previous remarks $\cdiff_{A/R}$ is also the $J$-adic completion of
$\diff_{A/R}$ when the latter is regarded as a $A$-module via $d_0$ or
$d_1$. We denote by $\hP^n_{A/R}$ the completion of $P^n_{A/R}$ with
respect to its topology as a quotient of $A\tens_RA$, i.e. the
$K$-adic topology.

If $R\to A$ and $A\to B$ are continuous homomorphisms of preadic rings,
the canonical exact sequence of relative 1-forms for the triple
$R\to A\to B$ induces a \textit{sequence}
\begin{equation}
  \label{eq:std-exact-sequence-completed-1}
  B\ctens_A\cdiff_{A/R}\to\cdiff_{B/R}\to\cdiff_{B/A}\to 0
\end{equation}
which is not necessarily exact. It is ``nearly exact'' in the sense
that $\cdiff_{B/R}\to\cdiff_{B/A}$ is surjective and the image of
$B\ctens_A\cdiff_{A/R}\to\cdiff_{B/R}$ is dense in the kernel of
$\cdiff_{B/R}\to\cdiff_{B/A}$, c.f. \cite[$0_{IV}$ 20.7.17.3]{EGA} and
the discussion there. If $A\to B$ is surjective with kernel $K$,
$\cdiff_{B/A}=0$ and there is similar sequence
\begin{equation}
  \label{eq:std-exact-sequence-completed-2}
  K/K^2\to B\ctens_A\cdiff_{A/R}\to\cdiff_{B/R}\to 0
\end{equation}
with the same ``near exactness'' property of
\ref{eq:std-exact-sequence-completed-1}, c.f. \cite[$0_{IV}$ 20.7.20]{EGA}

The sequence
\begin{displaymath}
  0\to I\to A\tens_RA\to A\to0
\end{displaymath}
is strict exact; in fact by construction $I$ has the induced topology,
and the image of $K^n\sset A\tens_RA$ in $A$ is $J^n$. Since the
completion of a strict exact sequence is strict exact \cite[Ch. III
\S2 no. 12 Lemme 2]{bourbaki-AC}, the sequence
\begin{equation}
  \label{eq:hatI-as-kernel}
  0\to\hI\to A\ctens_RA\to A\to 0
\end{equation}
is exact, in which $A\ctens_RA\to A$ is induced by
$a\ctens b\mapsto ab$. Our first goal is to show that when $A$ is a
universally noetherian $R$-algebra we may identify
$\cdiff_{A/R}\simeq\hI/\hI^2$ and
$\hat P^n_{A/R}\simeq A\ctens_RA/\hI^{r+1}$, c.f. proposition
\ref{prop:Omega-hat} below. 

\begin{prop}\label{prop:Omega-hat}
  With the hypotheses and notation of
  \S\ref{sec:differential-invariants}, suppose that $R$ is noetherian,
  and denote by $I\subset A\tens_RA$ the diagonal ideal of $R\to
  A$. If $A$ is a universally noetherian $R$-algebra, there are
  functorial isomorphisms
  \begin{displaymath}
    \cdiff_{A/R}\simeq\hI/\hI^2
    \quad\text{and}\quad
    \hP^n_{A/R}\simeq(A\ctens_RA)/\hI^{r+1}.
  \end{displaymath}
\end{prop}
\begin{demo}
  If $J\subset A$ is an ideal of definition, $J$ is finitely generated
  and thus $K=A\tens J+J\tens A$ is a finitely generated ideal of
  $R=A\tens_RA$. The sequence
  \begin{displaymath}
    0\to I^2\to I\to\diff_{A/R}\to 0
  \end{displaymath}
  is strict exact by definition of the topologies involved, so its
  completion
  \begin{displaymath}
    0\to (I^2)\,\widehat\relax\to\hI\to\cdiff_{A/R}\to 0
  \end{displaymath}
  is also strict exact. Since $A$ is a universally noetherian
  $R$-algebra, $A\ctens_RA$ is noetherian and corollary
  \ref{cor:completions-of-powers} shows that this exact sequence is
  \begin{displaymath}
    0\to\hI^2\to\hI\to\cdiff_{A/R}\to 0
  \end{displaymath}
  and the first assertion follows. The second is proven in the same way.
\end{demo}

\begin{cor}\label{cor:Omega-hat-finitely-generated}
  If $R$ is noetherian and $A$ is an universally noetherian
  $R$-algebra, $\cdiff_{A/R}$ is generated as a $A$-module by finitely
  many elements of the form $1\ctens x-x\ctens1$.
\end{cor}
\begin{demo}
  We know $I$ is generated by elements of the form $1\tens x-x\tens1$
  and $\hI$ is the $K$-adic completion of $I$, where as before
  $K=A\tens_RJ+J\tens_RA$. Then $\hI$ is generated by the
  $1\ctens x-x\ctens1$ for all $x\in A$, and thus by finitely many of
  them, since $A\ctens_RA$ is noetherian.
\end{demo}

For $x\in A$ we will use $\md x$ to denote both the image of
$1\tens x-x\tens1$ in $\diff_{A/R}$ and the image of
$1\ctens x-x\ctens1$ in $\cdiff_{A/R}$; this should not cause
confusion.

\begin{cor}\label{cor:standard-exact-sequences}
  If $R\to A\to B$ are homomorphisms of adic noetherian rings with
  $R\to A$ and $R\to B$ universally noetherian, the sequence
  \begin{displaymath}
    B\tens_A\cdiff_{A/R}\to\cdiff_{B/R}\to\cdiff_{B/A}\to0
  \end{displaymath}
  is exact. If $A\to B$ is surjective with kernel $K$, the sequence
  \begin{displaymath}
    K/K^2\to B\tens_A\cdiff_{A/R}\to\cdiff_{B/R}\to 0
  \end{displaymath}
  is exact.
\end{cor}
\begin{demo}
  By proposition \ref{prop:properties-of-loc-noeth2}, the hypotheses
  imply that $A\to B$ is universally noetherian. Therefore
  $\cdiff_{A/R}$ is a finitely generated $A$-module and $\cdiff_{B/R}$
  and $\cdiff_{B/A}$ are finitely generated $B$-modules. Thus inA
  \ref{eq:std-exact-sequence-completed-1} we may replace the completed
  tensor product by an ordinary one.  ince $B$ is noetherian, any
  submodule of the finitely generated modules $\cdiff_{B/R}$ and
  $\cdiff_{B/A}$ is closed, and exactness follows from the ``near
  exactness'' of the sequence
  \ref{eq:std-exact-sequence-completed-1}. The argument in the case of
  the second sequence is the same.
\end{demo}

Let $J$ be an ideal of definition $J\subset A$ and set
$A_n=A/J^{n+1}$. For $n'\ge n$ there is a natural $A$-module
homomorphism $\diff_{A_{n'}/R}\to\diff_{A_n/R}$, and the discussion of 
\cite[$O_{IV}$ 20.7.14]{EGA} shows that their inverse limit is the
separated completion of $\diff_{A/R}$, whence a canonical isomorphism 
\begin{equation}
  \label{eq:alternate-cdiff}
  \cdiff_{A/R}\simeq\liminv\diff_{A_n/R}
\end{equation}

\begin{example}
  Let $J$ be an ideal of definition of $R$, and let
  $A=R\{T_1,\ldots,T_d\}$ be the $J$-adic completion of the polynomial
  ring $R[T_1,\ldots,T_d]$. With $JA$ is an ideal of definition of
  $A$, $A$ is an $R$-algebra that is topologically of finite type, and
  therefore universally noetherian. Then \ref{eq:alternate-cdiff}
  shows that $\cdiff_{A/R}$ is free over $A$ with basis
  $\md T_1,\ldots,\md T_d$.
\end{example}

\begin{prop}\label{prop:cdiff-of-completions}
  Let $A$ be a universally noetherian $R$-algebra. If $M\subset A$ is
  an ideal containing an ideal of definition and $B$ is the $M$-adic
  completion of $A$, there is a natural and functorial isomorphism
  \begin{displaymath}
    B\tens_A\cdiff_{A/R}\isom\cdiff_{B/R}
  \end{displaymath}
\end{prop}
\begin{demo}
  Setting $K=M^n$ in the second exact sequence of corollary
  \ref{cor:standard-exact-sequences} and $B_n=A/M^n$ yields exact
  sequences
  \begin{displaymath}
    M^n/M^{2n}\to(A/M^n)\tens_A\cdiff_{A/R}\to\diff_{B_n/R}\to0
  \end{displaymath}
  for all $n\ge0$. Since the pro-object $\{M^n/M^{2n}\}_{n\ge0}$ is
  essentially zero, its image in
  $\{(A/M^n)\tens_A\cdiff_{A/R}\}_{n\ge0}$ is essentially zero and in
  particular Mittag-Leffler. Therefore the inverse limit over $n$ is
  an isomorphism
  \begin{displaymath}
    B\ctens_A\cdiff_{A/R}\isom\liminv_n\diff_{B_n/R}
  \end{displaymath}
  and we may replace the completed tensor product by an ordinary one
  since $\cdiff_{A/R}$ is finitely generated. The assertion then
  follows from \ref{eq:alternate-cdiff}.
\end{demo}

\begin{example}
  If $(R,J)$ is adic noetherian, we saw in the last example that for
  $R$-algebra $A=R\{X_1,\ldots,X_d\}$ with the $J$-adic topology,
  $\cdiff_{A/R}$ is free with basis $\md X_1,\ldots,\md X_d$. Then
  $B=R[[X_1,\ldots,X_d]]$ is completion of $A$ with respect to the
  ideal $M=J+(X_1,\ldots,X_d)$, and proposition
  \ref{prop:cdiff-of-completions} says that $\cdiff_{B/R}$ is the free
  $B$-module on $\md X_1,\ldots,\md X_d$; compare this with
  \cite[$0_{IV}$ Cor. 21.9.3]{EGA}.
\end{example}

\subsubsection{Deformations.}
\label{sec:deformations}

One more consequence of \cite[$0_{IV}$ 20.7.14]{EGA} will be
useful. Suppose $A$ is a topological $R$-algebra and $B$ is a discrete
topological $R$-algebra with an ideal $I\subset B$ such that
$I^2=0$. If a continuous $R$-homomorphism $u_0:A\to B$ is given, the
set of continuous $R$-homomorphisms $u:A\to B$ having the same
composite with $B\to B/I$ is principal homogenous under the $A$-module
of continuous derivations $A\to I$; the argument is the same as the
discrete case \cite[$0_{IV}$ Prop. 20.1.1]{EGA}. Since $u_0$ is
continuous and $B$ is discrete, $I$ is annihilated by an open ideal of
$A$, and it follows that the $A$-module of continuous $R$-derivations
$A\to I$ is the same as the set of continuous homomorphisms
$\cdiff_{A/R}\to I$ of $\hat A$-modules, c.f. \cite[$0_{IV}$
20.7.14.4]{EGA}.  When $A$ is adic, the topology of $\cdiff_{A/R}$ is
induced by the topology of $A$, and it follows that any $A$-linear
$\cdiff_{A/R}\to I$ is continuous. Therefore the set of $u:A\to B$
having the same composition with $B\to B/I$ as $u_0$ is principal
homogenous under the group $M=\Hom_A(\cdiff_{A/R},I)$. Since $I^2=0$,
the $A$-module structure of $I$ comes from a $B/I$-module structure,
and
\begin{equation}
  \label{eq:deformation-group}
  M\simeq\Hom_{B/I}((B/I)\tens_A\cdiff_{A/R},I).   
\end{equation}

\subsubsection{Globalization.}
\label{sec:1-forms-global}

If $A$ is a universally noetherian $R$-algebra and $S\subset A$ is a
multiplicative system, $\cdiff_{A/R}$ is a finitely generated
$A$-module, and therefore the completion of the canonical isomorphism
$\diff_{S^{-1}A/R}\simeq S^{-1}\diff_{A/R}$ is an isomorphism
\begin{displaymath}
  \cdiff_{A\{S^{-1}\}/R}\simeq A\{S^{-1}\}\tens_A\cdiff_{A/R}.
\end{displaymath}
of finitely generated $A\{S^{-1}\}$-modules, where as before
$A\{S^{-1}\}$ is the completion of $S^{-1}A$. From this we see that
the construction of $\cdiff_{A/R}$ globalizes: for any universally
noetherian morphism $f:\cX\to\cS$ there is a coherent $\O_\cX$-module
$\diff_{\cX/\cS}$ such that
\begin{equation}
  \label{eq:globalization-of-cdiff}
  \Gamma(\Spf{A},\diff_{\cX/\cS})=\cdiff_{A/R}
\end{equation}
for all open affines $\Spf{A}\sset\cX$, $\Spf{R}\sset\cS$ such that
$f(\Spf{A})\sset\Spf{R}$ and $A$ is a universally noetherian
$R$-algebra (note the absence of a ``hat,'' as per our convention on
completions in the introduction). A similar argument shows that for
any $r\ge0$ there is a sheaf $\cP^r_{\cX/\cS}$ of rings such that
\begin{equation}
  \label{eq:globalization-of-principal-parts}
  \Gamma(\Spf{A},\cP^r_{\cX/\cS})=\hat P^r_{A/R}.
\end{equation}
Like its affine counterpart $\cP^r_{\cX/\cS}$ has two $\O_\cX$-algebra
structures, with respect to both of which it is a finite
$\O_\cX$-algebra. 

In the previous discussion we did not need to assume that $\cX\to\cS$
was separated, but if it is the diagonal $\cX\to\cX\times_\cS\cX$ is a
closed immersion defined by an ideal $\cI\subset\O_{\cX\times_\cS\cX}$
with the property that in the above affine setting,
\begin{displaymath}
  \Gamma(\Spf{A},\cI)=\hat I
\end{displaymath}
where as before $I\subset A\ctens_RA$ is the kernel of
$A\ctens_RA\to A$. In this case proposition \ref{prop:Omega-hat}
yields canonical and functorial isomorphisms
\begin{equation}
  \label{eq:global-diff-and-principal-parts}
  \diff_{\cX/\cS}=\cI/\cI^2,\qquad
  \cP^n_{\cX/\cS}=\O_{\cX\times_\cS\cX}/\cI^n.
\end{equation}

The standard exact sequences globalize immediately: if $\cY\to\cS$ and
$\cX\to\cS$ are universally noetherian and $f:\cY\to\cX$ is a morphism
(necessarily universally noetherian), the sequence
\begin{equation}
  \label{eq:std-exact-sequence-completed-3}
  f^*\diff_{\cX/\cS}\to\diff_{\cY/\cS}\to\diff_{\cY/\cX}\to0
\end{equation}
is exact; if in addition $\cY\to\cX$ is a closed immersion with ideal
$\cK$, the sequence
\begin{equation}
  \label{eq:std-exact-sequence-completed-4}
  \cK/\cK^2\to f^*\diff_{\cX/\cS}\to\diff_{\cY/\cS}\to0
\end{equation}
is exact; these assertions follow immediately for corollary
\ref{cor:standard-exact-sequences}. Finally, if $J\subset\O_\cX$ is an
ideal of definition and $X_n$ is the closed subscheme of $\cX$ defined
by $J^{n+1}$, the isomorphism \ref{eq:alternate-cdiff} 
globalizes to
\begin{equation}
  \label{eq:alternate-diff-global}
  \diff_{\cX/\cS}\simeq\liminv_n\diff_{X_n/\cS}.
\end{equation}
From proposition \ref{prop:diff-of-completions} we get

\begin{prop}\label{prop:diff-of-completions}
  If $\cX\to\cS$ is separated and universally noetherian and
  $Y\subset\cX$ is a closed subscheme, the canonical
  morphism
  \begin{displaymath}
    i_Y^*\diff_{\cX/\cS}\to\diff_{\hat\cX_Y/\cS}
  \end{displaymath}
  is an isomorphism.\nodemo
\end{prop}

The deformation theory of section \ref{sec:deformations} globalizes in
the same way. If $f:\cX\to\cS$ is universally noetherian, $Z$ is an
affine scheme over $\cS$, $Z_0\inj Z$ is a closed immersion whose
ideal $I$ is such that $I^2=0$, and $g:Z_0\to\cX$ is a $\cS$-morphism,
the sheaf of liftings of $g$ to a morphism $u:Z\to\cX$ making the
diagram
\begin{equation}
  \label{eq:formal-smooth-diagram2}
  \xymatrix{
    Z_0\ar[r]^g\ar[d]&\cX\ar[d]\\
    Z\ar[r]\ar[ur]^u&\cS
  }  
\end{equation}
commutative is a pseudo-torsor under the sheaf
$M=Hom_{\O_{Z_0}}(g^*\diff_{\cX/\cS},I)$.


\subsection{Smooth and quasi-smooth morphisms.}
\label{sec:smooth-formal-case}

The standard definition of smoothness extends immediately to the
setting of formal schemes:

\begin{defn}\label{defn:smooth}
  A morphism $f:\cX\to\cS$ of locally noetherian formal schemes is
  \emph{smooth} (resp. \emph{unramified}, \emph{\'etale}) if it is of
  finite type and formally smooth (resp. formall unramified, formally
  \'etale) in the sense of definition \ref{defn:formally-smooth}.
\end{defn}
This definition is equivalent to the one given by Berthelot
\cite[2.1.5]{berthelot:1996} in the case of morphisms of locally
noethrian formal schemes over a complete discrete valuation ring.

\begin{defn}\label{defn:quasi-smooth}
  A morphism $f:\cX\to\cS$ of locally noetherian formal
  schemes is \emph{quasi-smooth} (resp. \emph{quasi-unramified},
  \emph{quasi-\'etale}) if it is universally noetherian and formally
  smooth (resp. formally unramified, formally \'etale)
\end{defn}

This definition of ``quasi-smooth'' conflicts with \cite[Ch. IV
1.5.1]{berthelot:1974}. As we will not use Berthelot's notion, this
will not be a problem.

It is clear that a smooth morphism is quasi-smooth, and conversely a
quasi-smooth morphism of finite type is smooth. By proposition
\ref{prop:formally-smooth-implies-flat}, a quasi-smooth morphism is
flat. In particular a quasi-\'etale morphism is flat and
quasi-unramified; I do not know if the converse is true.

The next proposition summarizes the basic properties of quasi-smooth,
quasi-unramified and quasi-\'etale morphisms. They follow from the
results on universally noetherian morphisms in section
\ref{sec:formal-schemes-finiteness} and basic properties of formally
smooth (resp. formally unramified, formally \'etale) morphisms whose
proofs are entirely parallel to the corresponding assertions for
morphisms of schemes, c.f. \cite[\S17]{EGA} propositions 17.1.3--5.

\begin{prop}\label{prop:quasi-smooth-sorites}
  \begin{enumerate}
  \item An immersion is quasi-unramified. An open immersion is
    quasi-\'etale.
  \item If $f:\cX\to\cS$ and $g:\cY\to\cX$ are quasi-smooth
    (resp. quasi-unramified, quasi-\'etale) then so is
    $g\circ f:\cY\to\cS$.
  \item If $\cX\to\cS$ is quasi-smooth (resp. quasi-unramified,
    quasi-\'etale) and $\cS'\to\cS$ is any morphism of locally
    noetherian schemes, $\cX\times_\cS\cS'\to\cS'$ is quasi-smooth
    (resp. quasi-unramified, quasi-\'etale).
  \item If $f:\cX\to\cS$ and $g:\cY\to\cS$ are quasi-smooth
    (resp. quasi-unramified, quasi-\'etale) then so is
    $f\times g:\cX\times_\cS\cY\to\cS$.
  \item If $f:\cX\to\cS$ and $g:\cY\to\cX$ are morphisms such that
    $f\circ g$ is quasi-unramified, then $g$ is quasi-unramified.
  \item If $f:\cX\to\cS$ is quasi-unramified, $g:\cY\to\cX$ is a
    morphism and $f\circ g$ is quasi-smooth (resp. quasi-\'etale) then
    $g$ is quasi-smooth (resp. quasi-\'etale).
  \item If $f:\cX\to\cS$ is quasi-\'etale and $g:\cY\to\cX$ is a
    morphism then $f\circ g$ is quasi-smooth (resp. quasi-\'etale) if
    and only if $g$ is quasi-smooth (resp. quasi-\'etale).\nodemo
  \end{enumerate}
\end{prop}
In statements (v) through (vii) note that $f\circ g$ being
quasi-smooth or quasi-\'etale implies that it is universally
noetherian, so that $g$ is universally noetherian by proposition
\ref{prop:properties-of-loc-noeth2}. 

We can now return to a question that was left open in section
\ref{sec:formal-flatness}:

\begin{prop}\label{prop:quasi-smooth-is-local}
  Let $f:\cX\to\cS$ be a universally noetherian morphism.
  \begin{enumerate}
  \item If $\{U_\alpha\}$ is an open cover of $\cX$ and
    $f_\alpha:U_\alpha\to\cS$ is the composite of $f$ with the open
    immersion $U_\alpha\to\cX$, then $f$ is quasi-smooth
    (resp. quasi-unramified, quasi-\'etale) if and only if all the
    $f_\alpha$ are quasi-smooth (resp. quasi-unramified,
    quasi-\'etale).
  \item If $\{V_\alpha\}$ is an open cover of $\cS$ then $f$ is
    quasi-smooth (resp. quasi-unramified, quasi-\'etale) if and only
    if the morphisms $f^{-1}(V_\alpha)\to V_\alpha$ are quasi-smooth
    (resp. quasi-unramified, quasi-\'etale).
  \end{enumerate}
\end{prop}
\begin{demo}
  Assertion (ii) follows from (i), and the quasi-unramified and
  quasi-\'etale cases of (i) are formal consequences of the
  definitions, as in the proof of \cite[IV Prop. 17.1.6]{EGA},
  together with the basic properties of universally noetherian
  morphisms. In the quasi-smooth case the main thing is to prove
  formal smoothness, and we can again follow the argument of
  \textit{loc.\ cit.}, the point being that with the assumptions of
  (i), the set of local liftings $u$ in the diagram
  \ref{eq:formal-smooth-diagram2} is a torsor under the sheaf
  $M=Hom_{\O_{Z_0}}(h^*\diff_{\cX/\cS},\cI)$, and $Z_0$ being an
  affine scheme, this torsor is trivial.
\end{demo}

\begin{prop}\label{quasi-unramified-criterion}
  A universally noetherian morphism $\cX\to\cS$ is quasi-unramified
  if and only if $\diff_{\cX/\cS}=0$.
\end{prop}
\begin{demo}
  We may assume that $\cX=\Spf{A}$ and $\cS=\Spf{R}$ are formally
  affine, so that $A$ is a universally noetherian $R$-algebra. It
  suffices to check that the $R$-algebra $A$ is formally unramified if
  and only if $\cdiff_{A/R}=0$, which is \cite[$0_{IV}$
  Prop. 20.7.4]{EGA}.
\end{demo}

\begin{cor}\label{inclusion-of-completion-is-quasi-unramified}
  If $Y\subset\cX$ is a closed subscheme then $i_Y:\hat\cX_\cY\to\cX$
  is quasi-\'etale.
\end{cor}
\begin{demo}
  We know that $i_Y$ is universally noetherian, so by proposition
  \ref{quasi-unramified-criterion} it suffices to show that
  $\diff_{\hat\cX_Y/\cX}=0$. We can assume $\cX=\Spf{A}$ is formally
  affine, in which case $\hat\cX_Y=\Spf{\hat A}$, and then
  $\cdiff_{\hat A/A}\simeq\cdiff_{\hat A/\hat A}=0$.
\end{demo}

For example, if $\cY\to\cS$ is quasi-smooth and $\cX$ is the
completion of $\cY$ along a closed subscheme then $\cX\to\cS$ is
quasi-smooth. Note that the corresponding morphism of reduced closed
subschemes need not be smooth. The next proposition shows that a
standard criterion for smoothness in the case of a morphism of finite
type remains true in the general case; as a consequence we get a
structure theorem for quasi-smooth morphisms analogous to the usual
one for morphisms of finite type.

\begin{lemma}\label{lemma:local-direct-factor}
  Suppose $B$ is adic and noetherian and $f:M\to N$ is a homomorphism
  of finitely generated $B$-modules.
  \begin{enumerate}
  \item $f$ is the inclusion of a direct summand in the category of
    $B$-modules if and only if $\Hom_B(N,L)\to\Hom_B(M,L)$ is
    surjective for all discrete $B$-modules $L$ annihilated by an open
    ideal of $B$.
  \item If $N$ is projective and $\fp\in\Spf{B}$ is such that the
    $B_\fp$-module homomorphism $M_\fp\to N_\fp$ is the inclusion of a
    direct summand, there is a $f\in B\setminus\fp$ such that the
    $B_f$-module homomorphism $M_f\to N_f$ is the inclusion of a
    direct summand.
  \end{enumerate}
\end{lemma}
\begin{demo}
  For (i), the condition is clearly necessary, and to show that it is
  sufficient it suffices to show that $\Hom_B(N,L)\to\Hom_B(M,L)$ is
  surjective for all finitely generated $B$-modules $L$. If $L$ is
  finitely generated and $J\subset B$ is an ideal of definition, the
  hypothesis implies that $\Hom_B(N,L/J^nL)\to\Hom_B(M,L/J^n)$ is
  surjective; since $M$ and $N$ are finitely generated this says that
  \begin{displaymath}
    \Hom_B(N,L)\tens_BB/J^n\to\Hom_B(M,L)\tens_BB/J^n
  \end{displaymath}
  is surjective for all $n$, and the assertion follows by the faithful
  flatness of the $J$-adic completion. Part (ii) is a consequence of
  \cite[Cor. 19.1.12]{EGA}.
\end{demo}

The proof of the next proposition uses properties of the functor
$\Exalcotop$ defined in \cite[$0_{IV}$ \S18.5]{EGA}. One property of
this functor apparently not stated in \textit{loc.\ cit.} is the
isomorphism
\begin{equation}
  \label{eq:Exalcotop-surjection}
  \Exalcotop_A(B,L)\simeq\Hom_B(K/K^2,L)
\end{equation}
which holds when $A\to B$ is a surjective homomorphism of linearly
topologized rings with kernel $K$. It is however an immediate
consequence of a similar isomorphism \cite[$0_{IV}$ \S18.5.4.1]{EGA}
for the related functor $\mathrm{Exantop}$; the case of discrete rings
is \cite[$0_{IV}$ \S18.3.8.1]{EGA} and \cite[$0_{IV}$
\S18.4.2.1]{EGA}.

\begin{prop}\label{prop:smoothness-criteria}
  Suppose that $\cY/\cS$ and $\cX/\cS$ are universally noetherian and
  $f:\cY\to\cX$ is an $\cS$-morphism.
  \begin{enumerate}
  \item If $\cY\to\cS$ is quasi-smooth, then $f$ is quasi-smooth
    if and only if the $\O_\cY$-module homomorphism
    $f^*\diff_{\cX/\cS}\to\diff_{\cY/\cS}$ is the inclusion of a local
    direct summand.
  \item (Jacobian criterion) If $f$ is a closed immersion with ideal
    $\cK$ and $\cX/\cS$ is quasi-smooth, then $\cY\to\cS$ is
    quasi-smooth if and only if the $\O_\cY$-module homomorphism
    $\cK/\cK^2\to f^*\diff_{\cX/\cS}$ is the inclusion of a local
    direct summand.
  \end{enumerate}
\end{prop}
\begin{demo}
  The hypotheses imply that $\cY\to\cX$ is universally noetherian
  (prop. \ref{prop:properties-of-loc-noeth}, (iii)), so we only have
  to prove formal smoothness in both cases. We may work locally
  everywhere, so we may assume that $\cY\to\cX\to\cS$ is the formal
  spectrum of $R\to A\to B$. By \cite[$0_{IV}$ 19.4.4]{EGA}, a
  morphism $A\to B$ of topological algebras is formally smooth if and
  only if $\Exalcotop_A(B,L)=0$ for any discrete $B$-module $L$
  annihilated by an open ideal of $B$.  In the situation of the
  proposition there is an exact sequence
  \begin{multline*}
    0\to\Der_A(B,L)\to\Der_R(B,L)\to\Der_R(A,L)\to\\
    \to\Exalcotop_A(B,L)\to\Exalcotop_R(B,L)\to\Exalcotop_R(A,L)
  \end{multline*}
  for any discrete $B$-module $L$ annihilated by an open ideal of $B$
  (c.f. \cite[$0_{IV}$]{EGA} Propositions 20.3.5 and 20.3.7).

  In case (i) we have $\Exalcotop_R(B,L)=0$ for all $L$ as above, and
  thus $\Exalcotop_A(B,L)=0$ if and only if
  $\Der_R(B,L)\to\Der_R(A,L)$ is surjective. Equivalently, $A\to B$ is
  formally smooth if and only if
  $\Hom_B(\cdiff_{B/R},L)\to\Hom_B(B\tens_A\cdiff_{A/R},L)$ is
  surjective for all discrete $L$ annihilated by an open ideal of
  $B$. By (i) of the lemma, this is equivalent to
  $B\tens_A\cdiff_{A/R}\to\cdiff_{B/R}$ being the inclusion of a
  direct summand.

  In the case of (ii) we have $\Exalcotop_R(A,L)=0$ and
  $\Exalcotop_A(B,L)$ is given by \ref{eq:Exalcotop-surjection}).
  Therefore $R\to B$ is formally smooth if and only if
  $\Hom_B(B\tens_A;\cdiff_{A/R},L)\to\Hom_B(K/K^2,L)$ is surjective
  for all $L$. As before, this condition is equivalent to
  $K/K^2\to B\tens_A\cdiff_{A/R}$ being the inclusion of a direct
  summand.
\end{demo}

\begin{cor}
  Let $\cX\to\cS$ be a morphism of formal schemes, $x\in\cX$ is a
  point and $\hat\O_x$ is the completion of the local ring $\O_x$ with
  respect to the adic topology of $\O_\cX$. Suppose that some
  neighborhood $\cV$ of $x$ has a closed embedding $\cV\inj\cY$ over
  $\cS$ into a quasi-smooth $\cS$-scheme $\cY$. If the composite
  morphism $\Spf{\hat\O_x}\to\cS$ is quasi-smooth, there is an open
  neighborhood $\cU$ of $x$ such that $\cU\to\cS$ is quasi-smooth.
\end{cor}
\begin{demo}
  We can assume $\cX=\cV$, and denote by $\cK$ the ideal of
  $f:\cX\to\cY$. Since $\cY\to\cS$ is universally noetherian, so is
  $\cX\to\cS$. By (ii) of lemma \ref{lemma:local-direct-factor} there
  is an open neighborhood $\cU\sset\cX$ of $x$ on which the
  $\O_\cX$-module homomorphism $\cK/\cK^2\to f^*\diff_{\cY/\cS}$ is
  the inclusion of a (local) direct summand, and the assertion follows
  from (ii) of proposition \ref{prop:smoothness-criteria}.
\end{demo}

The condition that $\cX\to\cS$ is locally embeddable into a
quasi-smooth $\cY\to\cS$ is satisfied when $\cX\to\cS$ is of finite
type, so this condition could be regarded as a weak finiteness
property.

The next proposition is a very special case of a very general (and
difficult) criterion of Grothendieck for a homomorphism of topological
rings to be formally smooth, c.f.  \textit{loc. cit.} Th. 19.5.3 and
Cor. 19.5.7, and more particularly \cite[$0_{IV}$ Rem. 19.5.8]{EGA}.
Its conclusion could be taken as a definition of
\textit{differentiably smooth} in the formal case:

\begin{prop}\label{prop:Grothendieck-formal-smoothness}
  Let $A\to B\to C$ be continuous homomorphisms of topological rings,
  and suppose that $A\to B$ formally smooth, $B\to C$ surjective with
  finitely generated kernel $I$, and $C$ is a Zariski ring. Then
  $I/I^2$ is a projective $C$-module and for all $n\ge0$ the natural
  map $\Sym_C^n(I/I^2)\to I^n/I^{n+1}$ is an isomorphism.\nodemo
\end{prop}

\begin{prop}\label{prop:formally-smooth-and-diff-smooth-formal-case}
  Suppose $f:\cX\to\cS$ is a quasi-smooth morphism of locally noetherian
  formal schemes, and let $\cI$ be the ideal of the diagonal of
  $f$. Then
  \begin{enumerate}
  \item $\diff_{\cX/\cS}$ is a locally free $\O_\cX$-module of finite
    type, and
  \item the natural morphism
    \begin{displaymath}
      \Sym^n_{\O_\cX}(\diff_{\cX/\cS})\to\cI^n/\cI^{n+1}
    \end{displaymath}
    is an isomorphism for all $n\ge0$.
  \end{enumerate}
\end{prop}
\begin{demo}
  We may assume that $\cX\to\cS$ is the formal spectrum of $A\to B$,
  so that $A$ and $B$ are noetherian and $A\to B$ is universally
  noetherian. Then $B\ctens_AB$ is noetherian, and in fact a Zariski
  ring, so we may apply to previous proposition with $B$, $C$ and $I$
  replaced by $B\ctens_AB$, $B$ and $\hI=\Ker(B\ctens_AB\to B)$
  respectively. The conclusion follows since $\cI$ is the sheaf of
  ideals associated to $\hI$, and $\diff_{\cX/\cS}$ is the module
  associated to $\hI/\hI^2$.
\end{demo}

\begin{cor}\label{cor:smooth-implies-regular-diagonal}
  If $f:\cX\to\cS$ is a quasi-smooth morphism of locally noetherian
  formal schemes, the diagonal $\cX\to\cX_\cS(r)$ is a regular
  immersion.
\end{cor}
\begin{demo}
  As assertion is local we may assume $\cX=\Spf{A}$ and $\cS=\Spf{R}$
  are affine, and then $\cX_\cS(r)=\Spf{A_R(r)}$. Both $A$ and
  $A_R(r)$ are formally smooth $R$-algebras, and $A$ is a Zariski
  ring. The ideal $I\subset A_R(r)$ of the diagonal is the kernel of
  the multiplication map $A_R(r)\to A$; as this is surjective,
  proposition \ref{prop:Grothendieck-formal-smoothness} shows that $I$
  is quasi-regular, i.e. locally generated by a quasi-regular sequence
  $(f_i)$. Since $B$ is noetherian it is $I$-adically separated, and
  $(f_i)$ is regular by \cite[$0_{IV}$ 15.1.9]{EGA}.
\end{demo}

If $\cX/\cS$ is quasi-smooth and $\diff_{\cX/\cS}$ has constant rank, we
call this rank the \textit{formal relative dimension} of $\cX/\cS$,
and denote it by $\fdim{\cX/\cS}$. If $\cX/\cS$ is of finite type,
proposition \ref{prop:diff-of-completions} shows that $\fdim{\cX/\cS}$
is the relative dimension of $X/S$, but this is not true in general.

When $\cX/\cS$ is quasi-smooth we will say that an open affine $U\sset\cX$
is \textit{parallelizable} if there are
$x_1,\ldots,x_d\in\Gamma(U,\O_\cX)$ such that
$\{\md x_1,\ldots,\md x_d\}$ is a basis of
$\Gamma(U,\diff_{\cX/\cS})$; we will also say that the
$x_1,\ldots,x_d$ are \textit{local coordinates} on $U$ for
$\cX/\cS$. It is clear that when $\cX/\cS$ is quasi-smooth, the topology of
$\cX$ has a basis consisting of parallelizable open sets.

As an application we get a structure theorem for quasi-smooth morphisms
similar to the one that obtains for formally smooth morphisms of
finite type:

\begin{cor}\label{cor:structure-of-smooth-morphisms}
  Suppose $f:\cX\to\cS$ is quasi-smooth and $d=\fdim{\cX/\cS}$. Locally on
  $\cX$ there is a factorisation $f=p\circ g$ where
  $g:\cX\to\bA^d_\cS$ is quasi-\'etale and $p:\bA^d_\cS\to\cS$ is the
  canonical projection. In fact this factorisation exists on any open
  parallelizable $U\sset\cX$.
\end{cor}
\begin{demo}
  We may assume $\cX=\Spf{B}$ and $\cS=\Spf{R}$, so that $B$ a
  quasi-smooth $R$-algebra, and by further shrinking $\cX$ we can
  assume $\diff_{B/R}$ is free with basis $\md x_1,\ldots,\md
  x_d$. The sections $x_1,\ldots,x_d$ give a factorisation
  $R\to R[T_1,\ldots,T_d]\to B$ with $T_i\mapsto x_i$, and since
  $R\to B$ is continuous this extends to a factorisation
  $R\to R\{T_1,\ldots,T_d\}\to B$. If $A=R\{T_1,\ldots,T_d\}$ is given
  the topology induced by $A$, proposition
  \ref{prop:smoothness-criteria} and the example after equation
  \ref{eq:alternate-cdiff} show that $\cdiff_{B/A}=0$. Therefore
  $A\to B$ is quasi-unramified by proposition
  \ref{quasi-unramified-criterion}, and it suffices to show that
  $A\to B$ is formally smooth. By proposition
  \ref{prop:smoothness-criteria} this is so if and only if
  $B\tens_A\cdiff_{A/R}\to\cdiff_{B/R}$ is the inclusion of a local
  direct summand in the category of $B$-modules, but by construction
  this map is actually an isomorphism.
\end{demo}


\begin{prop}\label{prop:etale-and-radicial}
  A morphism of locally noetherian formal schemes that is \'etale and
  radicial is an open immersion.
\end{prop}
\begin{demo}
  Suppose $f:\cY\to\cX$ is \'etale and radicial, and let
  $J\subset\O_\cX$ be an ideal of definition. Since $f$ is \'etale it
  is adic, and $f^*J\subset\O_\cY$ is an ideal of definition. If we
  set $X_n=V(J^n)$ and $Y_n=V(J^n\O_\cY)$ then $X_n$, $Y_n$ are
  ordinary schemes and for all $n$ the induced map $f_n:Y_n\to X_n$ is
  \'etale and radicial. It is therefore an open immersion by \cite[IV
  Th. 17.9.1]{EGA}, and since $f$ is the inductive limit of the $f_n$,
  it is an open immersion as well.
\end{demo}

\begin{remark}
  We cannot weaken the hypothesis ``\'etale'' to ``quasi-\'etale.'' In
  fact if $\cX$ is any formal scheme and $\hat\cX_Y$ is the completion
  of $\cX$ along a proper closed subscheme, the canonical morphism
  $i_Y:\hat\cX_Y\to\cX$ is quasi-\'etale and radicial, but not an open
  immersion.
\end{remark}

\begin{prop}\label{prop:flatness-of-Frobenius}
  Let $f:\cX\to\cS$ be a quasi-smooth morphism of formal schemes of
  characteristic $p$.
  \begin{enumerate}
  \item The relative $q$th power Frobenius
    $F_{\cX/\cS}:\cX\to\niv{\cX}{q}$ is flat.
  \item If $f$ is formally of finite type, $F_{\cX/\cS}$ is finite.
  \end{enumerate}

\end{prop}
\begin{demo}
  We first prove (i) and (ii) in the case when $f$ is of finite
  type. Then $f$ is smooth, and is the inductive limit of a sequence
  of smooth morphisms $f_n:X_n\to S_n$ of schemes; here we choose an
  ideal of definition $J$ of $\cS$, and have set $S_n=V(J^n)$ and
  $X_n=V(f^*J^n)$. Then $F_{X_n/S_n}$ is finite and flat for all
  $n>0$, and for $n'\ge n$ $f_n$ is the reduction modulo $J^n$ of
  $f_{n'}$. It follows that $F_{\cX/\cS}$ is finite and flat.

  We next prove (i) in the general case. Since the assertion to be
  proven is of local nature we may assume $f$ factors
  \begin{displaymath}
    \cX\Xto{g}\cY:=\bA^d_\cS\Xto{p}\cS
  \end{displaymath}
  in which $g$ is quasi-\'etale and $p$ is the canonical
  projection. There is a commutative diagram
  \begin{displaymath}
    \xymatrix{
      \cX\ar[rd]^h\ar@/^/[rrd]^{F_{\cX/\cS}}\ar@/_/[rdd]_g\\
      &\cY\times_{\niv{\cY}{q}}\niv{\cX}{q}\ar[r]_{\quad p_2}\ar[d]^{p_1}
      &\niv{\cX}{q}\ar[d]^{\niv{g}{q}}\\
      &\cY\ar[r]_{F_{\cY/\cS}}&\niv{\cY}{q}
    }
  \end{displaymath}
  in which the square is Cartesian. Since $p$ is smooth, $F_{\cY/\cS}$
  is finite and flat, and by base change the same is true for
  $p_2$. Since $F_{\cX/\cS}=h\circ p_2$ it suffices to show that $h$
  is flat. Again by base change $\niv{g}{q}$ and $p_1$ are
  quasi-\'etale, and in particular quasi-unramified; then since $g$ is
  quasi-\'etale, $h$ is quasi-\'etale as well by proposition
  \ref{prop:quasi-smooth-sorites} (vii). In particular $h$ is
  quasi-smooth, and therefore flat, as required.

  We assume finally that if $f$ is formally of finite type, and show
  that $F_{\cX/\cS}$ is finite. Again the assertion is local and we
  may assume $f=p\circ g$ with the same diagram as above. By the first
  part of the argument we know that $F_{\cY/\cS}$ is finite, so by
  base change $p_2$ is also finite. It thus suffices to show that $h$
  is finite. In fact it is formally of finite type since $g$ is, and
  its construction shows that it is an adic morphism; the claim then
  follows by \cite[I Prop. 10.13.1]{EGA}.
\end{demo}

I do not know if there are quasi-smooth morphisms $\cX\to\cS$ such
that $F_{\cX/\cS}$ is not finite.

\subsection{Ordinary Differential Operators.}
\label{sec:diff-ord}

In this article we are mainly concerned with arithmetic differential
operators, but for the sake of completeness we show that the
construction of the usual (Grothendieck) ring of differential
operators extends to the case of a quasi-smooth morphism.

When $\cX/\cS$ is quasi-smooth, the rings $\cP^n_{\cX/\cS}$ are coherent,
locally free $\O_\cX$-modules, as are the $\O_\cX$-modules
\begin{displaymath}
  \Diff^n_{\cX/\cS}=\Hom_{\O_\cX}(\cP^n_{\cX/\cS},\O_\cX)
\end{displaymath}
of differential operators of order $n$. The projection maps for the
$\cP^n_{\cX/\cS}$ induce injective maps
$\Diff^{n'}_{\cX/\cS}\to\Diff^n_{\cX/\cS}$ for $n'\ge0$ and the
$\O_\cX$-module of differential operators
\begin{displaymath}
  \D_{\cX/\cS}=\limdir_n\Diff^n_{\cX/\cS}.
\end{displaymath}
is the direct limit in the category of $\O_\cX$-modules.  On any
parallelizable open affine, $\cP^n_{\cX/\cS}$ has a free basis
$\xi^I$, $|I|\le n$ where $\xi=1\ctens x-x\ctens1+I^{n+1}$ (we use the
usual multi-index notation). As usual the basis of $\Diff^n_{\cX/\cS}$
dual to $(\xi^I)_{|I|\le n}$ will be denoted by
$(\dpe{\d_i}{I})_{|I|\le n}$.

The same construction as in the algebraic case (c.f. for example
\cite[$0_{IV}$]{EGA}) gives $\D_{\cX/\cS}$ an $\O_\cX$-ring structure;
recall that this is done by dualizing a family of morphisms
\begin{displaymath}
  \delta_{n,n'}:\cP^{n+n'}_{\cX/\cS}\to\cP^n_{\cX/\cS}\tens\cP^{n'}_{\cX/\cS}
  \qquad
  x\tens y\mapsto x\tens1\tens1\tens y.
\end{displaymath}
(c.f. \cite[IV]{EGA}). Note that no completed tensor products are
involved. 

Although the lack of a category of ``quasicoherent $\O_\cX$-modules''
makes itself felt at this point, the fact that $\D_{\cX/\cS}$ is an
inductive limit of coherent $\O_\cX$-modules means that the sections
of $\D_{\cX/\cS}$ on any open affine are easily described; in
particular if $U=\Spf{A}$ is a parallelizable open affine in $\cX$
mapping to $\Spf{R}\sset\cS$, the $A$-module of sections
$D_{A/R}=\Gamma(U,\D_{U/\cS})$ is simply the free $R$ module with basis
$(\dpe{\d_i}{I})$. Thus local computations in $\D_{\cX/\cS}$ may be
done just as in the algebraic case.

If
\begin{equation}
  \label{eq:random-commutative-diagram}
  \xymatrix{
    \cX'\ar[r]^f\ar[d]&\cX\ar[d]\\
    \cS'\ar[r]&\cS
  }
\end{equation}
is a commutative diagram of locally noetherian formal schemes with
$\cX\to\cS$ and $\cX'\to\cS'$ quasi-smooth, the functoriality morphisms
\begin{displaymath}
 \md f:\Diff^n_{\cX'/\cS'}\to f^*\Diff^n_{\cX/\cS}
\end{displaymath}
are defined as usual; here the $f^*$ must be understood in the sense
of $\O$-modules on ringed spaces. The direct limit
\begin{displaymath}
  f^*\D_{\cX\to\cS}=\limdir_n f^*\Diff^n_{\cX/\cS}
\end{displaymath}
has a $(\D_{\cX'/\cS'},f^{-1}\D_{\cX/\cS})$-bimodule structure, and
the morphism
\begin{displaymath}
  \md f:\D_{\cX'/\cS'}\to f^*\D_{\cX\to\cS}
\end{displaymath}
can be used to give a construction of the left $\D_{\cX'/\cS'}$-module
structure of the inverse image $f^*M$ of a left $\D_{\cX/\cS}$-module
$M$.

\subsubsection{Stratifications.}
\label{sec:stratifications-algebraic}

Stratifications of an $\O_\cX$-module relative to $\cX/\cS$ are
defined in the same way as in the algebraic case; we review this to
set notation and terminology; we will use the same notation and
terminology for analogous concepts to be discussed later:
$m$-PD-stratifications, $m$-HPD-stratifications, and level $m$
analytic stratifications.

If $f:\cX\to\cS$ is a noetherian morphism of locally noetherian
schemes we recall from the notation section of the introduction that
$\cX_\cS(r)$ denotes the $(r+1)$-fold iterated fiber product of $\cX$
over $\cS$, and for any ordered subset $L\sset[0,1,\ldots,r]$ there is
a standard projection morphism $p_L:\cX_\cS(r)\to\cX_\cS(r')$ (where
$r'=\#L$). When $L$ is a singleton we identify $\cX_\cS(0)$ with
$\cX$. Finally the $(r+1)$-fold relative diagonal morphism
$\cX\to\cX_\cS(r)$ is defined by the multiplication map
$m(r):A_R(r)\to A$ when $\cX=\Spf{A}$ and $\cS=\Spf{R}$, and by
patching in general. We denote by $I(r)\subset A_R(r)$ the kernel of
$m(r)$ (note that $A_R(r)$ is the \textit{completed} tensor product of
$r+1$ copies of $A$).

For $n\in\bN$ we denote by $\cX_\cS^n(r)\subset\cX_\cS(r)$ the $n$th
order infinitesimal neighborhood of the diagonal. When $\cX=\Spf{A}$
and $\cS=\Spf{R}$, $\cX^n_\cS(r)$ is the affine formal scheme
associated $\hat P^n_{A/R}:=A_R(r)/I(r)^{n+1}$. For any
$L\sset[0,1,\ldots,r]$ with $\#L=r'$ the projection
$p_L:\cX_\cS(r)\to\cX_\cS(r')$ induces morphisms 
\begin{equation}
  \label{eq:canonical-projections-of-nth-order}
  p^n_L:\cX_\cS^n(r)\to\cX_\cS^n(r').
\end{equation}
In particular when $\#L=1$ the $r+1$ projections
\begin{equation}
  \label{eq:r+1-O-module-structures}
  p_i:\cX_\cS^n(r)\to\cX
\end{equation}
give $\cX_\cS^n(r)$ $r+1$ structures of a formal $\cX$-scheme, and for
$0\le i\le r$ the morphism $p_i$ is finite and formally affine. In
fact the rings $P^n_{A/R}(r)$ of $(r+1)$-fold principal parts of order
$n$ sheafify to yield sheaves of rings $\cP^n_{\cX/\cS}(r)$ on $\cX$
with $r+1$ distinct $\O_\cX$-algebra structures, which we denote by
\begin{equation}
  \label{eq:r+1-O-module-structures2}
  d_i:\O_\cX\to\cP^n_{\cX/\cS}(r).
\end{equation}
When $r=1$ we drop the $(1)$, and $\cP^n_{\cX/\cS}$ is the ring of
principal parts of order $n$ that was defined earlier.

If $M$ is an $\O_\cX$-module, 
a series of isomorphisms 
\begin{equation}
  \label{eq:stratification-algebraic}
  \chi_n:p_1^{n*}M\isom p_0^{n*}M  
\end{equation}
for $n\ge0$ will be said to be \textit{compatible} if
\begin{enumerate}
\item $\chi_0=id_M$,
\item for $n'\ge n$, the restriction of $\chi_{n'}$ to $\cX_S^n$ is
  $\chi_n$,
\end{enumerate}
and a \textit{stratification of $M$ relative to $\cS$} if the cocycle
condition holds as well:
\begin{enumerate}
\item
  $p^{n*}_{01}(\chi_n)\circ p^{n*}_{12}(\chi_n)=p^{n*}_{02}(\chi_n)$
  for all $n\ge0$.
\end{enumerate}
An $\cS$-stratified $\O_\cX$-module is an $\O_\cX$-module with an
(unspecified) stratification. An $\O_\cS$-linear morphism
$(M,\chi_n)\to(M',\chi'_n)$ of $\cS$-stratified $\O_\cX$-modules is
\textit{horizontal} if it is compatible with the stratifications.


The data of an $\cS$-stratification $\chi_n$ of $M$ is equivalent to a
family of morphisms
\begin{equation}
  \label{eq:stratification-theta}
  \theta_n=p^{n}_{1*}(\chi_n):M\to p^{n}_{1*}p^{n*}_0(M)
  =M\tens_{\O_\cX}\cP^n_{\cX/\cS}
\end{equation}
for $n\ge0$, $\O_\cX$-linear for the right structure of
$\cP^n_{\cX/\cS}$, which are compatible with the canonical morphisms
$\cP^{n'}_{\cX/\cS}\to\cP^n_{\cX/\cS}$, the identity for $n=0$, and
such that the diagram
\begin{equation}
  \label{eq:cocycle-in-terms-of-theta}
  \xymatrix{
    M\ar[d]^{\theta_{n'}}\ar[r]^{\theta_{n+n'}}
    &M\tens\cP^{n+n'}_{X/S}\ar[d]^{\delta_{n,n'}}\\
    M\tens\cP^{n'}_{X/S}\ar[r]_{\theta_n\tens1\qquad}
    &M\tens\cP^n_{X/S}\tens\cP^{n'}_{X/S}
  }
\end{equation}
commutes for all $n$, $n'\ge0$. 

The same argument as in the algebraic case shows that the category of
left $\D_{\cX/\cS}$-modules is equivalent to the category of
$\cS$-stratified $\O_\cX$-modules. The essential point is the
commutativity of \ref{eq:cocycle-in-terms-of-theta} and the fact that
the product in $\D_{\cX/\cS}$ is defined, essentially, by dualizing
the morphism $\delta_{n,n'}$ in the diagram
\ref{eq:cocycle-in-terms-of-theta}; we refer the reader to
\cite[\S2]{berthelot-ogus:1978} for the details. Finally,
$\cS$-stratified $\O_\cX$-modules have the following ``crystalline''
property: if $M$ is an $\cS$-stratified $\O_\cX$-module and $f$,
$g:\cY\to\cX$ are $\cS$-morphisms of locally noetherian formal schemes
that restrict to the same morphism $Y\to\cX$ on the reduced closed
subscheme $Y\subset\cY$, there is a canonical isomorphism
$\chi(f,g):g^*M\isom f^*M$ of $\O_\cY$-modules; it is an isomorphism
of $\D_{\cY/\cS}$-modules if $\cY\to\cS$ is quasi-smooth. The system of maps
$\chi(f,g)$ is transitive in the sense that
\begin{displaymath}
  \chi(f,g)\circ\chi(g,h)=\chi(f,h)
\end{displaymath}
for any three $f$, $g$, $h:\cY\to\cX$ reducing to the same map in the
reduced closed subscheme of $\cY$.


%% file: Dmods2.tex
\section{Arithmetic Differential Operators}
\label{sec:arith-diff}

\subsection{$m$-PD-structures.}
\label{sec:m-PD-structures}

We are fortunate that the theory of divided power ideals was worked
out in \cite{berthelot:1974} and \cite{berthelot:1990} in a very
general setting, and few modifications are needed to adapt the theory
to the formal case. We will briefly review this theory and how it is
used to construct the arithmetic differential operator rings, pointing
out the few places where perhaps something needs to be said about the
formal case. We will assume the reader is familiar with the
terminology and notation of \cite{berthelot:1996} and
\cite{berthelot:2000}, but we will summarize some of the main points
first.

From now on we fix a prime $p$, and all formal schemes will be formal
schemes over $\bZ_p$ (in addition to being noetherian). All divided
powers will be assumed compatible with the canonical divided powers of
$(p)$.

\subsubsection{Partial Divided Powers.}
\label{sec:m-PD}

Let $R$ be a commutative ring and $I\subset R$ an ideal. Recall that a
\textit{partial divided power structure of level $m$ on $I$} or an
\textit{$m$-PD-structure on $I$} is a PD-ideal $(J,\gamma)$ in $R$
such that
\begin{displaymath}
  \niv{I}{p^m}+pI\sset J\sset I.
\end{displaymath}
We also say that $(I,J,\gamma)$ is an $m$-PD-structure on $R$, that
$(R,I,J,\gamma)$ is an $m$-PD-ring and that $(I,J,\gamma)$ is an
$m$-PD-ideal in $R$. A morphism $(R,I,J,\gamma)\to(R',I',J',\gamma')$
of $m$-PD-rings is a ring homomorphism $f:R\to R'$ such that
$f(I)\sset I'$ and $f$ induces a morphism $(J,\gamma)\to(J',\gamma')$
of PD-ideals. 

We will use the following form of the definition of compatibility of
$m$-PD-structures, as in \cite[\S1.2]{berthelot:1996}:

\begin{defn}
  Suppose $R\to A$ is a homomorphism and $(\fa,\fb,\alpha)$,
  $(I,J,\gamma)$ are $m$-PD-structures on $R$ and $A$
  respectively. The $m$-PD-structures $(\fa,\fb,\alpha)$ and
  $(I,J,\gamma)$ are \textit{compatible} if the following conditions
  hold, in which $\fb_1=\fb+pR$:
  \begin{enumerate}
  \item $\fb_1A+J$ has a PD-structure inducing the PD-structures
    $\alpha$, $\gamma$ and the canonical PD-structure of $(p)$;
  \item $\fb_1A\cap I\sset\fb_1A$ is a sub-PD-ideal.
  \end{enumerate}
\end{defn}
There are a number of ways to reformulate the first condition,
c.f. \cite[Lemme 1.2.1 and Def. 1.2.2]{berthelot:1996}. The second
condition is used in the construction of the $m$-PD-adic filtration,
see \S\ref{sec:m-PD-adic-filtration} below.  

We will need the following from \cite[1.3.2 and
1.3.4]{berthelot:1996}:

\begin{lemma}\label{lemma:quotient-m-PD-structure}
  \begin{enumerate}
  \item If $K\subset A$ is an ideal, the $m$-PD-structure
    $(I,J,\gamma)$ induces an $m$-PD-structure on $I(A/K)$ such that
    $A\to A/K$ is an $m$-PD-morphism if and only if $(J+pA)\cap K$ is
    a sub-PD-ideal of $J+pA$.
  \item The induced $m$-PD-structure on $A/K$ is compatible with
    $(\fa,\fb,\alpha)$ if and only if
    \begin{enumerate}
    \item $(J+\fb_1A)\cap K$ is a sub-PD-ideal of $J+\fb_1A$, and
    \item $\fb_1A\cap (I+K)$ is a sub-PD-ideal of $\fb_1A$
    \end{enumerate}
  \end{enumerate}
\end{lemma}
In the situation of (i) we say that the $m$-PD-structure
$(I,J,\gamma)$ \textit{descends to $A/K$.}

If $(I,J,\gamma)$ is an $m$-PD-ideal in $R$, the partially divided
powers $\dpabniv{x}{k}{m}$ of an element $x\in I$ are defined by
\begin{equation}
  \label{eq:powers-and-level-m-powers}
  \dpbrniv{x}{k}{m}=x^r\gamma_q(x^{p^m})
\end{equation}
where $q$ and $r$ are integers satisfying $k=p^mq+r$ and $0\le r<p^m$.
It follows from the definition and the equality $q!\gamma_q(x)=x^q$
that
\begin{equation}\label{eq:why-theyre-called-divided-powers}
  q!\dpbrniv{x}{k}{m}=x^k.
\end{equation}
Thus if $p$ is nilpotent in $R$, $I$ is a nilideal. The partially
divided powers have a large number of formal properties which we will
not bother to state here.

\subsubsection{The $m$-PD-adic Filtration.}
\label{sec:m-PD-adic-filtration}

Any $m$-PD-ring $(R,I,J,\gamma)$ has a canonical filtration by ideals
$\dpbrshort{I}{n}\sset R$ with the following properties:
\begin{itemize}
\item $\dpbrshort{I}{0}=R$ and $\dpbrshort{I}{I}=I$.
\item $\dpbrshort{I}{n}\dpbrshort{I}{n'}\sset\dpbrshort{I}{n+n'}$.
\item If $x\in\dpbrshort{I}{n}$ then
  $\dpbrniv{x}{k}{m}\in\dpbrshort{I}{kn}$. 
\item Let $J_1=J+pR$; then for all $n\ge0$, $\dpbrshort{I}{n}\cap J_1$
  is a sub-PD-ideal of $J_1$. In particular $\dpbrshort{I}{n}\cap J$
  is a sub-PD-ideal of $J$.
\end{itemize}
The construction is quite involved and we refer the reader to
\cite[App.]{berthelot:1996}. An $m$-PD-ideal $(I,J,\gamma)$ is
\textit{$m$-PD-nilpotent} if $\dpbrshort{I}{n}=0$ for some $n$.

\subsubsection{The $m$-PD-envelope of an ideal.}
\label{sec:m-PD-envelope}

The principal construction of the theory is the $m$-PD-envelope of an
ideal. Let $(R,\fa,\fb,\alpha)$ be a ring with $m$-PD-structure, $A$
an $R$-algebra and $I\subset A$ an ideal. There is an $m$-PD-ring
\begin{equation}\label{eq:m-PD-envelope}
  (P_{(m),\alpha}(I),I^\cani,I^\canj,[\ ])
\end{equation}
and an $R$-algebra homomorphism $A\to P_{(m),\alpha}$ such that
$(I^\cani,I^\canj,[\ ])$ is compatible with $(\fa,\fb,\alpha)$ and
having the following universal property: for any $A$-algebra $B$
with an $m$-PD-structure compatible with $(\fa,\fb,\alpha)$, the
structure morphism $A\to B$ factors uniquely through an
$m$-PD-morphism $P_{(m),\alpha}(I)\to B$. The ideal $I^\cani$
(resp. $I^\canj$) is called the \textit{canonical $m$-PD-ideal}
(resp. the \textit{canonical PD-ideal}) of the $m$-PD-envelope
$P_{(m),\alpha}(I)$. 

The quotient of $P_{(m),\alpha}$ by the $(n+1)$-st step of the
$m$-PD-adic filtration is denoted by $P^n_{(m),\alpha}(I)$; it has a
similar universal property with respect to homomorphisms of $A$ to
$m$-PD-nilpotent $m$-PD-rings.  From
\ref{eq:powers-and-level-m-powers} we see that the image of $I$ in
$P^n_{(m),\alpha}(I)$ is a nilideal.

The formation of $P_{(m),\alpha}$ and $P^n_{(m),\alpha}$ commutes with
flat base change: if $A\to A'$ is flat, the natural homomorphisms
$A'\tens_AP_{(m),\alpha}(I)\to P_{A',(m)}(IA')$ and
$A'\tens_AP^n_{(m),\alpha}(I)\to P^n_{A',(m)}(IA')$ are isomorphisms.

\subsubsection{The Regular Case.}
\label{sec:m-PD-regular}

The case when $I\subset A$ is a regular ideal (Zariski-locally
generated by a regular sequence) is particularly nice, and
particularly important. The algebras $P^n_{(m),\alpha}(I)$ are
independent of the $m$-PD-structure of $R$, flat over $R$ and their
formation commutes with arbitrary base change $R\to R'$. Suppose
furthermore that $I$ is generated by a regular sequence
$x_1,\ldots,x_d$, that $A/I$ is flat over $R$ and that the quotient
map $A\to A/I$ has a section $\sigma:A/I\to A$. Then via $\sigma$,
$P^n_{(m),\alpha}(I)$ is a free $A/I$-module on the $m$-PD-polynomials
$\dpbrniv{x}{K}{m}$ for $|K|\le n$ (in the usual multi-index
notation). Furthermore the image of $I^\cani$ in $P^n_{(m),\alpha}(I)$
is free on the $\dpbrniv{x}{K}{m}$ for $|K|>0$, and $I^\canj$ is
generated by $pP^n_{(m),\alpha}(I)$ and by the $\dpbrniv{x}{K}{m}$ for
those $K=(k_1,\ldots,k_d)$ for which at least one entry is $\ge p^m$.

If $p$ is nilpotent in $A$ these assertions hold for the full
$m$-PD-envelope $P_{(m),\alpha}(I)$, the only modification being that
when $A\to A/I$ has a section and $x_1,\ldots,x_d$ is a regular
sequence generating $I$, the $A/I$-module $P_{(m),\alpha}(I)$ is free
on the entire set of $\dpbrniv{x}{K}{m}$.

An important example is the $m$-PD-polynomial algebra
$R\<X_1,\ldots,X_d\>$, defined as the $m$-PD-envelope of the regular
ideal $(X_1,\ldots,X_d)\subset R[X_1,\ldots,X_d]$. Elements of
$R\<X_1,\ldots,X_d\>$ are called $m$-PD-polynomials; as an example of
their use in computation we recall, from the proof of
\cite[Prop. 4.2.1]{berthelot:1996} that for any natural number $r$
divisible by $p^{m+1}$ there is an $m$-PD-polynomial
$\niv{\varphi}{m}_r(X_1,X_2)$ such that for all $t_1$, $t_2$ in some
$m$-PD-ring $(R,I,J,\gamma)$ such that $t_1-t_2\in I$,
\begin{equation}
  \label{eq:varphi}
  t_2^r-t_1^r=p\niv{\varphi}{m}_r(t_1,t_2).
\end{equation}
We may work in $\bZ_{(p)}[X_1]\<X_2-X_1\>$, in which case, writing
$r=p^{m+1}q$, we see that
\begin{align*}
  X_2^r-X_1^r&=((X_1+(X_2-X_1))^{p^{m+1}})^q-X_1^r\\
  &=(X_1^{p^{m+1}}+p(*)+(X_2-X_1)^{p^{m+1}})^q-X_1^r\\
  &=(X_1^{p^{m+1}}+p(*)+p!\dpbrniv{(X_2-X_1)}{p^{m+1}}{m})^q-X_1^r\\
\end{align*}
from which the assertion follows. The identities
\begin{equation}
  \label{eq:varphi-identities}
  \niv{\varphi}{m}_r(X,X)=0,\qquad
  \niv{\varphi}{m}_r(X_1,X_2)+\niv{\varphi}{m}_r(X_2,X_3)
  =\niv{\varphi}{m}_r(X_1,X_3)
\end{equation}
can be proven by reduction to the case of
$\bZ_{(p)}[X_1]\<X_2-X_1\>$. Since the latter has no $p$-torsion,
these identities can be checked after multiplication by $p$, in which
case they are obvious consequences of \ref{eq:varphi}.

\subsubsection{Application to Formal Schemes.}
\label{sec:m-PD-formal-schemes}

The construction of $m$-PD-envelopes sheafifies on a scheme because it
commutes with flat base change, and localizations are flat. Thus if
$S$ is an $m$-PD-scheme and $X$ is an $S$-scheme, any ideal
$I\subset\O_X$ has an $m$-PD-envelope $\cP_{(m),\alpha}(I)$; it is a
quasicoherent $\O_X$-module satisfying the same universal property as
in the affine case. The same holds for the $\cP^n_{(m),\alpha}(I)$.

Suppose now $\cS$ is an locally noetherian formal scheme with an
$m$-PD-structure $(\fa,\fb,\alpha)$, and $\cX$ is a locally noetherian
formal $\cS$-scheme. Since formal localizations are also flat, one
might expect that the construction of $m$-PD-envelopes sheafifies in
the same way. However the lack of a category of quasicoherent
$\O_\cX$-modules makes itself felt at this point: there is no analogue
here of the sheafification procedure that is available in complete
generality for schemes. If $\cI\subset\O_\cX$ is an ideal one can of
course sheafify the presheaf of divided power envelopes of $\cI$ on
affines; the trouble starts when one tries to prove that the ring of
sections of this sheaf over an affine open is a divided power
envelope.

The situation is somewhat better for regular ideals, when one is
concerned only with the truncated divided power envelopes. Suppose
$\cI\subset\O_\cX$ is a regular ideal, and that the closed immersion
$\cY\to\cX$ defined by $\cI$ has a retraction $\cX\to\cY$ (i.e. the
quotient map $\O_\cX\to\O_\cX/\cI$ has a section). Let
$U=\Spf{A}\sset\cX$ and $V=\Spf{R}\sset\cS$ be open affines such that
$\cX\to\cS$ maps $U$ into $V$, and $U\to V$ is parallelizable. If we
set $I=\Gamma(U,\cI)$, then by hypothesis $A\to A/I$ has a section
$\sigma$, and we know that $P^n_{(m),\alpha}(I)$ is a free
$A/I$-module of finite rank via the section $\sigma$. For any
$f\in\Gamma(U,\O_\cX)$ the morphism $A\to A_{\{f\}}$ is flat and the
canonical morphism
\begin{displaymath}
  A_{\{f\}}\tens_AP^n_{(m),\alpha}(I)\to P^n_{(m),\alpha}(IA_{\{f\}})
\end{displaymath}
is an isomorphism. Since $P^n_{(m),\alpha}(I)$ is finitely generated
the tensor product may be replaced by a completed tensor
product. Finally, the sections of $\O_\cX\to\O_\cX/\cI$ being used all
come from a single global section. We conclude that there is a
coherent $\O_\cY$-module $\cP^n_{(m),\alpha}(\cI)$ with the property
that
\begin{displaymath}
  \Gamma(U,\cP^n_{(m),\alpha}(\cI))\simeq P^n_{(m),\alpha}(I)
\end{displaymath}
when $U$ is affine and $I=\Gamma(U,\cI)$. Since $\cX$ is locally
noetherian, the ideal $\cI$ is locally nilpotent in
$\cP^n_{(m),\alpha}(\cI)$, and therefore $\cP^n_{(m),\alpha}(\cI)$ is
supported on the closed formal subscheme of $\cX$ defined by $\cI$.

Like its affine counterpart, the $\O_\cX$-algebra
$\cP^n_{(m),\alpha}(\cI)$ has a universal property, best expressed by
introducing the formal scheme
\begin{displaymath}
  \cX^n_{(m),\alpha}(\cI)=\Spfalg{\O_\cX}{\cP^n_{(m),\alpha}(\cI)}.
\end{displaymath}
Suppose $f:\cX'\to\cX$ is an $\cS$-morphism of an adic noetherian
formal schemes, and $\cX'$ has an $m$-PD-structure $(\cI',\cJ',\gamma')$
compatible with $(\fa,\fb,\alpha)$ and nilpotent of order $n$. If
$f^*\cI\sset\cI'$, $f$ has a unique factorization
\begin{displaymath}
  \cX'\Xto{g}\cX_{(m),\alpha}(\cI)\Xto{p}\cX
\end{displaymath}
for some $m$-PD-morphism $g$, where $p$ is the morphism corresponding
to the structure map of the coherent $\O_\cX$-algebra
$\cP^n_{(m),\alpha}(\cI)$.

The case of the full $m$-PD-envelope is more difficult, and requires
further hypotheses and no small amount of technicalities. We will deal
with it in section \ref{sec:P_m(I)}.

\begin{lemma}\label{lemma:passing-m-PD-to-quotient}
  Let $A$ be an noetherian ring with $m$-PD-structure
  $(I,J,\gamma)$. For any ideal $K\subset A$, the $m$-PD-structure
  $(I,J,\gamma)$ descends to $A/K^n$ for all sufficiently large $n$.
\end{lemma}
\begin{demo}
  By \ref{sec:m-PD-structures}.3 we need that
  $K^n\cap(J+pA)\sset J+pA$ is a sub-PD-ideal for all $n\gg0$,
  i.e. $\gamma_k(x)\in K^n\cap(J+pA)$ for all $x\in K^n\cap(J+pA)$ and
  $k>0$, and since $\gamma_1(x)=x$ we may assume $k>1$. By Artin-Rees
  there is an integer $c$ such that
  \begin{displaymath}
    K^n\cap(J+pA)\sset K^{n-c}(K^c\cap(J+pA))
  \end{displaymath}
  for $n>c$. Then for all $x\in K^n\cap(J+pA)$,
  \begin{displaymath}
    \gamma_k(x)\in K^{k(n-c)}(J+pA)
  \end{displaymath}
  by the basic PD-identities. We are done if $k(n-c)\ge n$, and since
  $k\ge2$ this holds when $n\ge 2c$.
\end{demo}

\begin{remark}
  The conclusion of the lemma can be restated as follows: $A$ has an
  cofinal set of ideals of definition $K$ such that the
  $m$-PD-structure of $A$ descends to $A/K^n$ for \textit{all} $n$; it
  suffices to replace $K$ by $K^N$ for all sufficiently large $N$.
\end{remark}

\subsection{The ring $\niv{\D}{m}_{\cX/\cS}$.}
\label{sec:diff-arith}

Suppose now $\cX\to\cS$ is quasi-smooth.  The considerations of the
last paragraph apply to the sheaf of rings $\O_{\cX_\cS(r)}$ and its
diagonal ideal $\cI(r)$.

\subsubsection{Principal parts of level $m$}
\label{sec:principal-parts-level-m}

We denote by $\cP^n_{\cX/\cS,(m)}(r)$ the $m$-PD-envelope of order $n$
of $\cI(r)$ in $\O_{\cX_\cS(r)}$; it is supported on the diagonal of
$\cX_\cS(r)$ and may be regarded as a sheaf of rings on $\cX$. It has
$r+1$ $\O_\cX$-algebra structures, with respect to each of which it is
a coherent locally free $\O_{cX_\cS(r)}$-algebra. As before we drop
the $(r)$ when $r=1$. We reuse the notation
\begin{equation}
  \label{eq:canonical-projections-mPD}
    d^n_K:\cP^n_{\cX/\cS,(m)}(r)\to\cP^n_{\cX/\cS,(m)}(r')
\end{equation}
of \S\ref{sec:stratifications-algebraic} for the canonical
projections; their existence follows from the universal property of
the truncated $m$-PD-envelopes. If we define
\begin{equation}
  \label{eq:level-m-diagonal-order-n}
  \cX^n_{\cS,(m)}(r)=\Spfalg{\O_\cX}{\cP^n_{\cX/\cS,(m)}(r)}
\end{equation}
then $\cX^n_{\cS,(m)}(r)$ is a formally affine formal $\cX$-scheme for
each of its $r+1$ $\O_\cX$-algebra structures. The $r+1$ structure
morphisms are finite and the homomorphisms
\ref{eq:canonical-projections-mPD} induce morphisms
\begin{equation}
  \label{eq:canonical-projections-and-inclusions}
  p^n_K:\cX_{\cS,(m)}^n(r')\to\cX_{\cS,(m)}^n(r)
\end{equation}
As before when $r=1$ we drop the $(1)$.

The $\O_\cX$-module of level $m$ operators of order $\le n$ is defined
to be
\begin{displaymath}
  \Diff^n_{\cX/\cS,(m)}=Hom_{\O_\cX}(\cP^n_{\cX/\cS,(m)},\O_\cX).
\end{displaymath}
As in \cite[\S2.2]{berthelot:1990} the canonical projections
$\cP^{n'}_{\cX/\cS,(m)}\to\cP^n_{\cX/\cS,(m)}$ for
$n'\ge n$ induce injections
$\Diff^n_{\cX/\cS,(m)}\to\Diff^{n'}_{\cX/\cS,(m)}$ and the ring of
\textit{arithmetic differential operators of level $m$} is
\begin{displaymath}
  \niv{\D}{m}_{\cX/\cS}=\limdir_n\Diff^n_{\cX/\cS,(m)}
\end{displaymath}
where again the limit is to be understood in the sense of
$\O_\cX$-modules on the ringed space $(|\cX|,\O_\cX)$. When
$\cX=\Spf{B}$ and $\cS=\Spf{A}$ are affine, we define
$\Gamma(\cX,\Diff^n_{\cX/\cS,(m)})=\Diff^n_{A/R,(m)}$ and
$\Gamma(\cX,\niv{\D}{m}_{\cX/\cS})=\niv{D}{m}_{B/A}$, and then
\begin{displaymath}
  \niv{D}{m}_{B/A}=\limdir_n\Diff^n_{A/R,(m)}
\end{displaymath}
when $\cX$ is noetherian (so that taking global sections commutes with
the inductive limit). When $\cX\to\cS$ is parallelizable with local
coordinates $x_1,\ldots,x_d$, $\cP^n_{\cX/\cS,(m)}$ is the free module
with basis $(\dpniv{\xi}{I}{m})_{|I|\le n}$, where as usual
$\xi_i=1\ctens x_i-x_i\ctens1$. The dual basis of
$\Diff^n_{\cX/\cS,(m)}$ is denoted by
$\{\dpabniv{\d}{K}{m}\}_{|K|\le n}$.

We make $\niv{\D}{m}_{\cX/\cS}$ into a ring by means of a map
\begin{displaymath}
  \delta^{n,n'}_{(m)}:\cP^{n+n'}_{\cX/\cS,(m)}\to
  \cP^n_{\cX/\cS,(m)}\tens_{\O_\cX}\cP^{n'}_{\cX/\cS,(m)} 
\end{displaymath}
arising from $\delta^{n,n'}$ via the universal property of
$m$-PD-envelopes. The formal properties of the ring
$\niv{\D}{m}_{\cX/\cS}$ are the same as in
\cite[\S2.2]{berthelot:1996}, and are proven in the same way; we will
not bother to state them here. 

Once again we should point out that in the situation of the diagram
\ref{eq:random-commutative-diagram}, the
$(\niv{\D}{m}_{\cX'/\cS'},f^{-1}\niv{\D}{m}_{\cX/\cS})$-bimodule
$f^*\niv{\D}{m}_{\cX/\cS}$ is to be understood in the sense of ringed
spaces. This will not be an issue, as we will see later.

For $m'\ge m$ and all $n\ge0$ there is a canonical $m$-PD-morphism
\begin{equation}\label{eq:Pm-change-of-level}
  \iota_{m',m}^n:\cP^n_{\cX/\cS,(m')}\to\cP^n_{\cX/\cS,(m)}
\end{equation}
arising from the universal property by regarding the canonical
$m'$-PD-ideal of $\cP^n_{\cX/\cS,(m')}$ as an $m$-PD-ideal. Since
$\cP^n_{\cX/\cS,(m)}$ is generated by the $\dpbrniv{x}{k}{m}$ for
$0\le k<n$, $\iota_{m',m}^n$ is characterized by the formula
\begin{equation}\label{eq:Pm-change-of-level-explicit}
  k=qp^m+r=q'p^{m'}+r'\implies
  \iota_{m',m}^n(\dpbrniv{x}{k}{m})=\frac{q!}{q'!}\dpbrniv{x}{k}{m'}
\end{equation}
which follows from \ref{eq:why-theyre-called-divided-powers} by
reduction to the universal case. Dualizing
\ref{eq:Pm-change-of-level-explicit} and taking the inductive limit in
$n$ results in a ring homomorphism
\begin{equation}
  \label{eq:Dm-change-of-level}
  \rho_{m',m}:\niv{\D}{m}_{\cX/\cS}\to \niv{\D}{m'}_{\cX/\cS}
\end{equation}
for all $m'\ge m$. In local coordinates it is given by the formula
\begin{equation}\label{eq:Dm-change-of-level-explicit}
  \rho_{m',m}(\dpabniv{\d}{K}{m})=\frac{Q!}{Q'!}\dpabniv{\d}{K}{m'}
\end{equation}
where $Q$, $Q'$ are defined by
\begin{displaymath}
  K=p^mQ+R=p^{m'}Q'+R'
\end{displaymath}
with $0\le R<p^m$ and $0\le R'<p^{m'}$.

\subsubsection{Base change.}
\label{sec:base-change}

Let
\begin{equation}
  \label{eq:base-change}
  \xymatrix{
    \cX'\ar[r]^f\ar[d]&\cX\ar[d]\\
    \cS'\ar[r]&\cS
  }
\end{equation}
be a commutative diagram of locally noetherian formal schemes, with
$\cX'\to\cS'$ and $\cX\to\cS$ quasi-smooth and $\cS'\to\cS$ an
$m$-PD-morphism. For any left $\niv{\D}{m}_{\cX/\cS}$-module $M$ there
are, as in \cite[2.2.2]{berthelot:1996} and
\cite[\S2.1]{berthelot:2000} two equivalent ways of placing left
$\niv{\D}{m}_{\cX'/\cS'}$-module structure on $f^*M$. One is via the
natural homomorphism
\begin{equation}
  \label{eq:base-change-level-m}
  \md:\niv{\D}{m}_{\cX'/\cS'}\to\niv{\D}{m}_{\cX'\to\cX}=f^*\niv{\D}{m}_{\cX/\cS}
\end{equation}
which is deduced from the natural homomorphisms
\begin{displaymath}
  f^*\cP^n_{\cX/\cS}\to\cP^n_{\cX'/\cS'}
\end{displaymath}
by duality and passage to the limit; the induced
$(\niv{\D}{m}_{\cX'/\cS'},f^{-1}\niv{\D}{m}_{\cX/\cS})$-bimodule
structure on $\niv{\D}{m}_{\cX/\cS}$ yields a left
$\niv{\D}{m}_{\cX'/\cS'}$-module structure via the canonical
isomorphism
\begin{displaymath}
  \niv{\D}{m}_{\cX'\to\cX}\tens_{f^{-1}\niv{\D}{m}_{\cX/\cS}}M  
  \isom f^*M
\end{displaymath}
(c.f. \cite[2.1.3]{berthelot:2000}). The other method is to use the
commutative diagrams
\begin{displaymath}
  \xymatrix{
    (\cX')^n_\cS\ar[d]_{p_i}\ar[r]^{f\times f}&\cX^n_\cS\ar[d]^{p_i}\\
    \cS'\ar[r]&\cS
  }
\end{displaymath}
for $i=0$, $1$
to show that an $m$-PD-stratification of $M$ relative to $\cS$ pulls
back to an $m$-PD-stratification of $f^*M$ relative to $\cS'$. The
latter method is perhaps more convenient for proving the transitivity
formula $(fg)^*M\simeq g^*f^*M$, c.f. \cite[2.1.1]{berthelot:2000}. On
this point nothing needs to be added to the treatment of
\cite[2.1.3]{berthelot:2000} and \cite[2.1.1]{berthelot:2000}. 

\subsubsection{$m$-PD-stratifications.}
\label{sec:m-PD-stratifications}

An $m$-PD-stratification relative to $\cS$ of an $\O_\cX$-module $M$ is
defined just as before, but with the $\cX^n_{\cS,(m)}(r)$ in place of
the $\cX^n_{\cS}(r)$: it is a series of isomorphisms
\begin{equation}
  \label{eq:m-PD-stratification1}
  \chi_n:p^{n*}_1(M)\isom p^{n*}_0(M)
\end{equation}
satisfying the conditions \ref{sec:stratifications-algebraic}.4--6
(same conditions on different maps!). More generally, a sequence
$\{\chi_n\}_{n\ge0}$ is \textit{compatible} if it satisfies
\ref{sec:stratifications-algebraic}.4--5. The isomorphisms
\ref{eq:m-PD-stratification1} can also be given as a series of
isomorphisms
\begin{equation}
  \label{eq:m-PD-stratification2}
  \chi_n:\cP^n_{\cX/\cS,(m)}\tens_{\O_\cX}M\isom
  M\tens_{\O_\cX}\cP^n_{\cX/\cS,(m)}
\end{equation}
of $\O_\cX$-modules with analogous properties, or via the adjunction
isomorphism as a series of morphisms
\begin{equation}
  \label{eq:m-PD-stratification3}
  \theta_n:M\to M\tens_{\O_\cX}\cP^n_{\cX/\cS,(m)}
\end{equation}
that are $\O_\cX$-linear for the right structure of
$\cP^n_{\cX/\cS,(m)}$, compatible with the canonical morphisms
$\cP^{n'}_{\cX/\cS,(m)}\to\cP^n_{\cX/\cS,(m)}$ for $n'\ge n$, the
identity for $n=0$, and making commutative a diagram analogous to
\ref{eq:cocycle-in-terms-of-theta}.

If $P$ is a level $m$ differential operator of order $\le n$, we can
compose $\theta_n$ with $1\tens P$ to get a map
\begin{displaymath}
  M\Xto{\theta_n} M\tens_{\O_\cX}\cP^n_{\cX/\cS,(m)}\Xto P M
\end{displaymath}
and the usual argument shows that this induces a left $\niv\D
m_{\cX/\cS}$-module structure on $M$. As in \cite{berthelot:1996}, one
proves: 

\begin{prop}\label{prop:m-PD-stratifications-and-Dm-module-structure}
  The category of left $\niv\D m_{\cX/\cS}$-modules is equivalent to
  the category the category of $\O_\cX$-modules endowed with an
  $m$-PD-stratification.
\end{prop}

\subsection{The ring $\niv{\hD}{m}_{\cX/\cS}$.}
\label{sec:completions}

In the setting of \cite{berthelot:1996}, the next step in the theory
is to form the $p$-adic completion $\niv{\hD}{m}_{\cX/\cS}$ of
$\niv{\D}{m}_{\cX/\cS}$, and then take the inductive limit of the
$\niv{\hD}{m}_{\cX/\cS\bQ}$ to get the full ring of arithmetic
differential operators. We have explained in the introduction why this
is not the thing to do here, and the reader will see a case of this in
the proof of theorem \ref{thm:coherence-finite-level}. Instead we must
complete $\niv{\hD}{m}_{\cX/\cS}$ with respect to an ideal of
definition; this does not obviously result in a sheaf of rings, and
what makes this idea workable is the fact that rings like
$\niv{\D}{m}_{\cX/\cS}$ are particularly rich in two-sided ideals, as
we see from corollary \ref{cor:lots-of-bilateralising-ideals}.

We    first    recall    some    definitions    and    results    from
\cite{nastasescu-van-oystaeyen:1982}             (c.f.            also
\cite[\S3.2]{berthelot:1996}). An  (left or right) ideal  $I\subset R$
in a ring is  \textit{central} if it is generated by  a set of central
elements,  and   \textit{centralising}  if   it  is  generated   by  a
centralising sequence, i.e. a sequence $x_1,\ldots,x_n\in R$ such that
for all $i$ the image of $x_i$ in $R/(x_1,\ldots,x_{i-1})$ lies in the
center. A  centralising ideal  is evidently 2-sided,  and when  $R$ is
noetherian the  standard results  from commutative  algebra concerning
$I$-adic  topologies   and  completions   extend  to  the   case  when
$I\subset  R$  is  centralising.  We   refer  the  reader  to  \cite[D
III]{nastasescu-van-oystaeyen:1982}  for   proofs  of   the  following
assertions, in  which $R$ is any  left and right noetherian  ring with
unit and $I$ is a centralising ideal:
\begin{itemize}
\item The Artin-Rees lemma holds in the following form: if $M$ is a
  finitely generated left $R$-module and $N$ is submodule of $M$,
  there is a function $f:\bN\to\bN$ such that
  \begin{displaymath}
    N\cap I^{f(n)}M\sset I^nN.
  \end{displaymath}
\item The $I$-adic completion functor functor $M\mapsto\hat M$ is
  exact on the category of finitely generated left (or right)
  $R$-modules.
\item The completion $\hat R$ is left and right flat over $R$, and 
\item The natural map $\hat R\tens_RM\to\hat M$ is an isomorphism for
  any finitely generated left $R$-module $M$.
\end{itemize}
If $I$ is a central ideal, the function $f$ may be taken to be of the
form $f(n)=c+n$ for some constant $c$, and the proof of the Artin-Rees
lemma is the same as in the commutative case. The general case is more
complicated, c.f. \cite[D V]{nastasescu-van-oystaeyen:1982}. The
remaining statements follow from the first by the usual arguments.

We can relativize the notion of a centralising ideal. Suppose $R\to A$
is a ring homomorphism with $R$ commutative; we say that an ideal
$J\subset R$ is \textit{centralising in $A$} if it is generated by a
sequence whose image in $A$ is centralising. Thus if $J\subset R$ is
centralising in $A$, $JA=AJ$ is a centralising ideal of $A$.

\begin{lemma}\label{lemma:constructing-centralising-ideals}
  Let $R\to A$ be a ring homomorphism with $R$ commutative and
  noetherian. Let $I\subset R$ be an ideal such that $IA$ is a central
  ideal in $A$. Let $J\subset R$ be any ideal and denote by $f:A\to
  A/IA$ the canonical homomorphism. The ideal $J'=J\cap
  f^{-1}(Z(A/I))$ is centralising in $A$.
\end{lemma}
\begin{demo}
  This is an immediate consequence of the definitions; the noetherian
  hypothesis is there to ensure that $J'$ is finitely generated, which
  is implicit in the definition.
\end{demo}

\begin{prop}\label{prop:existence-of-centralising-ideals}
  Suppose $\cX/\cS$ is quasi-smooth and $m\ge0$. Any ideal of
  definition containing the prime $p$ contains an ideal of definition
  that is centralising in $\niv{\D}{m}_{\cX/\cS}$.
\end{prop}
\begin{demo}
  Suppose first that $\cX$ and $\cS$ are affine. Let $J\subset\O_\cX$
  be an ideal of definition and write $J=(p,f_1,\ldots,f_n)$. The
  ideal $J'$ constructed in lemma
  \ref{lemma:constructing-centralising-ideals} with $I=(p)$ is
  topologically nilpotent since it is contained in $J$ and open since
  it contains $(p,f_1^{p^{m+1}},\ldots,f_n^{p^{m+1}})$, by
  \cite[prop. 2.2.6]{berthelot:1996}, c.f. also the remark at the end
  of \cite[\S3.2.3]{berthelot:1996}. It is therefore an ideal of
  definition.

  If $J=(p,f_1,\ldots,f_n)=(p,g_1,\ldots,g_r)$ then
  \begin{displaymath}
  (p,f_1^{p^{m+1}},\ldots,f_n^{p^{m+1}})=(p,g_1^{p^{m+1}},\ldots,g_r^{p^{m+1}}),
  \end{displaymath}
  so this construction globalizes.
\end{demo}

\begin{defn}\label{defn:bilateralising-ideal}
  Suppose $R$ is a commutative ring and $R\to A$ is an $R$-ring. We will
  say that an ideal $I\subset R$ is \textit{bilateralising in $A$} if
  $IA=AI$. 
\end{defn}
If the reference to $A$ is clear we will simply say that $I$
\textit{bilateralising}. The following assertions are immediate:
\begin{itemize}
\item If $I$ is bilateralising in $A$, $IA\subset A$ is a 2-sided ideal.
\item Sums and products of bilateralising ideals are
  bilateralising. In particular, powers of bilateralising ideals are
  bilateralising. 
\item A centralising ideal in $A$ is bilateralising.
\item If $M$ is a left $A$-module and $I\subset R$ is bilateralising,
  $M/IM$ is a left $A/IA$-module.
\end{itemize}

\begin{cor}\label{cor:lots-of-bilateralising-ideals}
  The set of ideals of definition of $\O_\cX$ that are bilateralising
  in $\niv{\D}{m}_{\cX/\cS}$ is cofinal in the set of all ideals of
  definition.
\end{cor}
\begin{demo}
  In fact the lemma says that centralising ideals of definition exist,
  and if $J$ is one such, $\{J^n\}_{n\ge0}$ is a cofinal system of
  bilateralising ideals of definition.
\end{demo}

\begin{defn}\label{defn:m-bilateralising}
  Suppose $\cS$ is a formal scheme with an $m$-PD-structure and
  $\cX\to\cS$ is quasi-smooth. An ideal $J\sset\O_\cX$ is
  \textit{$m$-bilateralising} if it is bilateralising in
  $\niv{\D}{m}_{\cX/\cS}$
\end{defn}
As before, if $m$ is understood we will simply say that $J$ is
bilateralising.

It is easy to characterize the $m$-bilateralising ideals of
$\O_\cX$. We first recall the level $m$ Leibnitz identity
\begin{equation}
  \label{eq:commutation-in-Dm}
  \dpabniv{\d}{K}{m}f=\sum_{I+J=K}\abinom{K}{I}{m}
  \dpabniv{\d}{I}{m}(f)\dpabniv{\d}{J}{m}
\end{equation}
which is \cite[Prop. 2.2.4, (iv)]{berthelot:1996}.

\begin{prop}\label{prop:bilateral-iff-horizontal}
  An ideal $J\subset\O_\cX$ is $m$-bilateralising if and only if it is
  horizontal for $\niv{\D}{m}_{\cX/\cS}$, i.e.
  $\niv{\D}{m}_{\cX/\cS}J\sset J$.
\end{prop}
\begin{demo}
  As usual we may reduce to a parallelizable affine situation. If the
  $\dpabniv{\d}{K}{m}$ are the basic differential operators
  corresponding to a choice of local coordinates, the assertion that
  $J$ is $m$-bilateralising is equivalent to the containments
  \begin{displaymath}
    f\dpabniv{\d}{K}{m}\sset\niv{D}{m}J
    \qquad
    \dpabniv{\d}{K}{m}f\sset J\niv{D}{m}
  \end{displaymath}
  for all $f\in J$ and $K\in\bN^d$. Thus if $J$ is $m$-bilateralising
  and $f\in J$,
  \begin{displaymath}
    \dpabniv{\d}{K}{m}f-f\dpabniv{\d}{K}{m}\in J\niv{D}{m}
  \end{displaymath}
  for all $K\in\bN^d$. Since $\niv{D}{m}$ is free on the
  $\dpabniv{\d}{K}{m}$, \ref{eq:commutation-in-Dm} shows that
  $\dpabniv{\d}{K}{m}(f)\in J$. Conversely if $\niv{D}{m}J\sset J$,
  the containment $\dpabniv{\d}{K}{m}(f)\sset J\niv{D}{m}$ is an
  immediate consequence of \ref{eq:commutation-in-Dm}, and the
  containment $\dpabniv{\d}{K}{m}(f)\sset\niv{D}{m}J$ follows from
  \ref{eq:commutation-in-Dm} by induction on $|K|$.
\end{demo}

Since locally $\niv{\D}{m}_{\cX/\cS}$ is generated locally by the
$\dpabniv{\d}{K}{m}$ for $|K|\le p^{m}$, we deduce:

\begin{cor}\label{cor:bilateralising-descends}
  If $J\sset\O_\cX$ is $m'$-bilateralising and $m'\ge m$ then $J$ is
  $m$-bilateralising.\nodemo
\end{cor}


\begin{cor}\label{cor:bilaterals-and-envelopes}
  If $J\sset\O_\cX$ is $m$-bilateralising then
  \begin{displaymath}
    J\cP^n_{\cX/\cS,(m)}=\cP^n_{\cX/\cS,(m)}J
  \end{displaymath}
  for all $n\ge0$.
\end{cor}
\begin{demo}
  For $f\in\O_\cX$ we set $\delta^n(f)=d^n_1(f)-d^n_0(f)$, and write
  $\cP^n$ for $\cP^n_{\cX/\cS,(m)}$. It suffices to show that
  $\delta^n(J)\sset\cP^nJ\cap J\cP^n$ for all $n$. The Taylor formula
  says that
  \begin{displaymath}
    \delta^n(f)=\sum_{0<|K|\le n}d_0(\dpabniv{\d}{K}{m}(f))\dpniv{\xi}{K}{m}
  \end{displaymath}
  which by proposition \ref{prop:bilateral-iff-horizontal} yields
  $\delta^n(J)\sset J\cP^n$ for all $n$. It also shows that
  \begin{displaymath}
    \delta^n(f)=\sum_{0<|K|\le
      n}d_1(\dpabniv{\d}{K}{m}(f))\dpniv{\xi}{K}{m}
    -\sum_{0<|K|\le n}\delta^n(\dpabniv{\d}{K}{m}(f))\dpniv{\xi}{K}{m}
  \end{displaymath}
  from which we deduce that
  \begin{displaymath}
    \delta^n(J)\sset\cP^nJ+I^\cani\delta^n(J)
  \end{displaymath}
  and the result follows by iteration since $I^\cani$ is nilpotent in
  $\cP^n$. 
\end{demo}

If $J\subset\O_\cX$ is any open ideal we denote by $X_J$ the
(ordinary) scheme $(|\cX|,\O_\cX/J)$.  If $J$ is $m$-bilateralising we
set
\begin{equation}
  \label{eq:nivD-XJ}
  \niv{\D}{m}_{X_J/\cS}=\niv{\D}{m}_{\cX/\cS}/J\niv{\D}{m}_{\cX/\cS}
\end{equation}
(the notation is purely formal, since $X_J$ is not quasi-smooth over
$\cS$) which may be regarded indifferently as a sheaf of
$\O_\cX$-rings, or of $\O_{X_J}$-rings on $X_J$. Via the latter
structure, it is clearly a quasicoherent $\O_{X_J}$-module; in fact on
any parallelizable open $U\sset\cX$, $\niv{\D}{m}_{\cX/\cS}$ is free
$\O_\cX$-module on the $\dpabniv{\d}{I}{m}$, so that
$\niv{\D}{m}_{X_J/\cS}$ is a free $\O_{X_J}$-module on the images of
the $\dpabniv{\d}{I}{m}$ (for which we use the same notation).

If $J'\sset J$ is another open $m$-bilateralising ideal there is an
evident homomorphism
\begin{displaymath}
  \niv{\D}{m}_{X_{J'}/\cS}\to\niv{\D}{m}_{X_J/\cS}
\end{displaymath}
of $\O_\cX$-rings, inducing an isomorphism
\begin{displaymath}
  \niv{\D}{m}_{X_{J'}/\cS4}\tens_{\O_{J'}}\O_{X_J}
  \isom\niv{\D}{m}_{X_J/\cS}
\end{displaymath}
of $\O_{X_J}$-rings. The $\O_\cX$-ring $\niv{\hD}{m}_{\cX/\cS}$ is the
inverse limit
\begin{equation}
  \label{eq:Dhatm}
  \niv{\hD}{m}_{\cX/\cS}=\liminv_J\niv{\D}{m}_{X_J/\cS}
\end{equation}
where $J$ runs through any cofinal set of ideals of definition of
$\O_\cX$ bilateralising in $\niv{\D}{m}_{\cX/\cS}$. On any
parallelizable open affine $\Spf{A}=U\sset\cX$, elements of
\begin{displaymath}
  \niv{\hat D}{m}_{A/R}=\Gamma(U,\niv{\D}{m}_{\cX/\cS})
\end{displaymath}
may be identified with formal series
$\sum_{I\in\bN^d}a_I\dpabniv{\d}{I}{m}$ with $a_I\to 0$ in the adic
topology of $A$. 

\begin{remark}
  The proof of corollary \ref{cor:lots-of-bilateralising-ideals} shows
  that the inverse limit in \ref{eq:Dhatm} is actually a $J$-adic
  completion for any centralising ideal of definition of $\O_\cX$. In
  particular the properties of such completions summarized in section
  \ref{sec:completions} apply in this case.
\end{remark}

It follows from corollary \ref{cor:bilateralising-descends} that
the natural morphism $\niv{\D}{m}_{\cX/\cS}\to\niv{\D}{m'}_{\cX/\cS}$
extends uniquely to a morphism
\begin{equation}
  \label{eq:Dmhat-change-of-level}
  \hat\rho_{m',m}:\niv{\hD}{m}_{\cX/\cS}\to \niv{\hD}{m'}_{\cX/\cS}.
\end{equation}
In fact in \ref{eq:Dhatm} we can use a set of $J$ bilateralising for
$\niv{\D}{m'}_{\cX/\cS}$ to compute both $\niv{\D}{m}_{\cX/\cS}$ and
$\niv{\D}{m'}_{\cX/\cS}$. The uniqueness of the extensions shows that
this system of morphisms is transitive.

\begin{thm}\label{thm:coherence-finite-level}
  Suppose $\cX\to\cS$ is a quasi-smooth morphism of adic locally
  noetherian schemes.
  \begin{enumerate}
  \item For any open ideal $J\subset\O_\cX$ bilateralising in
    $\niv{\D}{m}_{\cX/\cS}$, the ring $\niv{\D}{m}_{X_J/\cS}$ is left
    and right coherent.
  \item The ring $\niv{\hD}{m}_{\cX/\cS}$ is left and right coherent.
  \end{enumerate}
\end{thm}
\begin{demo}
  For (i) it suffices, by \cite[Prop. 3.1.3]{berthelot:1996} to show
  that (a) for the canonical injection
  $\O_\cX\to\niv{\D}{m}_{X_J/\cS}$, $\niv{\D}{m}_{X_J/\cS}$ is
  quasicoherent for the $\O_{X_J}$-module structures given by left
  and right multiplication, and (b) for any open affine $U\sset\cX$,
  $\Gamma(U,\niv{\D}{m}_{X_J/\cS}$ is a left and right noetherian. We
  have already seen that (a) is true, and (b) is proven in the same
  way as in \cite[Cor. 2.2.5 (ii)]{berthelot:1996}, i.e. by showing
  the the graded algebra for the filtration by order is finitely
  generated. Since the completion $\niv{\hD}{m}_{\cX/\cS}$ may be
  taken to be the $J$-adic completion for some $J\sset\O_\cX$
  centralising in $\niv{\D}{m}_{\cX/\cS}$, part (ii) follows from the
  facts about centralising ideals recalled in \S\ref{sec:completions},
  c.f.\ also \cite[\S3.2.3]{berthelot:1996} and the last remark in
  that section.
\end{demo}

From \cite[Prop. 3.1.3]{berthelot:1996} we also get:

\begin{prop}\label{prop:coherent-modules-level-m}
  For $\cX\to\cS$ and $J$ as in theorem
  \ref{thm:coherence-finite-level}, a left (resp. right)
  $\niv{\D}{m}_{X_J/\cS}$-module $M$ is coherent if and only if it is
  quasicoherent as an $\O_\cX/J$-module, and for every $U\sset\cX$
  belonging to an open cover of $\cX$, the left (resp. right)
  $\Gamma(U,\niv{\D}{m}_{X_J/\cS})$-module of sections $\Gamma(U,M)$ is
  of finite type.
\end{prop}

The description of coherent left or right
$\niv{\hD}{m}_{\cX/\cS}$-modules is a little more involved, but very
little needs to be added to the treatment of
\cite[\S3.3]{berthelot:1996}. For the reader's convenience we recall
some results regarding completions from \cite[\S3.3]{berthelot:1996},
slightly reformulated for the present purposes. In what follows $\D$
is a sheaf of rings on $\cX$ endowed with a homomorphism
$\O_\cX\to\D$, and we assume that $\D$ satisfies the following
conditions:
\begin{equation}
  \begin{minipage}[t]{0.8\linewidth}
    $\O_\cX$ has an ideal of definition centralising in $\D$; in
    particular $\O_\cX$ has a fundamental system of ideals of
    definition that are bilateralising in $\D$.
  \end{minipage}
\end{equation}
\smallskip
\begin{equation}
  \begin{minipage}[t]{0.8\linewidth}
    For any open affine $U\sset\cX$, $\Gamma(U,\D)$ is left
    noetherian.  \phantom{|}
  \end{minipage}
\end{equation}
\smallskip
\begin{equation}
  \begin{minipage}[t]{0.8\linewidth}
    As a left $\O_\cX$-module, $\D$ is a filtered inductive limit of
    $\O_\cX$-modules $\D_\lambda$ such that for all $\lambda$,
    $\D_\lambda\simeq\liminv_J\D_\lambda/J\D_\lambda$ (where $J$ runs
    through the set of bilateralising ideals of definition), and for
    all $\lambda$ and bilateralising open $J\sset\O_\cX$,
    $\D_\lambda/J\D_\lambda$ is a quasicoherent $O_{X_J}$-module.
  \end{minipage}
\end{equation}

These hypotheses apply in particular to
\begin{displaymath}
  \D=\niv{\D}{m}_{\cX/\cS},\quad
  \D_J=\niv{\D}{m}_{X_J/\cS},\quad
  \hD=\niv{\hD}{m}_{\cX/\cS}=\liminv_J\D_J
\end{displaymath}
and when $\cX$ is affine we write
\begin{displaymath}
  D=\Gamma(\cX,\niv{\D}{m}_{\cX/\cS}),\quad
  D_J=\Gamma(\cX,\niv{\D}{m}_{X_J/\cS}),\quad
  \hat D=\Gamma(\cX,\niv{\hD}{m}_{\cX/\cS})=\liminv_JD_J.
\end{displaymath}
In the inverse limits we can restrict $J$ to run of the powers of a
centralising ideal; then the results cited at the beginning of
\S\ref{sec:completions} show:

\begin{prop}
  Suppose $\D$ satisfies conditions
  \ref{prop:coherent-modules-level-m}.1--3.
  \begin{enumerate}
  \item For any open affine $U\sset\cX$, the ring
    $\Gamma(U,\hat\D)$ is left noetherian.
  \item For any pair of open affines $U'\sset U$, the homomorphism
    \begin{displaymath}
      \Gamma(U,\hat\D)\to\Gamma(U',\hat\D)
    \end{displaymath}
    is right flat.
  \end{enumerate}
\nodemo
\end{prop}

When $\cX$ is affine we denote by $M\mapsto M^\triangle$ the functor
on $D$-modules defined by
\begin{equation}
  \label{eq:associated-sheaf}
  M^\triangle=\liminv_J(M/JM)\tilde{\relax}
\end{equation}
where the tilde denotes sheaf associated to a $O_{X_J}$-module. If we
identify $M$ with the constant presheaf with value $M$, there is a
natural homomorphism $\D\tens_DM\to M^\triangle$. Arguing as in
\cite[3.3.7--8]{berthelot:1996} we obtain:

\begin{prop}\label{prop:triangle-functor}
  With the above hypotheses and notation,
  \begin{enumerate}
  \item The canonical homomorphism $\hD\to\hat D^\triangle$ is an
    isomorphism.
  \item The functor $M\mapsto M^\triangle$ is exact on the category of
    $\hat D$-modules of finite type.
  \item For any $\hat D$-module $M$ of finite type, the canonical
    homomorphism $M\to\Gamma(\cX,M^\triangle)$ is an isomorphism.
  \item For all $\hat D$-modules $M$, $N$ of finite type, the canonical
    homomorphism
    \begin{displaymath}
      \Hom_{\hat D}(M,N)\to\Hom_{\hD}(M^\triangle,N^\triangle)
    \end{displaymath}
    is an isomorphism.
  \end{enumerate}
\end{prop}

As in \cite[3.3.8]{berthelot:1996}, the essential point is to show
that the canonical homomorphism $\Gamma(\cX,\D_J)\to \hat D/J\hat D$
is an isomorphism for all ideals of definition that are bilateralising
in $\D$, and here it is important that $\D_J$ is a quasicoherent
$O_{X_J}$-module. The argument is basically that of \cite[I
10.10.2]{EGA}. Continuing as in \cite[\S3.3]{berthelot:1996} and
\cite[I 10.10]{EGA}, we obtain the following ``theorem A'':

\begin{thm}\label{thm:theorem-A}
  Suppose $\cX$ is affine and $\D$ satsifies conditions
  \ref{prop:coherent-modules-level-m}.1--3. The following are
  equivalent, for any $\hD$-module $\cM$:
  \begin{enumerate}
  \item For every ideal of definition $J\subset\O_\cX$ bilateralising
    in $\D$, the $\D_J$-module $\cM/J\cM$ is coherent, and the
    canonical homomorphism $\cM\to\liminv_J\cM/J\cM$ is an
    isomorphism, where the limit is over ideals of definition
    bilateralising in $\D$.
  \item There is an isomorphism $\cM\isom\liminv_J\cM_J$ where for all
    $J$ as before $\cM_J$ is a coherent $\D_{X_J}$-module, and for
    $J\sset K$, the canonical homomorphism
    $\cM_K\tens_{\O_K}O_{X_J}\to \cM_J$ is an isomorphism.
  \item There is a $\hat D$-module $M$ and an isomorphism $\cM\isom
    M^\triangle$.
  \item The $\hat D$-module $\Gamma(\cX,\cM)$ is of finite type and
    the canonical homomorphism
    $\hD\tens_{\hat D}\Gamma(\cX,\cM)\to\cM$ is an isomorphism.
  \item $\cM$ is a coherent $\hD$-module.    
  \end{enumerate}
\end{thm}

\begin{cor}\label{cor:thmA}
  Suppose $\cX$ is affine. With the above notation, the functors
  \begin{displaymath}
    \cM\mapsto\Gamma(\cX,\cM),
    \qquad
    M\mapsto M^\triangle
  \end{displaymath}
  are inverse equivalences between the category of coherent
  $\hD$-modules and the category of $\hat D$-modules of finite type.
\end{cor}

We also get ``theorem B'':

\begin{thm}\label{thm:thmB}
  Suppose $\cX$ is affine. With the above notation
  \begin{displaymath}
    H^q(\cX,\cM)=0
  \end{displaymath}
  for any coherent $\hD$-module $\cM$.
\end{thm}

For emphasis we repeat that corollary \ref{cor:thmA} and theorem
\ref{thm:thmB} apply to the ring $\niv{\hD}{m}_{\cX/\cS}$ for
quasi-smooth $\cX\to\cS$.

\subsubsection{$m$-PD-stratifications again.}
\label{sec:m-PD-again}

We have seen (proposition
\ref{prop:m-PD-stratifications-and-Dm-module-structure}) that a left
$\niv\D m_{\cX/\cS}$-module structure on an $\O_\cX$-module $M$ is
equivalent to an $m$-PD-stratification on $M$. To get a similar
equivalence for left $\niv\hD m_{\cX/\cS}$-module structures it seems
necessary to put restrictions on $M$.

\begin{defn}\label{defn:pro-qcoh}
  An $\O_\cX$-module is \emph{pro-quasicoherent} if the following
  conditions hold:
  \begin{equation}
    \begin{minipage}[t]{0.7\linewidth}{
        For every ideal of definition $J\sset\O_\cX$; $M/JM$ is a
        quasicoherent sheaf of
        $O_{X_J}$-modules;}\label{enum:module-coeff-ring-B1}
    \end{minipage}
  \end{equation}
  \begin{equation}
    \begin{minipage}[t]{0.7\linewidth}
      {$M\simeq\liminv_JM/JM$ where the inverse limit is over
        ideals of definition of $\cX$.\label{enum:module-coeff-ring-B2}}
    \end{minipage}
  \end{equation}
\end{defn}
We do not claim that this is a particularly useful class of
$\O_\cX$-modules. As a full subcategory of the category of left
$\O_\cX$-modules, it is additive and kernels of morphisms are
representable, but not necessarily cokernels, at least not with
further restrictions. It does guarantee that $M$ can be recovered from
the $M/JM$ for all $J$, and that local calculations are
possible. Furthermore there is a large supply of such modules:

\begin{prop}\label{prop:coh-hD-implies-pro-qcoh}
  A coherent left $\niv\hD m_{\cX/\cS}$-module is pro-quasicoherent as
  an $\O_\cX$-module.
\end{prop}
\begin{demo}
  This follows immediately from theorem \ref{thm:theorem-A} applied
  with $\D=\niv\D m_{\cX/\cS}$, which is applicable by the discussion
  after proposition \ref{prop:coherent-modules-level-m}. 
\end{demo}

In what follows we will say ``pro-quasicoherent $\niv\hD
m_{\cX/\cS}$-module'' to mean ``$\niv\hD m_{\cX/\cS}$-module that is
pro-quasicoherent as an $\O_\cX$-module.'' In other words
``pro-quasicoherent'' refers to the $\O_\cX$-module structure, and not
to any other.

\begin{prop}\label{prop:m-PD-stratifications-and-Dmhat}
  The category of pro-quasicoherent left
  $\niv\hD m_{\cX/\cS}$-modules is equivalent to the category of
  pro-quasicoherent $\O_\cX$-modules endowed with an
  $m$-PD-stratification.
\end{prop}
\begin{demo}
  Suppose $M$ is pro-quasicoherent. It is clear that a left
  $\niv\hD m_{\cX/\cS}$-module structure on $M$ can be identified with
  a set of left $\niv\D m_{\cX/\cS}$-module structures on the $M/JM$
  as $J$ runs through the set of $m$-bilateralising ideals of
  definition, where we impose that these left module structures are
  compatible with the quotient maps $M/J'M\to M/JM$ for $J'\sset
  J$. This is equivalent to a set of $m$-PD-stratifications on each of
  the $M/JM$, again with the same compatibility condition. The limit
  of these $m$-PD-stratifications is an isomorphism
  \begin{displaymath}
    \chi_n:\cP^n_{\cX/\cS,(m)}\ctens_{\O_\cX}M\isom
    M\ctens_{\O_\cX}\cP^n_{\cX/\cS,(m)}
  \end{displaymath}
  satisfying the cocycle condition. The pullback of this isomorphism
  by the diagonal is the identity of $M$, as we see from condition
  \ref{enum:module-coeff-ring-B2} and the fact that it is the identity
  modulo all $J$. Finally, we can remove that hat on the tensor
  product since $\cP^n_{\cX/\cS,(m)}$ is a locally free
  $\O_\cX$-module of finite type. This construction is clearly
  reversible, whence the equivalence.
\end{demo}

\subsection{The ring $\niv{\hD}{m}_{\cX/\cS\bQ}$.}
\label{sec:Dm-Q}

\subsubsection{} Note that corollary \ref{cor:thmA} and theorem
\ref{thm:thmB} also apply to the ring
\begin{equation}
  \label{eq:Dm-Q}
  \niv{\hD}{m}_{\cX/\cS\bQ}:=\niv{\hD}{m}_{\cX/\cS}\tens\bQ
\end{equation}
In fact if $\cX$ is affine and $M$ is a coherent
$\niv{\hD}{m}_{\cX/\cS\bQ}$-module there is a coherent
$\niv{\hD}{m}_{\cX/\cS}$-module $M_0$ such that $M\simeq M_0\tens\bQ$
(we will later prove a more general result of this sort, namely
proposition \ref{prop:Ogus-integrality-thm}). The analogue of
corollary \ref{cor:thmA} is an immediate consequence of this
observation; for theorem \ref{thm:thmB} it suffices to observe that
since $\cX$ is noetherian, cohomology commutes with inductive limits.

The most important result concerning these rings is:

\begin{thm}\label{thm:flatness-change-of-m}
  For all $m\le m'$, the canonical homomorphism
  \begin{displaymath}
    \niv{\hD}{m}_{\cX/\cS\bQ}\to\niv{\hD}{m'}_{\cX/\cS\bQ}
  \end{displaymath}
  is left and right flat.
\end{thm}

By induction one can assume $m'=m+1$. We will give the proof as a
series of lemmas; as the argument is roughly parallel to that of
\cite[Th. 3.5.3]{berthelot:1996} we will concentrate the few
modifications that are needed. It begins with reductions to the case
where $\cS=\Spf{R}$ is affine and $\cX=\Spf{A}$ is affine and
parallelizable relative to $\cS$. The next reduction is to the case
where $\O_\cS$ is $p$-torsion free. If $\cS'\inj\cS$ is the closed
formal subscheme defined by the ideal of $p$-torsion elements of $R$
and $\cX'=\cX\times_\cS\cS'$, then $\cX'\to\cS'$ is quasi-smooth and
therefore flat. Then $\cX'$ is flat over $\bZ_p$ since this is the
case for $\cS'$. Finally if we write $\cS'=\Spf{R'}$ and
$\cX'=\Spf{A'}$ then the natural morphism
$\niv{\hat D}{m}_{A/R}\to\niv{\hat D}{m}_{A'/R'}$ induces an
isomorphism $\niv{\hat D}{m}_{A/R\bQ}\to\niv{\hat D}{m}_{A'/R'\bQ}$.

In what follows we denote the $p$-adic completion of a $\bZ_p$-algebra
$A$ by $\tilde{A}$; the $J$-adic completion of an $\O_\cX$-ring $A$
will be denoted by $\hat A$ as before (we assume $J$ is bilateralising
in $A$). We abbreviate $\niv{D}{m}=\niv{D}{m}_{A/R}$ and similarly for
$\niv{\tilde D}{m}$, $\niv{\hat D}{m}$ etc. Note that since $\O_\cX$
is $\bZ_p$-flat and $\niv{D}{m}$ is free on the $\dpabniv{\d}{K}{m}$,
the maps $\niv{\tilde D}{m}\to\niv{\hat D}{m}$ and
$\niv{\tilde D}{m}_\bQ\to\niv{\hat D}{m}_\bQ$ are injective. We define
\begin{align*}
  D_1&=\text{the subring of $\niv{\hat D}{m+1}_\bQ$ generated by
       $\niv{\tilde D}{m}$ and $\niv{D}{m+1}$},\\
  D_2&=\text{the subring of $\niv{\hat D}{m+1}_\bQ$ generated by
       $\niv{\hat D}{m}$ and $\niv{D}{m+1}$}\\
\end{align*}
and note that
\begin{equation}
  \label{eq:description-of-D1,D2}
  \begin{minipage}{0.75\linewidth}
    $D_1$ (resp. $D_2$) is generated as a left
    $\niv{\tilde D}{m}$-module (resp. left $\niv{\hat D}{m}$-module) by
    the $(\dpe{\d}{p^{m+1}})^K$ for all $K\in\bN^d$.
  \end{minipage}
\end{equation}
%
%
From this it follows that
\begin{equation}
  \label{eq:D(m+1)fcapD2}
  \niv{D}{m+1}\cap JD_2\sset J\niv{D}{m+1}.
\end{equation}
We now consider the injective homomorphisms
\begin{equation}
  \label{eq:flatness-inclusion}
  \niv{D}{m+1}\inj D_1\inj D_2.
\end{equation}

\begin{lemma}\label{lemma:flatness-1st-isom}
  (i) For any $n\ge0$,
  \begin{displaymath}
    D_1=\niv{D}{m+1}+p^n\niv{\tilde D}{m}.
  \end{displaymath}
  (ii) The first map in \ref{eq:flatness-inclusion} induces an
  isomorphism
  \begin{displaymath}
    \niv{\hat D}{m+1}\isom\hat D_1.
  \end{displaymath}
\end{lemma}
\begin{demo}
  The argument of \cite[Th. 3.5.3]{berthelot:1996} can be used without
  change to prove (i) in the case $n=0$ and to deduce from this that
  $\niv{D}{m+1}\inj D_1$ induces an isomorphism
  $\niv{\tilde D}{m+1}\isom\tilde D_1$ (note that one needs $A$ to be
  $\bZ_p$-flat to show that the map
  $\niv{D}{m+1}/p^i\niv{D}{m+1}/\to D_1/p^iD_1$ is injective). Since
  $\niv{\tilde D}{m}=\niv{D}{m}+p^n\niv{\tilde D}{m}$, the general
  case of (i) follows from the case $n=0$.

  Any open $(m+1)$-bilateralising ideal $J$ contains a power of $p$, so
  the $J$-adic completions of $\niv{\tilde D}{m+1}$ and $\tilde D_1$
  may be identified with $\niv{\hat D}{m+1}$ and $\hat D_1$
  respectively. The isomorphism $\niv{\tilde D}{m+1}\isom\tilde D_1$
  thus induces an isomorphism $\niv{\hat D}{m+1}\isom\hat D_1$.
\end{demo}

Note that the inclusion $\niv{D}{m+1}+J^n\niv{\hat D}{m}\subset D_2$
is strict, which is what necessitates the introduction of the ring
$D_1$. 

\begin{lemma}\label{lemma:flatness-2nd-isom}
  The second map in \ref{eq:flatness-inclusion} induces an isomorphism
  \begin{displaymath}
    \hat D_1\isom\hat D_2.
  \end{displaymath}
\end{lemma}
\begin{demo}
  It suffices to show that form any open $(m+1)$-bilateralising ideal
  $J\subset A$, the induced map
  \begin{displaymath}
    D_1/JD_1\to D_2/JD_2
  \end{displaymath}
  is a bijection. That it is surjective follows from
  \ref{eq:description-of-D1,D2}. By (i) of lemma
  \ref{lemma:flatness-1st-isom}, injectivity is equivalent to
  \begin{displaymath}
    (\niv{D}{m+1}+\niv{\tilde D}{m})\cap J\D_2\sset
    J\niv{D}{m+1}+J\niv{\tilde D}{m}.
  \end{displaymath}
  Suppose $P\in (\niv{\tilde D}{m}+\niv{D}{m+1})\cap J\D_2$; again by
  (i) in lemma \ref{lemma:flatness-1st-isom} we can write $P=Q+p^nR$
  with $Q\in\niv{D}{m+1}$ and $R\in\niv{\tilde D}{m}$ for any given
  $n\ge0$. If we choose $n$ so that $p^n\in J$,
  \begin{displaymath}
    P-Q=p^nR\in J\niv{\tilde D}{m}\sset JD_2
  \end{displaymath}
  and since $P\in JD_2$ this implies $Q\in\niv{D}{m+1}\cap JD_2$. Then
  $Q\in J\niv{D}{m+1}$ by \ref{eq:D(m+1)fcapD2} and consequently
  $P=Q+p^nR\in J\niv{D}{m+1}+J\niv{\tilde D}{m}$, as required.
\end{demo}

The last bit of the proof of \cite[Th. 3.5.3]{berthelot:1996} can be
used without change to prove:

\begin{lemma}\label{lemma:flatness-D2-noetherian}
  The ring $D_2$ is left noetherian.
\end{lemma}

The proof of theorem \ref{thm:flatness-change-of-m} is concluded as
follows: since $D_2$ is noetherian, $\hat D_2$ is a left flat
$D_2$-module, and thus $\hat D_{2\bQ}$ is a left flat
$D_{2\bQ}$-module. Now lemmas \ref{lemma:flatness-1st-isom} and
\ref{lemma:flatness-2nd-isom} show that
$\hat D_2\simeq\niv{\hat D}{m+1}$ and
$\hat D_{2\bQ}\simeq\niv{\hat D}{m+1}_\bQ$. On the other hand
\ref{eq:description-of-D1,D2} yields $D_{2\bQ}=\niv{\hat D}{m}_\bQ$,
and it follows that $\niv{\hat D}{m+1}_\bQ$ is a left flat
{$\niv{\hat D}{m}_\bQ$-module. The assertion for any $m'\ge m$ follows
by induction.\nodemo}

\subsection{The ring $\Ddag_{\cX/\cS\bQ}$.}
\label{sec:Ddag}

As in \cite[\S2.5]{berthelot:1996} we define
\begin{equation}
  \label{eq:Ddag}
  \Ddag_{\cX/\cS\bQ}=\limdir_m\niv{\hD}{m}_{\cX/\cS\bQ}
\end{equation}
with the inductive limit over the canonical morphisms
\ref{eq:Dmhat-change-of-level}. When $\cX=\Spf{A}$ and $\cS=\Spf{R}$
are formally affine we put
\begin{displaymath}
  \Ddag_{A/R}=\Gamma(\cX,\Ddag_{\cX/\cS}).
\end{displaymath}
It follows from theorem \ref{thm:flatness-change-of-m} by taking
inductive limits that the canonical inclusion
\begin{displaymath}
  \niv{\hD}{m}_{\cX/\cS\bQ}\to\Ddag_{\cX/\cS\bQ}
\end{displaymath}
is also flat.

\begin{thm}\label{thm:coherence-of-Ddag}
  If $\cX/\cS$ is quasi-smooth, $\Ddag_{\cX/\cS\bQ}$ is a left and
  right coherent sheaf of rings.\nodemo
\end{thm}

As in \cite[\S3.5]{berthelot:1996}, this follows from theorem
\ref{thm:flatness-change-of-m} by standard arguments. Also standard is
the deduction of theorems A and B from the corresponding theorems for
the rings $\niv{\hD}{m}_{\cX/\cS}$:

\begin{thm}\label{thm:thmA-Ddag}
  Suppose $\cX=\Spf{A}$ and $\cS=\Spf R$ are noetherian and formally
  affine $\cX\to\cS$ is quasi-smooth. The functors
  \begin{displaymath}
    \cM\mapsto\Gamma(\cX,\cM),
    \qquad
    M\mapsto M^\triangle
  \end{displaymath}
  are inverse equivalences between the category of coherent
  $\Ddag_{\cX/\cS}$-modules and the category of
  $D^\dagger_{A/R}$-modules of finite type.\nodemo
\end{thm}

\begin{thm}\label{thm:thmB-Ddag}
  Suppose $\cX\to\cS$ is quasi-smooth and $\cX$ is affine. With the
  above notation
  \begin{displaymath}
    H^q(\cX,\cM)=0
  \end{displaymath}
  for any coherent $\Ddag_{\cX/\cS}$-module $\cM$.\nodemo
\end{thm}

\section{Stratifications}
\label{sec:stratifications}

We can now return to a question that was left open in section
\ref{sec:m-PD-formal-schemes}, that of the sheafification of the full
$m$-PD-envelope of a regular ideal. As in the classical theory, these
are used to characterize topologically quasi-nilpotent left $\niv\hD
m$-modules. 

\subsection{The ring $\cP_{(m)}(I)$.}
\label{sec:P_m(I)}

As always we begin with the affine case, so let $R$ be an adic
noetherian $\bZ_p$-algebra with $m$-PD-structure $(\fa,\fb,\alpha)$
and $A$ an adic noetherian $R$-algebra. By hypothesis any ideal of
definition of $A$ contains a power of $p$. For $n\ge0$ we set
$R_n=R/p^{n+1}R$ and $A_n=A/p^{n+1}A$.
\begin{defn}\label{defn:split-regular-ideal}
  An ideal $I\subset A$ is \emph{split-regular} if
  \begin{enumerate}
  \item the ideal $IA_0\subset A_0$ is regular;
  \item the quotient homomorphism $\pi:A\to A/I$ has a section
    $\sigma:A/I\to A$;
  \item $A/I$ is a flat $R$-algebra.
  \end{enumerate}
\end{defn}
In what follows it will be convenient to set $A'=A/I$.  If
$(\bar f_1,\ldots,\bar f_r)$ is a regular sequence generating $IA_0$
we pick $f_1,\ldots,f_r$ in $I$ lifting $\bar f_1,\ldots,\bar
f_r$. Then $I=(f_1,\ldots,f_r)$ and $(p,f_1,\ldots,f_r)$ are regular,
as is $(p^{n+1},f_1,\ldots,f_r)$ for all $n\ge0$. It follows that
$IA_n$ is regular for all $n\ge0$.

The $m$-PD-structure $(\fa,\fb,\alpha)$ descends to $R_n$ and we
denote by $P_{(m)}(IA_n)$ the $(\fa,\fb,\alpha)$-compatible
$m$-PD-envelope of $IA_n\subset A_n$ (we will never deal with the full
$m$-PD-envelope $P_{(m)}(I)$ of $I\subset A$).  If we set
$A'_n=A'\tens_AA_n$ then $P_{(m)}(IA_n)$ with the $A'_n$-module
structure determined by $\sigma$ is a free $A'_n$-module, with basis
the $m$-PD-monomials in the regular generators of $IA_n$. We note, for
use in the arguments that follow that $p$ is nilpotent in $A_n$, so
that $P_{(m)}(IA_n)$ has the usual base-change properties.

For any $n'\ge n$ there is an isomorphism
\begin{displaymath}
  P_{(m)}(IA_{n'})\tens_{R_{n'}}R_n\isom P_{(m)}(IA_n)
\end{displaymath}
since the formation of $P_{(m)}(IA_n)$ is compatible with arbitrary
base-change in $R_n$ (c.f. the remark in the last paragraph). Since
this merely says that
\begin{displaymath}
  P_{(m)}(IA_n)\simeq P_{(m)}(IA_{n'})/p^nP_{(m)}(IA_{n'})
\end{displaymath}
we may rewrite it as an isomorphism
\begin{equation}
  \label{eq:P_m-base-change}
  P_{(m)}(IA_{n'})\tens_{A_{n'}}A_n\isom P_{(m)}(IA_n).
\end{equation}
Let $J\subset A$ be an ideal of definition and choose $n\ge0$ such
that $p^{n+1}\in J$. Then $A/J$ is a $A_n$-algebra and the isomorphism
\ref{eq:P_m-base-change} shows that
\begin{equation}
  \label{eq:P_A/R,J}
  P_{J,(m)}(I):=P_{(m)}(IA_n)\tens_{A_n}(A/J)
  \simeq P_{(m)}(IA_n)/JP_{(m)}(IA_n)
\end{equation}
is independent of the choice of $n$ (and justifies the notation).
For $m'\ge m$ there is a natural morphism
\begin{equation}
  \label{eq:PmJ-change-of-level}
  P_{J,(m')}(I)\to P_{J,(m)}(I)
\end{equation}
arising from the fact that an $m$-PD-structure is automatically an
$m'$-PD-structure. Finally, we set
\begin{equation}
  \label{eq:P_A/R-hat}
  \hat P_{(m)}(I)=\liminv_J P_{J,(m)}(I)
\end{equation}
where $J$ runs through the set of ideals of definition of $A$. We
denote by 
\begin{equation}
  \label{eq:completion-of-canonical-m-PD-structure}
  \hat I^\canj\sset\hat I^\cani\subset\hat P_{(m)}(I)
\end{equation}
the closures of $I^\canj$ and $I^\cani$ in $\hat P_{(m)}(I)$ (or
equivalently, the $J$-adic completions). The change-of-level morphisms
\ref{eq:PmJ-change-of-level} induce ring homomorphisms
\begin{equation}
  \label{eq:hatPmJ-change-of-level}
  \hat P_{J,(m')}(I)\to \hat P_{J,(m)}(I)
\end{equation}
for all $m'\ge m$. By construction the change-of-level formula
\ref{eq:Pm-change-of-level-explicit} holds for any $x\in I$.

\begin{lemma}\label{lemma:flatness-of-hatP}
  Suppose $I$ is split-regular. If $R$ is $\bZ_p$-flat, so is
  $\hat P_{(m)}(I)$.
\end{lemma}
\begin{demo}
  We know that for all $n\ne0$, $P_{(m)}(IA_n)$ is flat over
  $R_n$. Tensoring the exact sequence
  \begin{displaymath}
    0\to p^nR/p^{n+1}R\to R/p^{n+1}R\Xto{p} R/p^{n+1}R
  \end{displaymath}
  with $P_{(m)}(IA_n)$ yields an exact sequence
  \begin{displaymath}
    0\to P_{(m)}(IA_n)\tens_{R_n}p^nR/p^{n+1}R
    \to P_{(m)}(IA_n)\Xto{p}P_{(m)}(IA_n).
  \end{displaymath}
  Since the inverse system $\{p^nR/p^{n+1}R\}$ is essentially null,
  passing to the inverse limit in $n$ shows that multiplication by $p$
  is injective.
\end{demo}

What is not obvious in this construction is whether the canonical
$m$-PD-structure of $P_{(m)}(IA_n)$ descends to $P_{J,(m)}(I)$ when
$p^{n+1}\in J$, or extends to the completion $\hat P_{(m)}(I)$, and if
so, whether the extensions are compatible with $(\fa,\fb,\alpha)$. We
will restrict our attention to the case where $J$ satisfies
\begin{equation}
  \label{eq:special-J1}
  \begin{minipage}[t]{0.7\linewidth}
    $JP_{(m)}(IA_n)=\sigma(J')P_{(m)}(IA_n)$ for some ideal $J'\subset
    A'=A/I$ and some $n$ such that $p^{n+1}\in J'$.
  \end{minipage}
\end{equation}
When $p$ is nilpotent in $A$, $P_{(m)}(I)/I^\cani\simeq A/I\simeq A'$,
and thus \ref{eq:special-J1} implies that $J'=p(J)$ where
$\pi:A\to A'$ is the canonical projection, and we see that
\ref{eq:special-J1} is equivalent to
\begin{equation}
  \label{eq:special-J2}
  \begin{minipage}[t]{0.7\linewidth}
    $\sigma(\pi(J))P_{(m)}(IA_n)=JP_{(m)}(IA_n)$ for some $n$ such
    that $p^{n+1}\in J$.
  \end{minipage}
\end{equation}
For any given $n$, $p^{n+1}\in J$ if and only if $p^{n+1}\in J'$; it
follows that the conditions \ref{eq:special-J1}, \ref{eq:special-J2}
are independent of the particular value of $n$.  
\begin{defn}\label{defn:adapted-to-split-regular-ideal}
  Suppose $I\sset A$ is a split-regular ideal. An ideal $J\sset A$ is
  \emph{adapted to $I$} if the equivalent conditions
  \ref{eq:special-J1}, \ref{eq:special-J2} hold.
\end{defn}
The following is immediate from either \ref{eq:special-J1} or
\ref{eq:special-J2}:
\begin{lemma}\label{lemma:power-of-adapted-is-adapted}
  Suppose $I\sset A$ is split-regular. If $J$ is adapted to $I$ then
  so is $J^n$ for any $n$.
\end{lemma}

When $I$ is generated by a regular sequence $(x_1,\ldots,x_d)$ and
$J\subset A'$ satisfies \ref{eq:special-J1}, $P_{J,(m)}(I)$ is the
free $A'/J$-module on the $m$-PD-monomials $\dpbrniv{x}{K}{m}$; this
follows from the description of the full $m$-PD-envelope in section
\ref{sec:m-PD-regular}.

If $A$ has an ideal of definition adapted to $I$ we can give a similar
description of the completion $\hat P_{(m)}(I)$: elements of
$\hat P_{(m)}(I)$ can be identified with series
$\sum_Ka_K\dpbrniv{x}{K}{m}$ with $a_K\in A'$ and $a_K\to0$ as
$|K|\to\infty$. In fact if $J\subset A$ is an ideal of definition
adapted to $I$ then $J^n$ satsifies \ref{eq:special-J1} as well, and
furthermore $J'=p(J)$ is an ideal of definition of $A'$. The above
description of $\hat P_{(m)}(I)$ follows from the previous description
of the $\hat P_{J^n,(m)}(I)$. In particular
$\hat P_{(m)}(I)/\hat I^\cani\simeq A'$ is a flat $R$-algebra.

\begin{prop}\label{prop:descent-of-m-PD-str-to-P/J}
  Let $R$ be a ring with $m$-PD-structure $(\fa,\fb,\alpha)$, $A$ an
  $R$-algebra. Suppose $I\subset A$ is a split-regular ideal, and
  $J\subset A$ is adapted to $I$.
  \begin{enumerate}
  \item For any $n$ such that $p^{n+1}\in J$, the canonical
    $m$-PD-structure of the $m$-PD-envelope $P_{(m)}(A_nI)$ descends
    to $P_{J,(m)}(I)$.
  \item If $\fb_1=\fb+pR$, the $m$-PD-structure on $P_{J,(m)}(I)$ is
    compatible with $(\fa,\fb,\alpha)$ if and only if $J'\cap \fb_1A'$
    is a sub-PD-ideal of $\fb_1A'$.
  \end{enumerate}

\end{prop}
\begin{demo}
  If $p^{n+1}\in J$ we can replace $A$ and $I$ by $A_n$ and $IA_n$,
  which is a regular ideal in $A_n$; furthermore $A_n/IA_n$ is a flat
  $R_n$-algebra and the section $\sigma$ induces a section of
  $A_n\to A_n/IA_n$. We may therefore assume that $p$ is nilpotent in
  $A$ and set $P=P_{(m)}(I)$; from lemma
  \ref{lemma:quotient-m-PD-structure} (c.f. also
  \cite[1.3.4]{berthelot:1996}) we see that the conditions to be
  checked are that the following are sub-PD-ideals:
  \begin{align}
    \label{eq:descent-of-m-PD-to-P/J1}
    (I^\canj+pP)\cap J&\sset I^\canj+pP\\
    \label{eq:descent-of-m-PD-to-P/J2}
    (I^\canj+\fb_1P)\cap J&\sset I^\canj+\fb_1P\\
    \label{eq:descent-of-m-PD-to-P/J3}
    \fb_1P+(I^\cani+J)&\sset\fb_1P
  \end{align}
  where \ref{eq:descent-of-m-PD-to-P/J1} guarantees that the
  $m$-PD-structure of $P$ descends to $P/J$, and
  \ref{eq:descent-of-m-PD-to-P/J2}, \ref{eq:descent-of-m-PD-to-P/J3}
  guarantee that it is compatible with $(\fa,\fb,\alpha)$.  The
  question is Zariski-local so we may assume that $I$ is generated by
  a regular sequence $(x_1,\ldots,x_n)$. Then $P$ is a free
  $A'$-module on the $m$-PD-monomials $\dpbrniv{x}{K}{m}$ and the
  ideals in
  \ref{eq:descent-of-m-PD-to-P/J1}--\ref{eq:descent-of-m-PD-to-P/J3}
  have the following descriptions, where
  $x=\sum_Ka_K\dpbrniv{x}{K}{m}$ and $K<p^m$ means that every entry of
  $K$ is less than $p^m$:
  \begin{align*}
    x\in(I^\canj+pP)\cap J&\iff a_K\in J',\text{\ and\ }
                             K<p^m\implies a_K\in pA'\\
    x\in(I^\canj+\fb_1P)\cap J&\iff a_K\in J',\text{\ and\ }
                             K<p^m\implies a_K\in \fb_1A'\\
    x\in\fb_1P\cap(I^\cani+J)&\iff a_0\in J',\text{\ and\ }a_K\in \fb_1A'.
  \end{align*}
  Thus \ref{eq:descent-of-m-PD-to-P/J1} is a sub-PD-ideal because
  $pA'\cap J'$ is a sub-PD-ideal of $pA'$. If
  \ref{eq:descent-of-m-PD-to-P/J3} is a sub-PD-ideal then
  $\fb_1A'\cap J'$ is a sub-PD-ideal of $\fb_1A'$, and conversely this
  implies that \ref{eq:descent-of-m-PD-to-P/J2} and
  \ref{eq:descent-of-m-PD-to-P/J3} are sub-PD-ideals.
\end{demo}

With the hypotheses of the lemma, $(\fa,\fb,\alpha)$ extends to $A'$,
and the condition in (ii) is equivalent to the assertion that the
$m$-PD-structure $(\fa A',\fb A',\alpha)$ descends to $A'/J'$. Note
that this is automatic if $\fb_1$ is principal, or if $A'/J'$ is a
flat $R$-algebra. 

If $J$ is adapted to $I$, $P_{J,(m)}(I)$ has an $m$-PD-adic
filtration, and we denote by $P^n_{J,(m)}(I)$ the quotient of
$P^n_{J,(m)}(I)$ by the $n+1$-st step of that filtration.  On the
other hand the truncations $P^n_{(m)}(I)$ of the full $m$-PD-envelope
of the diagonal commutes with arbitrary base change in $R$, and in
particular with the base change $R\to R_n$. From the construction we
see that $P^n_{J,(m)}(I)$ is iso\-morphic to $P^n_{J,(m)}(I)/JP^n_{J,(m)}(I)$.

\begin{cor}\label{cor:descent-of-m-PD-str-to-P/J^n}
  With the assumptions of proposition
  \ref{prop:descent-of-m-PD-str-to-P/J}, for any $J\subset A$ adapted
  to $I$ the canonical $m$-PD-structure of $P_{J^n,(m)}(I)$ is
  compatible with the $m$-PD-structure $(\fa,\fb,\alpha)$ of $R$ for
  all sufficiently large $n$.
\end{cor}
\begin{demo}
  With the notation of the proposition and its proof, it suffices to
  show that $(J')^n\cap\fb_1A'$ is a sub-PD-ideal of $\fb_1A'$ for
  $n\gg0$, but this follows from lemma
  \ref{lemma:passing-m-PD-to-quotient}. 
\end{demo}

\begin{thm}\label{thm:descent-of-m-PD-str-to-Phat}
  Suppose $I\subset A$ is a split-regular ideal and $A$ has an ideal
  of definition adpated to $I$. The canonical $m$-PD-structure of
  $P_{(m)}(I)$ extends to an $m$-PD-structure
  $(\hat I^\cani,\hat I^\canj,\hat\gamma)$ on $\hat P_{(m)}(I)$ with
  $\hat I^\cani$ and $\hat I^\canj$ as in
  \ref{eq:completion-of-canonical-m-PD-structure}, and this
  $m$-PD-structure is compatible with $(\fa,\fb,\alpha)$.
\end{thm}
\begin{demo}
  Let $J$ be an ideal of definition adapted to $I$ and denote by
  $(I^\cani_n,I^\canj_n,\gamma_n)$ the quotient $m$-PD-structure of
  $P_{J^n,(m)}(I)$. By construction
  \begin{displaymath}
    \hat I^\cani=\liminv_n I^\cani_n,\qquad \hat I^\canj=\liminv_n I^\canj
  \end{displaymath}
  and the containments
  $\niv{(I_n^\cani)}{p^m}+p I_n^\cani\sset I_n^\canj$ for all $n$ show
  that $\niv{(\hat I^\cani)}{p^m}+p\hat I^\cani\sset\hat I^\canj$.  If
  $\gamma_n=\{\gamma_{n,k}\}_{k>0}$, the functions
  $\hat\gamma_k=\liminv_n\gamma_{n,k}$ define a PD-structure on
  $\hat I^\canj$, and $(\hat I^\cani,\hat I^\canj,\hat\gamma)$ an
  $m$-PD-structure on $\hat P_{(m)}(I)$. We must show that
  $(\hat I^\cani,\hat I^\canj,\hat\gamma)$ is compatible with
  $(\fa,\fb,\alpha)$; this means that
  $\fb_1\hat P_{(m)}(I)+\hat I^\canj$ has a PD-structure extending the
  PD-structures $\bar\alpha$ of $\fb_1$ and $\hat\gamma$ of
  $\hat I^\canj$, and that $\fb_1\hat P_{(m)}(I)\cap\hat I^\cani$ is a
  sub-PD-ideal of $\fb_1\hat P_{(m)}(I)$. By construction
  $\fb_1P_{(m)}(I)+I^\canj$ has a PD-structure $\{\delta_k\}_{k>0}$
  extending $\bar\alpha$ and $\gamma$; on the other hand corollary
  \ref{cor:descent-of-m-PD-str-to-P/J^n} says that
  $(I^\cani_n,I^\canj_n,\gamma_n)$ of is compatible with
  $(\fa,\fb,\alpha)$ for all $n\gg0$, which implies that for all
  $k>0$, $\delta_k$ is $J$-adically continuous on
  $\fb_1P_{(m)}(I)+I^\canj$. Therefore $\delta$ extends by continuity
  to the closure of $\fb_1P_{(m)}(I)+I^\canj$, and in particular to
  $\fb_1\hat P_{(m)}(I)+\hat I^\canj$. Finally we observed earlier
  that $\hat P_{(m)}(I)/\hat I^\cani\simeq A'$ is a flat $R$-algebra,
  which implies that
  $\fb_1\hat P_{(m)}(I)\cap\hat I^\cani=\fb_1\hat I^\cani$ is a
  sub-PD-ideal of $\fb_1\hat P_{(m)}(I)$.
\end{demo}

Passing to the limit over $J$ in \ref{eq:hatPmJ-change-of-level}
yields ring homomorphisms
\begin{equation}
  \label{eq:eq:hatPm-change-of-level}
  \hat P_{(m')}(I)\to \hat P_{(m)}(I)
\end{equation}
for all $m'\ge m$. By construction the change-of-level formula
\ref{eq:Pm-change-of-level-explicit} holds for any $x\in I$.

The universal properties of these rings are as follows:

\begin{prop}\label{prop:univ-property-of-m-PD-envelope}
  Suppose that $A$ has an ideal of definition adapted to $I$.
  \begin{enumerate}
  \item Let $A'$ be a \textit{discrete} topological $R$-algebra
    with an $m$-PD-structure $(I',J',\gamma)$ compatible with
    $(\fa,\fb,\alpha)$, and suppose $f:A\to A'$ is a continuous
    $R$-algebra homomorphism such that $f(I)\sset I'$. For any
    $m$-bilateralising ideal of definition $K\subset A$ such that
    $f(K)=0$, $f$ has a unique factorisation
    \begin{displaymath}
      A\to P_{K,(m)}(I)\Xto{f_K}A'
    \end{displaymath}
    in which $f_K$ is an $m$-PD-homomorphism over $R$.
  \item Suppose $A'$ is an adic noetherian $R$-algebra with an
    $m$-PD-structure $(I',J',\gamma)$ compatible with
    $(\fa,\fb,\alpha)$. Any continuous $R$-algebra homomorphism
    $f:A\to A'$ such that $f(I)\sset I'$ has a unique factorization
    \begin{displaymath}
      A\to\hat P_{(m)}(I)\Xto{g}A'
    \end{displaymath}
    in which $g$ is an some $m$-PD-homomorphism $g$ over $R$.
  \end{enumerate}
\end{prop}
\begin{demo}
  (i) Pick $n$ such that $p^{n+1}\in J$; then $f$ factors through a
  morphism $f_n:A_n\to A'$, and the $m$-PD-structure $(I',J',\gamma)$
  descends to an $m$-PD-structure $(I'A_n,J'A_n,\bar\gamma)$
  compatible with $(\fa,\fb,\alpha)$. It follows that $f_n$ factors
  through a unique $m$-PD-morphism $f':P_{(m),\alpha}(IA_n)\to A'$,
  and since $JA'=0$, $f'$ factors through an $m$-PD-morphism
  $f':P_{(m),J,\alpha}(I)\to A'$ which is unique since $f'$ is.

  (ii) The same argument as before shows that for all $n\ge0$ the
  reduction $f_n:A_n\to A'_n$ of $f$ factors uniquely through an
  $m$-PD-morphism $P_{(m),\alpha}(IA_n)\to A'_n$ for all $n\ge0$. By
  lemma \ref{lemma:passing-m-PD-to-quotient} we know that $A'$ has an
  ideal of definition $K'$ such that the the $m$-PD-structure of $A'$
  descends to $A/(K')^n$ for all $n$ (c.f. the remark after the
  lemma). We can then find an ideal of definition $K\subset A$ such
  that $f(K)\sset K'$ and the $m$-PD-structure of $A$ descends to
  $A/K^n$ for all $n$. For any particular $n$ we can choose an $n'$
  such that $p^{n'+1}\in K^n$; then the morphism
  $f_{n'}\to P_{(m)}(IA_{n'})\to A'_{n'}$ induces a morphism
  $g_n:P_{(m),K^n}(I)\to A'/(K')^n$. The latter morphism is
  necessarily an $m$-PD-morphism since $f_{n'}$ is, and since
  $P_{(m)}(IA_{n'})\to P_{(m),K^n}(IA_{n'})$ is an
  $m$-PD-morphism. Since $A'$ is $K'$-adically complete, the inverse
  limit of the $g_n$ is an $m$-PD-morphism
  $\hat g:\hat P_{(m)}(I)\to A'$, and the construction shows that it
  is the unique morphism that factors $f$.
\end{demo}

We now globalize these constructions. Let $\cX\to\cS$ be a universally
noetherian morphism of locally noetherian adic formal $\bZ_p$-schemes
and $(\fa,\fb,\alpha)$ is an $m$-PD-structure on $\cS$. 
\begin{defn}\label{defn:split-regular-sheaf}
  An ideal $\cI\sset\O_\cX$ is \emph{split-regular} if
  \begin{enumerate}
  \item The image of $\cI$ in $\O_\cX\tens_{\bZ_p}\bF_p$ is regular
    ideal;
  \item If $\cY=V(\cI)$, the canonical closed immersion $\cY\to\cX$
    locally has a retraction $\cX\to\cY$ over $\cS$.
  \item $\cY\to\cS$ is flat.
  \end{enumerate}
  If $\cI\sset\O_\cX$ is split-regular, and ideal $\cJ\sset\O_\cX$ is
  \emph{adapted to $\cI$} if it locally satisfies the equivalent
  conditions \ref{eq:special-J1}, \ref{eq:special-J2}, in which
  $\sigma$ is the homomorphism associated to the (local) retraction
  $\cX\to\cY$. 
\end{defn}
Fix a split-regular ideal $\cI\sset\O_\cX$ and and ideal
$\cJ\sset\O_\cX$ adapted to $\cI$. To explain how the preceding
constructions patch together we can assume that $\cS=\Spf{R}$ is
affine. Let $U=\Spf{A}\sset\cX$ be an affine open and set
$I=\Gamma(U,\cI)$, $J=\Gamma(U,\cJ)$. For any $f\in A$ the
$A_n$-algebra $(A_n)_f$ is flat, and the natural morphism
\begin{displaymath}
  P_{(m)}(IA_n)_f\to P_{(m)}(I(A_n)_f)
\end{displaymath}
is an isomorphism. Thus if $J$ is an open ideal such that
$p^{n+1}\in J$, the natural morphism
\begin{displaymath}
  A_f\tens_AP_{J,(m)}(I)\to P_{J,(m)}(IA_f)
\end{displaymath}
is an isomorphism. It follows that there is a quasicoherent sheaf of
$\O_{X_J}$-algebras $\cP_{\cJ,(m)}(\cI)$ with an $m$-PD-structure such
that for affine opens $U=\Spf{A}\sset\cX$, $V=\Spf{R}\sset\cS$ such
that $\cX\to\cS$ sends $U\to V$,
\begin{displaymath}
  \Gamma(U,\cP_{\cJ,(m)}(\cI))=P_{J,(m)}(I)
\end{displaymath}
where $I=\Gamma(U,\cI)$ and $J=\Gamma(U,\cJ)$. The
$\O_{X_\cJ}$-algebra $\cP_{\cJ,(m)}(\cI)$ gives us an affine morphism
$\cX_{\cS,(m)}^\cJ(\cI)\to\cX$ of formal schemes, where
\begin{equation}
  \label{eq:X_(m),SJ}
  \cX_{\cS,(m)}^\cJ(\cI):=Spec_{\O_{X_\cJ}}(\cP_{\cJ,(m)}(\cI)).   
\end{equation}

Suppose now $\O_\cX$ has an ideal of definition adapted to $\cI$. By
lemma \ref{lemma:power-of-adapted-is-adapted} it has a cofinal set of
ideals of definition adapted to $\cI$, and we define
\begin{equation}
  \label{eq:P_A/R-global}
  \cP_{(m)}(\cI)=\liminv_\cJ\cP_{\cJ,(m)}(\cI)
\end{equation}
where $\cJ$ runs through the set of ideals of definition adapted to
$\cI$.  If $U=\Spf{A}\sset\cX$ is an open affine mapping to an open
affine $\Spf{A}\sset\cS$ then
\begin{displaymath}
  \Gamma(U,\cP_{(m)}(\cI))\simeq\hat P_{(m)}(I)
\end{displaymath}
where as before $I=\Gamma(U,\cI)$. By definition $\cP_{(m)}(\cI)$ is a
sheaf of $\O_\cX$-algebras whose reduction modulo $J\subset\O_\cX$ for
any $J$ adapted to $\cI$ is the sheaf $\cP_{J,(m)}(\cI)$. We will not
attach a formal scheme to $\cP_{(m)}(\cI)$ since this would take us
out of the category of locally noetherian formal schemes. As before
there is a change-of-level morphism
\begin{equation}
  \label{eq:Pm-change-of-level}
  \cP_{(m')}(\cI)\to \cP_{(m)}(\cI)
\end{equation}
for $m'\ge m$.

We can now state the universal properties of $\cX_{(m),\cS}^J\to\cX$
and $\O_\cX\to\cP_{(m)}(\cI)$ when $J$ is adapted to $\cI$.
Suppose, first, that $X'$ is an $\cS$-\textit{scheme} with an
$m$-PD-structure $(I',J',\gamma')$ compatible with $(\fa,\fb,\alpha)$,
and $f:X'\to\cX$ is an $\cS$-morphism such that $f^*\cI\sset I'$.
There is a cofinal set of ideals of definition $J\subset\O_\cX$
satisfying \ref{eq:special-J1} such that $f$ factors
\begin{displaymath}
  X'\to \cX_{(m),\cS}^J(\cI)\Xto{f_J}\cX  
\end{displaymath}
for some unique $m$-PD-morphism $f_J$. Suppose, on the other hand that
$\cX'$ is a formal $\cS$-scheme with an $m$-PD-structure
$(\cI',\cJ',\gamma')$ compatible with $(\fa,\fb,\alpha)$, and
$f:X'\to\cX$ is an $\cS$-morphism such that $f^*\cI\sset\cI'$. Then
the canonical morphism $f^*\O_\cX\to\O_{\cX'}$ has a unique
factorization
\begin{displaymath}
  f^*\O_\cX\to f^*\cP_{(m)}(\cI)\Xto{g}\O_{\cX'}
\end{displaymath}
where $g$ is an $m$-PD-morphism (here the $f^*$ is understood in the
sense of ringed spaces).

\subsection{The sheaf $\cP_{\cX/\cS,(m)}(r)$.}
\label{sec:P_X/S,m}

When $R\to A$ is a quasi-smooth homomorphism of adic noetherian
$\bZ_p$-algebras we can apply the preceding constructions to the
diagonal ideal $I(r)$ of the completed tensor product $\hat A_R(r)$ of
$r+1$ copies of $A$ over $R$. Recall $I(r)$ is the kernel of the
multiplication map $A_R(r)\to A$, which has $r+1$ sections, namely the
maps $d_i:A\to A_R(r)$ for $0\le i\le r$. Set $R_n=R/p^{n+1}R$ and
$A_n=A/p^{n+1}$ as before, and set $I_n=I(A_n)_{R_n}(r)$. Since
$R_0\to A_0=A\tens_{\bZ_p}\bF_p$ is quasi-smooth, $A_n$ is a flat
$R_n$-algebra and $I_n$ is a regular ideal. It follows that the ideal
$I(r)\sset\hat A_R(r)$ is split-regular.

If $J\subset A$ is any ideal of definition, 
\begin{equation}
  \label{eq:Jr-A}
  J(r)=\sum_{0\le i\le r}d_i(J)A_R(r)
\end{equation}
is an ideal of definition of $A_R(r)$; that it is adapted to $I(r)$ 
with $\sigma=d_i$ for any $i$ follows from:

\begin{lemma}\label{lemma:bilaterals-and-envelopes}
  Suppose $R\to A$ is quasi-smooth and $(\fa,\fb,\alpha)$ is an
  $m$-PD-structure on $R$. If $J\sset A$ is an $m$-bilateralising
  ideal,
  \begin{displaymath}
    J(r) P_{A_n/R_n,(m)}=d_i(J) P_{A_n/R_n,(m)}
  \end{displaymath}
  for $0\le i\le r$ and any $n$ such that $p^n\in J$.
\end{lemma}
\begin{demo}
  As before we reduce to the case where $p^n=0$ in $A$. 
  The canonical isomorphisms
  \begin{displaymath}
    P_{A/R,(m)}(r)\ctens_A P_{A/R,(m)}(r')
    \simeq P_{A/R,(m)}(r+r')
  \end{displaymath}
  show that it suffices to treat the case $r=1$, which may be
  rephrased as an equality
  \begin{equation}
    \label{eq:J2}
    JP_{A/R,(m)}=P_{A/R,(m)}J
  \end{equation}
  for all $m$-bilateralising $J$. Set $P=P_{A/R,(m)}$; in the notation
  of corollary \ref{cor:bilaterals-and-envelopes} we must show that
  $\delta(J)\sset JP\cap PJ$, and the same Taylor series argument
  shows that $\delta(J)\sset JP$. By Zariski localization we may
  assume that $A$ has local coordinates relative to $R$, which we may
  use to get a basis $\{\dpniv{\xi}{K}{m}\}_{K\ge0}$ of $P$ as a
  $B$-module via $d_1:B\to P$. The corollary tells us that the image
  of $\delta(J)$ under the natural projection $P\to P^n=P^n_{A/R,(m)}$
  is contained in $P^nJ$ for all $n$. Thus if $x\in\delta(J)$ is
  $\sum_Kd_1(a_K)\dpbrniv{\xi}{I}{m}$ in terms of the basis we have
  $a_K\in J$ for all $K$, and thus $x\in PJ$.
\end{demo}

\begin{remark}
  The equality \ref{eq:J2} looks like the definition of
  ``bilateralising'' but has nothing to do with it, since in fact
  $P_{B/A,(m)}$ is a commutative ring. The ideals $JP_{B/A,(m)}$,
  $P_{B/A,(m)}J$ are the ideals generated by the image of $J$ under
  the two ring homomorphisms $d_0$, $d_1:B\to B_{B/A,(m)}$.
\end{remark}

Suppose now $\cS$ is an adic locally noetherian formal $\bZ_p$-scheme
with $m$-PD-structure $(\fa,\fb,\alpha)$ and $\cX\to\cS$ is
quasi-smooth. We may apply the results of \S\ref{sec:P_m(I)}, with the
result that for any bilateralising $J\subset\O_\cX$ there is a sheaf
$\cP_{\cX/\cS,J,(m)}(r)$ with $r+1$ $\O_\cX$-module structures,
quasicoherent for any one of them. The sheaf $\cP_{\cX/\cS,(m)}(r)$
is the inverse limit
\begin{equation}
  \label{eq:P_X/S,m}
  \cP_{\cX/\cS,(m)}(r)=\liminv_J\cP_{\cX/\cS,J,(m)}(r)
\end{equation}
where $J$ runs through all bilateralising ideals of definition (this
is another case of dropping the hat in a geometric context). By
construction it is an algebra over the structure sheaf of the
completion $\hat\cX_\cS(r)$ of the $r+1$-fold fiber product
$\cX_\cS(r)$ with respect to the diagonal ideal $\cI(r)$. When
$U=\Spf{A}\sset\cX$ is an open affine lying over $\Spf{R}\sset\cS$ we
have
\begin{equation}
  \label{eq:P_X/S,m}
  \Gamma(U,\cP_{\cX/\cS,(m)}(r))\simeq\hat P_{B/A,(m)}(r)
\end{equation}
As before we omit the $(r)$ when $r=1$


The scheme
\begin{displaymath}
  \cX^J_{\cS,(m)}=\Spalg{\O_{X_J}}{\cP_{\cX/\cS,J,(m)}}
\end{displaymath}
is relatively affine over $(X_J)_{\cS}(1)$. We denote by $q_0$,
$q_1:\cX^J_{\cS,(m)}\to X_J$ the composites of the structure morphism
$\cX^J_{\cS,(m)}\to (X_J)_{\cS}(1)$ with the natural projections
$p_0$, $p_1:(X_J)_{\cS}(1)\to X_J$. The inductive system of
$\cX^J_{\cS,(m)}$ has the following universal property. Let $Y$ be a
\textit{scheme} over $\cS$ such that the $m$-PD-structure
$(\fa,\fb,\alpha)$ of $\cS$ extends to $Y$, and let $f_0$,
$f_1:Y\to\cX$ be $\cS$-morphisms congruent modulo $\fa$ in the sense
that if $Y_0\subset Y$ is the closed subscheme defined by $\fa\O_Y$,
the two composite morphisms
\begin{displaymath}
  \xymatrix{
    Y_0\ar[r]&Y\ar@/^/[r]^{f_0}\ar@/_/[r]_{f_1}&\cX
  }
\end{displaymath}
are equal. There is a cofinal set of $J\subset\O_\cX$ such that there
is a unique morphism $g_J:Y\to\cX^J_{\cS,(m)}$ such
that the morphism $(f_0,f_1):Y\to\cX\times_\cS\cX$ factors
\begin{displaymath}
  Y\Xto{g_J}\cX^J_{\cS,(m)}\to X_J\times_\cS X_J
  \to\cX\times_\cS\cX.
\end{displaymath}

\subsection{$m$-HPD-stratifications.}
\label{sec:m-HPD-stratifications}

As with $m$-PD-stratifications, we approach the topic of topologically
quasi-nilpotent $\niv\D m$-modules and $m$-HPD-stra\-ti\-fi\-cation by a
limiting procedure. As before, if $J\subset\O_\cX$ is an ideal of
definition we denote by $X_J\subset\cX$ the closed subscheme
corresponding to $J$.  Denote by $p_0$, $p_1:\cX_\cS(1)\to\cX$ the
natural projections.

\begin{defn}\label{defn:m-HPD-strat-modJ}
  An \textit{$m$-HPD-stratification} of an $\O_{X_J}$-module
  $M$ relative to $\cS$ is an isomorphism
  \begin{equation}
    \label{eq:m-HPD-stratification}
    \chi:\cP_{\cX/\cS,J,(m)}\tens_{\O_{\cX_\cS(1)}}p_1^*M
    \isom p_0^*M\tens_{\O_{\cX_\cS(1)}}\cP_{\cX/\cS,J,(m)}
  \end{equation}
  restricting to the identity on the diagonal and satisfying the cocycle
  condition. 
\end{defn}

If $\chi$ is an $m$-HPD-stratification, extending scalars
by $\cP_{\cX/\cS,(m)}\to\cP^n_{\cX/\cS,(m)}$ for all $n$ results in an
$m$-PD-stratification of $M$ relative to $\cS$, and $\chi$ is
determined by this $m$-PD-stratification. We may then say that a left
$\niv{\D}{m}_{X_J/\cS,J}$-module $M$ is \textit{quasi-nilpotent} if
its associated $m$-PD-stratification extends to an
$m$-HPD-stratification. The argument of
\cite[Prop. 2.3.7]{berthelot:1996} with $\cP_{\cX/\cS,J,(m)}$ in place
of $\cP_{\cX/\cS}$ then shows:

\begin{prop}\label{prop:nilpotence-criterion}
  A left $\niv{\D}{m}_{X_J/\cS}$-module $M$ is quasi-nilpotent if and
  only if for every local section $x$ of $M$ and some system of local
  coordinates (defined in the same neighborhood as $m$),
  $\dpabniv{\d}{I}{m}(x)=0$ for $|I|\gg0$. If this is so, then in fact
  $\dpabniv{\d}{I}{m}(x)=0$ for any system of local coordinates and
  $|I|\gg0$. \nodemo
\end{prop}

\begin{cor}\label{cor:nilpotence-criterion}
  If $M$ is a quasi-nilpotent left $\niv{\D}{m}_{X_J/\cS}$-module,
  then so is any submodule or quotient module of $M$, and conversely
  if $M$ is a left $\niv{\D}{m}_{X_J/\cS}$-module and $N\sset M$ is
  a submodule such that $N$ and $M/N$ are quasi-nilpotent, then $M$ is
  quasi-nilpotent. If $M$ and $N$ are quasi-nilpotent left
  $\niv{\D}{m}_{X_J/\cS}$-modules then so are $M\tens_{\O_\cX}N$ and
  $Hom_{\O_\cX}(M,N)$.\nodemo
\end{cor}

For example, if $J\subset\O_\cX$ is $m$-bilateralising then the
standard $\niv{\D}{m}_{\cX/\cS}$-module structure of $\O_\cX$ induces
a quasi-nilpotent $\niv{\D}{m}_{X_J/\cS}$-module structure
$\O_\cX/J$. Applying this to the $\niv{\D}{m}_{\cX/\cS,J^n}$-module
$\O_\cX/J^n$, we see from the corollary that $J^n/J^{n+1}$ is a
quasi-nilpotent $\niv{\D}{m}_{X_J/\cS}$-module for all $n\ge0$.

The argument of \cite{berthelot:2000} shows:

\begin{prop}\label{prop:qu-nilp-and-m-HPD}
  The category of quasi-nilpotent left $\niv\D m_{X_J/\cS}$-modules is
  equivalent to the category of $\O_{X_J}$-modules endowed with an
  $m$-HPD-stra\-ti\-fi\-cation. \nodemo
\end{prop}

We now pass to the limit; as before we restrict to the case of
pro-quasicoherent modules.

\begin{defn}\label{defn:topologically-quasi-nilpotent}
  A pro-quasicoherent left $\niv{\hD}{m}_{X/\cS}$-module $M$ is
  \emph{topologically quasi-nilpotent} for every bilateralising ideal
  of definition $J\subset\O_\cX$ the reduction $M/JM$ is a
  quasi-nilpotent $\niv{\hD}{m}_{X_J/\cS}$-module.
\end{defn}


\begin{prop}\label{prop:top-nilpotence-criterion}
  A pro-quasicoherent left $\niv{\D}{m}_{X/\cS}$-module $M$ is
  topologically quasi-nilpotent if and only if $M/JM$ is a
  quasi-nilpotent $\niv{\D}{m}_{X_J/\cS}$-module for \emph{some}
  $m$-bilateralising ideal of definition $J\subset\O_\cX$.
\end{prop}
\begin{demo}
  The condition is evidently necessary, and for the converse it
  suffices to show that if $M/JM$ is quasi-nilpotent for $J$ in the
  proposition then $M/J^nM$ is a quasi-nilpotent
  $\niv{\D}{m}_{X_{J^n}/\cS}$-module for all $n>0$. We have seen that
  $J^k/J^{k+1}$ is a quasi-nilpotent $\niv{\D}{m}_{X_J/\cS}$-module
  for all $k$, so by the hypothesis and corollary
  \ref{cor:nilpotence-criterion} the same holds for the tensor product
  $J^k/J^{k+1}\tens_{\O_\cX}M/JM$ and for its quotient
  $J^kM/J^{k+1}M$. Then for $k<n$, $J^kM/J^{k+1}M$ is quasi-nilpotent
  as a $\niv{\D}{m}_{X_{J^n}/\cS}$-module, and corollary
  \ref{cor:nilpotence-criterion} shows that $M/J^n$ is a
  quasi-nilpotent $\niv{\D}{m}_{X_{J^n}/\cS}$-module as well.
\end{demo}

\begin{cor}\label{cor:quasi-nilpotence-criterion}
  If $M$ is a topologically quasi-nilpotent left
  $\niv{\hD}{m}$-module, so is any submodule of $M$, or any quotient
  module of $M$ that is pro-quasicoherent. If $M$ is a left
  $\niv{\hD}{m}$-module and $N\sset M$ is a submodule such that $N$
  and $M/N$ are topologically quasi-nilpotent, then so is $M$. If $M$
  and $N$ are topologically quasi-nilpotent $\niv{\hD}{m}$ then so are
  $M\ctens_{\O_\cX}N$ and $Hom_{\O_\cX}(M,N)$ if they are
  pro-quasicoherent.
\end{cor}
\begin{demo}
  This follows from proposition \ref{prop:top-nilpotence-criterion}
  and corollary \ref{cor:nilpotence-criterion}.
\end{demo}

We remark that the pro-quasicoherence assumptions are automatic when
$M$ and $N$ are coherent $\niv{\hD}{m}$-modules.

\begin{cor}\label{cor:forcing-quasi-nilpotence}
  If $M$ is a pro-quasicoherent left $\niv{\hD}{m+1}_{X/\cS}$-module,
  the $\niv{\hD}{m}_{X/\cS}$-module induced by restriction is
  topologically quasi-nilpotent.
\end{cor}
\begin{demo}
  Topological quasi-nilpotence can be detected locally, so we can work
  in local coordinates. It suffices to observe that in the ring
  $\niv{\hD}{m+1}_{X/\cS}$, the differential operators
  $\dpabniv{\d_i}{p^m}{m+1}$ are multiples of $p$; since $p$ is
  topologically nilpotent in $\O_\cX$ we can apply proposition
  \ref{prop:nilpotence-criterion}.
\end{demo}

Passing to the limit in the definition \ref{defn:m-HPD-strat-modJ}
results in:

\begin{defn}\label{defn:m-HPD-pcoh}
  An $m$-HPD-structure on a pro-quasicoherent 
  $\O_\cX$-module $M$ is an isomorphism
  \begin{equation}
    \label{eq:m-HPD-stratification-limit}
    \chi:\cP_{\cX/\cS,(m)}\ctens_{\O_\cX}M\isom
    M\ctens_{\O_\cX}\cP_{\cX/\cS,(m)} 
  \end{equation}
  reducing to the identity on the diagonal and satisfying the cocycle
  condition.
\end{defn}
Note that the completed tensor product is necessary. 

Similarly, passing to the limit in proposition
\ref{prop:qu-nilp-and-m-HPD} yields:

\begin{prop}\label{prop:top-qu-nilp-and-m-HPD}
  The category of topologically quasi-nilpotent left
  $\niv\D m_{X_J/\cS}$-modules is equivalent to the category of
  $\O_{X_J}$-modules endowed with an $m$-HPD-stratification. \nodemo
\end{prop}
\begin{demo}
  If $M$ is topologically quasi-nilpotent (and in particular
  pro-quasi\-coherent), $M/JM$ is quasi-nilpotent for all $J$, whence an
  $m$-HPD-stra\-ti\-fi\-cation $\chi_J$ for all $J$. Since a left
  $\niv{\D}{m}_{X/\cS}$-module structure arises from at most one
  $m$-HPD-stratification on an $\O_{X_J}$-module, the various $\chi_J$
  of the $M/JM$ for variable $J$ must all be compatible, and then
  $\chi=\liminv_J\chi_J$ has the required properties. Conversely an
  isomorphism \ref{eq:m-HPD-stratification-limit} induces
  $m$-HPD-stratifications of $M/JM$ for all $J$, all of which
  correspond to the same left $\niv{\D}{m}_{X/\cS}$-module structure.
\end{demo}

As in the smooth case, the main interest of (topological)
quasi-nilpotence is the following invariance property.

\begin{prop}\label{prop:invariance-of-base-change1}
  Suppose $\cX\to\cS$ is quasi-smooth and $(\fa,\fb,\alpha)$ is an
  $m$-PD-structure on $\cS$, and $M$ is a left
  $\niv{\hD}{m}_{\cX/\cS}$-module. Suppose $\cX'$ is a formal
  $\cS$-scheme such that $(\fa,\fb,\alpha)$ extends to $\cX'$, and let
  $f$, $f':\cX'\to\cX$ be two $\cS$-morphisms having the same
  restriction to the closed formal subscheme $X'_0\subset\cX'$ defined
  by $\fa\O_{\cX'}$.
  \begin{enumerate}
  \item If the $m$-PD-ideal $(\fa,\fb,\alpha)$ is $m$-PD-nilpotent,
    there is a canonical isomorphism
    \begin{displaymath}
      \tau_{f,f'}:(f')^*M\isom f^*M
    \end{displaymath}
    of left $\O_{\cX'}$-modules, such that $\tau_{f,f}=id_{f^*M}$, and
    the system of $\tau_{f,f'}$ is transitive in the sense that if
    $f'':\cX'\to\cX$ is a third such morphism then
    $\tau_{f,f'}\circ\tau_{f',f''}=\tau_{f,f''}$.
  \item If $(\fa,\fb,\alpha)$ is not assumed to be $m$-PD-nilpotent,
    but if $M$ is a topologically quasi-nilpotent left
    $\niv{\hD}{m}_{\cX/\cS}$-module, the system of $\tau_{f,f'}$ also
    exists and has the same properties.
  \item In either case, in the situation of the diagram
    \ref{eq:base-change} when $\cX'/\cS'$ is quasi-smooth, the
    morphism $\tau_{f,f'}$ is $\niv{\D}{m}_{\cX'/\cS'}$-linear.
  \end{enumerate}
\end{prop}

The argument is the same as in \cite[Prop. 2.1.5]{berthelot:2000}, and
is a direct consequence of the universal property of the system of
$\cX^J_{\cS,(m)}$. In (iii), the proof of $\niv{\D}{m}$-linearity uses
the compatibility of the canonical ideal of $\cX^n_{\cS,(m)}$ with the
$m$-PD-structure of $\fa$, which relies on the structure theorem
\cite[Prop. 1.5.3]{berthelot:1996} for $m$-PD-envelopes of a regular
ideal.

\subsubsection{Coefficient Rings.}
\label{sec:coefficients}

As in \cite[2.3.4]{berthelot:1996}, a left
$\niv{\D}{m}_{\cX/\cS}$-module structure on a commutative
$\O_\cX$-algebra $\cB$ is \textit{compatible with its algebra
  structure} if the given $\O_\cX$-module structure of $\cB$ coincides
with the one derived from its left $\niv{\D}{m}_{\cX/\cS}$-module
structure, and if the isomorphisms \ref{eq:m-PD-stratification2} are
isomorphisms of $\cP^n_{\cX/\cS,(m)}$-algebras. An equivalent
condition is that the multiplication map $\cB\tens_{\O_\cX}\cB\to\cB$
is $\niv{\D}{m}_{\cX/\cS}$-linear. In local coordinates, this is
equivalent to the level $m$ Leibnitz rule
\begin{equation}
  \label{eq:Leibnitz-coefficient-rings}
  \dpabniv{\d}{K}{m}(ab)
  =\sum_{I+J=K}\bbinom{K}{I}{m}\dpabniv{\d}{I}{m}(a)\dpabniv{\d}{J}{m}(b)
\end{equation}
holding for all local sections $a$, $b$ of $\cB$. 

When $\cB$ is an $\O_\cX$-algebra with a compatible left
$\niv{\D}{m}_{\cX/\cS}$-module structure, the $\O_\cX$-module
$\cB\tens_{\O_\cX}\niv{\D}{m}_{\cX/\cS}$ has a unique ring structure
such that the canonical homomorphisms
\begin{align*}
  \cB\to\cB\tens_{\O_\cX}\niv{\D}{m}_{\cX/\cS}\qquad &b\mapsto b\tens1\\
  \niv{\D}{m}_{\cX/\cS}\to\cB\tens_{\O_\cX}\niv{\D}{m}_{\cX/\cS}
  \qquad&P\mapsto 1\tens P\\
\end{align*}
are ring homomorphisms. The product is defined as follows: with the
identification
\begin{equation}
  \label{eq:DmB/s-mult1}
  \cB\tens_{\O_\cX}\Diff^n_{\cX/\cS,(m)}
  \simeq\Hom_{\O_\cX}(\cP^n_{\cX/\cS,(m)},\cB)
\end{equation}
the product of local sections
$P\in\cB\tens_{\O_\cX}\Diff^{n'}_{\cX/\cS,(m)}$,
$Q\in\cB\tens_{\O_\cX}\Diff^n_{\cX/\cS,(m)}$ is
\begin{equation}
  \label{eq:DmB/s-mult2}
  \begin{split}
    \cP^{n+n'}_{\cX/\cS,(m)}
    &\Xto{\delta^{n,n'}}\cP^{n'}_{\cX/\cS,(m)}\tens\cP^n_{\cX/\cS,(m)}\\
    &\Xto{1\tens Q}\cP^{n'}_{\cX/\cS,(m)}\tens\cB\\
    &\Xto{\chi_{n'}}\cB\tens\cP^{n'}_{\cX/\cS,(m)}\\
    &\Xto{P}\cB\tens\cB\to\cB
  \end{split}
\end{equation}
where $\chi$ is the stratification of $\cB$.

Since on occasion we will be considering
several $\cB$ at once we will use the notation
\begin{equation}
  \label{eq:D-with-coefficients}
  \niv{\D}{m}_{\cB/\cS}=\cB\tens_{\O_\cX}\niv{\D}{m}_{\cX/\cS}
\end{equation}
in preference to that of \cite{berthelot:1996}. For any open affine
$U\sset\cX$ the canonical homomorphism
\begin{equation}
  \label{eq:coefficient-ring}
  \Gamma(U,\cB)\tens_{\Gamma(U,\O_\cX)}\Gamma(U,\niv{\D}{m}_{\cX/\cS})\to
  \Gamma(U,\niv{\D}{m}_{\cB/\cS})
\end{equation}
is an isomorphism; argument is the same as that of
\cite[Prop. 2.3.6]{berthelot:1996} and depends mainly on the fact that
$\niv{\D}{m}_{\cX/\cS}$ is an inductive limit of locally free
$\O_\cX$-modules. 

From now on the following conditions will be imposed in $\cB$,
without explicit mention to the contrary:
\begin{equation}\label{enum:coeff-ring-B1}
  \begin{minipage}[t]{0.7\linewidth}for every open affine $U\sset\cX$,
    $\Gamma(U,\cB)$ is noetherian;
  \end{minipage}
\end{equation}
\begin{equation}\label{enum:coeff-ring-B2}
  \begin{minipage}[t]{0.7\linewidth}{
      $\cB$ is pro-quasicoherent.}
  \end{minipage}
\end{equation}
We define
\begin{equation}
  \label{eq:coefficient-ring-complete}
  \begin{split}
    \niv{\hD}{m}_{\cB/\cS}&=\liminv_J(\cB/J\cB)
                            \tens_{\O_\cX}\niv{\D}{m}_{\cX/\cS}\\
    &=\cB\ctens_{\O_\cX}\niv{\hD}{m}_{\cX/\cS}
  \end{split}
\end{equation}
where the inverse limit is over $m$-bilateralising ideals of
definition. Then $\niv{\hD}{m}_{\cB/\cS}$ is a ring, and the previous
discussion shows that a left $\niv{\hD}{m}_{\cB/\cS}$-module is the
same as $\cB$-module with a compatible left
$\niv{\hD}{m}_{\cX/\cS}$-module structure in the previous sense. With
the hypotheses \ref{enum:coeff-ring-B1}--\ref{enum:coeff-ring-B2} on
$\cB$, theorem \ref{thm:coherence-finite-level}, \ref{thm:theorem-A}
and its corollary, theorem \ref{thm:thmB} and propositions
\ref{prop:coherent-modules-level-m}, \ref{prop:triangle-functor} hold
for $\niv{\hD}{m}_{\cB/\cS}$ without modification; we will not bother
to restate them.

For homomorphisms of pro-quasicoherent $\O_\cX$-algebras
we can define pullbacks as in rigid geometry,
i.e.\ by reducing modulo powers of an ideal of definition and passing
to a limit. The same applies to tensor products; we will explicitly
note completion of the tensor products unless one factor is coherent,
in which case completion is unnecessary. We can then reformulate the
condition that an $\O_\cX$-algebra $\cB$ has an $m$-PD-stratification
of an $\O_\cX$ compatible with its algebra structure. Suppose
$\cP^n_{\cB/\cS,(m)}$ is an $\O_{\cP^n_{\cX/\cS,(m)}}$-algebra and
\begin{equation}
  \label{eq:cP_B^n}
  \begin{split}
    \alpha^n_0:p^*_0(\cB)\tens_{\O_{\hat\cX(1)}}\cP^n_{\cX/\cS,(m)}&\isom\cQ\\
    \alpha^n_1:\cP^n_{\cX/\cS,(m)}\tens_{\O_{\hat\cX(1)}}p^*_1(\cB)&\isom\cQ
  \end{split}
\end{equation}
are isomorphisms of $\cP^n_{\cX/\cS,(m)}$-algebras. We will say that
$\{\alpha^n_0\}_{n\ge0}$ and $\{\alpha^n_1\}_{n\ge0}$ are compatible
if for $n'\ge n$ they define the same morphism
$\cP^{n'}_{\cB/\cS,(m)}\to\cP^n_{\cB/\cS,(m)}$ and for $n=0$ they
yield the same identification $\cP^0_{\cB/\cS,(m)}\simeq\cB$. If
$\{\alpha^n_0\}_{n\ge0}$ and $\{\alpha^n_1\}_{n\ge0}$ are compatible,
the isomorphisms
\begin{equation}
  \label{eq:alpha-to-chi}
  \chi_n=(\alpha^n_0)^{-1}\circ\alpha^n_1:
  \cP^n_{\cX/\cS,(m)}\tens_{\O_{\hat\cX(1)}}\cB
  \to\cB\tens_{\O_{\hat\cX(1)}}\cP^n_{\cX/\cS,(m)}  
\end{equation}
are compatible in the previous sense. Conversely if $\chi_n$ is given
we set $\cP^n_{\cB/\cS,(m)}=\cB\tens_{\O_\cX}\cP^n_{\cX/\cS,(m)}$; then
$\alpha^n_0=id$ and $\alpha^n_0=\chi_n$ are compatible. When
$\alpha^n_0$ and $\alpha^n_1$ are compatible we will take the
isomorphisms \ref{eq:cP_B^n} to define $\cP^n_{\cB/\cS,(m)}=\cQ$.

In any case the isomorphisms \ref{eq:cP_B^n} give
$\cP^n_{\cB/\cS,(m)}$ a $(p^*_0\cB,p^*_1\cB)$-bimodule structure, the
left (resp. right) one arising from $\alpha^n_0$
(resp. $\alpha^n_1$). These structures are exchanged by the
isomorphism $\chi_n=(\alpha^n_0)^{-1}\circ\alpha^n_1$, and by
construction $\chi_n$ is $\cP^n_{\cX/\cS,(m)}$-linear.

Given \ref{eq:cP_B^n} one can define two ring homomorphisms
\begin{displaymath}
  p^*_{02}\cP^{n+n'}_{\cB/\cS,(m)}\to
  p^*_{01}\cP^n_{\cB/\cS,(m)}\tens_{p^*_1\cB} p^*_{12}\cP^{n'}_{\cB/\cS,(m)}
\end{displaymath}
where $p_{ij}:\hat\cX(2)\to\hat\cX(1)$ and $p_{i}:\hat\cX(2)\to\cX$
are the usual projections. We first remark that we can identify
\begin{align*}
  &p^*_{01}\cP^n_{\cX/\cS,(m)}\tens_{\O_{\hat\cX(2)}}p^*_1\cB
  \tens_{\O_{\hat\cX(2)}}p^*_{12}\cP^{n'}_{\cX/\cS,(m)}\\
  &\qquad\isom(p^*_{01}\cP^n_{\cX/\cS,(m)}\tens_{\O_{\hat\cX(2)}}p^*_1\cB)
    \tens_{p^*_1\cB}(p^*_1\cB\tens_{\O_{\hat\cX(2)}}p^*_{12}\cP^{n'}_{\cX/\cS,(m)})\\
  &\qquad\underset{\sim}{\Xto{\alpha^n_1\tens\alpha^{n'}_0}}
  (p^*_{01}\cP^n_{\cB/\cS,(m)})\tens_{p^*_1\cB}(p^*_{12}\cP^{n'}_{\cB/\cS,(m)}).
\end{align*}
The first homomorphism is the composite
\begin{equation}
  \label{eq:P_B-delta1}
  \begin{split}
    \delta^{n,n'}_{\cB,0}:p^*_{02}\cP^{n+n'}_{\cB/\cS,(m)}
    &\Xto{(\alpha^{n+n'}_0)^{-1}}p^*_0\cB
      \tens_{\O_{\hat\cX(2)}}p^*_{02}\cP^{n+n'}_{\cX/\cS,(m)}\\
    &\Xto{\delta^{n,n'}\tens1}p^*_0\cB\tens_{\O_{\hat\cX(2)}}
      p^*_{01}\cP^n_{\cX/\cS,(m)}
    \tens_{\O_{\hat\cX(2)}}p^*_{12}\cP^{n'}_{\cX/\cS,(m)}\\
    &\Xto{\chi_n^{-1}\tens 1}
      p^*_{01}\cP^n_{\cX/\cS,(m)}\tens_{\O_{\hat\cX(2)}}p^*_1\cB
      \tens_{\O_{\hat\cX(2)}}p^*_{12}\cP^{n'}_{\cX/\cS,(m)}\\
    &\simeq
      (p^*_{01}\cP^n_{\cB/\cS,(m)})\tens_{p^*_1\cB}(p^*_{12}\cP^{n'}_{\cB/\cS,(m)})  
  \end{split}
\end{equation}
and the second is the composite
\begin{equation}
  \label{eq:P_B-delta2}
  \begin{split}
    \delta^{n,n'}_{\cB,1}:p^*_{02}\cP^{n+n'}_{\cB/\cS,(m)}
    &\Xto{\alpha^{n+n'}_1}
      p^*_{02}\cP^{n+n'}_{\cX/\cS,(m)}\tens_{\O_{\hat\cX(2)}}p^*_2\cB\\
    &\Xto{\delta^{n,n'}\tens1}p^*_{01}\cP^n_{\cX/\cS,(m)}
      \tens_{\O_{\hat\cX(2)}}p^*_{12}\cP^{n'}_{\cX/\cS,(m)}
      \tens_{\O_{\hat\cX(2)}}p^*_2\cB\\    
    &\Xto{\chi_{n'}\tens 1}p^*_{01}\cP^n_{\cX/\cS,(m)}
      \tens_{\O_{\hat\cX(2)}}p^*_1\cB
    \tens_{\O_{\hat\cX(2)}}p^*_{12}\cP^{n'}_{\cX/\cS,(m)}\\
    &\simeq
      (p^*_{01}\cP^n_{\cB/\cS,(m)})\tens_{p^*_1\cB}(p^*_{12}\cP^{n'}_{\cB/\cS,(m)})  
  \end{split}
\end{equation}
Note that $\cP^n_{\cB/\cS,(m)}\tens_\cB\cP^{n'}_{\cB/\cS,(m)}$ has a
$(p^*_0\cB,p^*_2\cB)$-bimodule structure, coming from $d_0\tens1$ on
the left and and $1\tens d_1$ on the right. By construction,
$\delta^{n,n'}_{\cB,0}$ is $\cB$-linear for the left structure and
$\delta^{n,n'}_{\cB,1}$ is $\cB$-linear for the right structure.

\begin{prop}\label{prop:m-PD-stratified-algebra}
  For any sheaf $\cB$ of $\O_\cX$-algebras $\cB$ satisfying conditions
  \ref{enum:coeff-ring-B1}-\ref{enum:coeff-ring-B2}, and 
  for any compatible system of isomorphisms \ref{eq:cP_B^n}, the
  following are equivalent:
  \begin{enumerate}
  \item The isomorphisms \ref{eq:alpha-to-chi} define an
    $m$-PD-stratification of $\cB$ compatible with its
    $\O_\cX$-algebra structure.
  \item For all $n$, $n'\ge0$,
    $\delta^{n,n'}_{\cB,0}=\delta^{n,n'}_{\cB,1}$. 
  \item There is a ring homomorphism
    \begin{equation}
      \label{eq:delta^B1}
      \delta^{n,n'}_\cB:p^*_{02}\cP^{n+n'}_{\cB/\cS,(m)}
      \to (p^*_{01}\cP^n_{\cB/\cS,(m)})
      \tens_\cB (p^*_{12}\cP^{n'}_{\cB/\cS,(m)})
    \end{equation}
    that is a homomorphism of $(p^*_0\cB,p^*_2\cB)$-bimodules and
    semilinear for the homomorphism
    $\delta^{n,n'}:p^*_{02}\cP^{n+n'}_{\cX/\cS,(m)}\to
    p^*_{01}\cP^n_{\cX/\cS,(m)}
    \tens_{\O_\cX}p^*_{12}\cP^{n'}_{\cX/\cS,(m)}$.
  \end{enumerate}
\end{prop}
\begin{demo}
  In the following calculations we drop the $p^*_{ij}$, $p^*_i$, $(m)$
  and some of the the tensor product subscripts. Consider the diagram
  \begin{displaymath}
    \xymatrix{
      \cP^{n+n'}_{\cX/\cS}\tens{\cB}\ar[rr]^{\chi_{n+n'}}
      \ar[d]_{\delta^{n,n'}\tens1}
      &&
      \cB\tens{\cP^{n+n'}_{\cX/\cS}}\ar[d]^{1\tens\delta^{n,n'}}\\
      \cP^n_{\cX/\cS}\tens{\cP^{n'}_{\cX/\cS}}\tens{\cB}
      \ar[r]^{1\tens\chi_{n'}}
      \ar[d]_{1\tens\chi_{n'}}
      &\cP^n_{\cX/\cS}\tens{\cB}\tens{\cP^{n'}_{\cX/\cS}}
      \ar[r]^{\chi_n\tens1}\ar@{=}[d]
      &\cB\tens{\cP^n_{\cX/\cS}}\tens{\cP^{n'}_{\cX/\cS}}
      \\
      \cP^n_{\cX/\cS}\tens\cB\tens{\cP^n_{\cX/\cS}}
      \ar@{=}[r]
      &(\cP^n_{\cX/\cS}\tens{\cB})\tens_\cB(\cB
      \tens{\cP^n_{\cX/\cS}})\ar@{=}[r]
      &\cP^n_{\cX/\cS}\tens{\cB}\tens{\cP^n_{\cX/\cS}}
      \ar[u]_{\chi_n\tens1}
    }
  \end{displaymath}
  where the equalities denote the appropriate canonical
  isomorphisms. By \cite[Prop. 2.3.2]{berthelot:1996} the cocycle
  condition for $\chi_n$ is equivalent to the commutativity of the top
  rectangle. The equality
  $\delta^{n,n'}_{\cB,0}=\delta^{n,n'}_{\cB,1}$ is equivalent to the
  commutativity of the outside square. Since all morphisms in lower
  part of the diagram are isomorphisms, (i) and (ii) are equivalent.

  Since $\chi_n$ (resp. $\chi_{n'}$) is $\cP^n_{\cX/\cS,(m)}$-linear
  (resp. $\cP^{n'}_{\cX/\cS,(m)}$-linear) the morphisms
  $\delta^{n,n'}_{\cB,0}$ and $\delta^{n,n'}_{\cB,1}$ are semilinear
  for $\delta^{n,n'}$. Since $\delta^{n,n'}_{\cB,0}$
  (resp. $\delta^{n,n'}_{\cB,1}$) $\cB$-linear for the left
  (resp. right) structure, (ii) implies (iii) with
  $\delta^{n,n'}_\cB= \delta^{n,n'}_{\cB,0}=\delta^{n,n'}_{\cB,1}$.
  Suppose conversely that (iii) holds. Since
  $\cP^{n+n'}_{\cB/\cS,(m)}$ is generated as a $\cB$-module (for
  either structure) by the image of
  $\cP^{n+n'}_{\cX/\cS,(m)}\to\cP^{n+n'}_{\cB/\cS,(m)}$,
  $\delta^{n,n'}_{\cB,0}$ (resp. $\delta^{n,n'}_{\cB,1}$) is the
  unique morphism that is semilinear for $\delta^{n,n'}$ and
  $\cB$-linear for the left (resp. rignt) $\cB$-structure. Thus (iii)
  implies (ii)
\end{demo}

We will let the reader check that the definition of the product given
by \ref{eq:DmB/s-mult1} and \ref{eq:DmB/s-mult2} can be rephrased in
terms of $\cP^n_{\cB/\cS,(m)}$ as follows. The isomorphism
\ref{eq:DmB/s-mult1} may be rewritten
\begin{equation}
  \label{eq:DmB/s-mult3}
  \cB\tens_{\O_\cX}\Diff^n_{\cX/\cS,(m)}\simeq\Hom_\cB(\cP^n_{\cB/\cS,(m)},\cB)
\end{equation}
and with this identification the product of
$P\in\Hom_\cB(\cP^{n'}_{\cB/\cS,(m)},\cB)$ and
$Q\in\Hom_\cB(\cP^n_{\cB/\cS,(m)},\cB)$ is
\begin{equation}
  \label{eq:DmB/s-mult4}
  \begin{split}
    \cP^{n+n'}_{\cB/\cS,(m)}
    &\Xto{\delta^{n,n'}_\cB}\cP^{n'}_{\cB/\cS,(m)}\tens_\cB\cP^n_{\cB/\cS,(m)}\\
    &\Xto{1\tens Q}\cP^{n'}_{\cX/\cB,(m)}\Xto{P}\cB
  \end{split}
\end{equation}

We can treat $m$-HPD-stratified $\O_\cX$-algebras in the same way. A
pair of ring isomorphisms
\begin{equation}
  \label{eq:cP_B}
  \begin{split}
    \alpha_0:p^*_0\cB\ctens_{\O_\cX}\cP_{\cX/\cS,(m)}&\isom\cQ\\
    \alpha_1:\cP_{\cX/\cS,(m)}\ctens_{\O_\cX}p^*_1\cB&\isom\cQ
  \end{split}
\end{equation}
is \textit{compatible} if they induce the same surjective homomorphism
$\cP_{\cB/\cS,(m)}:=\cQ\to\cB$. Note the use of completed tensor
products in place of ordinary ones. As before they induce a
$(p^*_0\cB,p^*_1\cB)$-module structure on $\cP_{\cB/\cS,(m)}$ via
$d_0$ and $d_1$, and on
$(p^*_{01}\cP_{\cB/\cS,(m)})\ctens_{p^*_1\cB}(p^*_{12}\cP_{\cB/\cS,(m)})$
via $d_0\tens1$ and $1\tens d_1$. As before, one shows:

\begin{prop}\label{prop:m-HPD-stratified-algebra}
  A pair of compatible ring isomorphisms \ref{eq:cP_B} defines an
  $m$-HPD-stratification of $\cB$ compatible with its ring structure
  if and only there is a ring homomorphism
  \begin{equation}
    \label{eq:delta^B3}
    \delta_{\cB,(m)}:p^*_{02}\cP_{\cB/\cS,(m)}\to(p^*_{01}\cP_{\cB/\cS,(m)})
    \ctens_{p^*_1\cB} (p^*_{12}\cP_{\cB/\cS,(m)})
  \end{equation}
  that is $(p^*_0\cB,p^*_2\cB)$-bilinear and semilinear for the homomorphism
  \begin{displaymath}
    \delta:p^*_{02}\cP_{\cX/\cS,(m)}\to
    (p^*_{01}\cP_{\cX/\cS,(m)})\ctens_{\O_{\hat\cX(2)}}(p^*_{12}\cP_{\cX/\cS,(m)}).
  \end{displaymath}
    \nodemo
\end{prop}
In this connection one should recall the canonical
isomorphism
\begin{displaymath}
  p^*_{01}\cP_{\cX/\cS,(m)}\ctens_{\O_{\hat\cX(2)}}p^*_{12}\cP_{\cX/\cS,(m)}\simeq
  \cP_{\cX/\cS,(m)}(2)
\end{displaymath}
so that it would make sense to define
\begin{equation}
  \label{eq:P_B/S(2)}
  p^*_{01}\cP_{\cB/\cS,(m)}\ctens_{p^*_1\cB} p^*_{12}\cP_{\cB/\cS,(m)}\simeq
  \cP_{\cB/\cS,(m)}(2).
\end{equation}
Then \ref{eq:delta^B3} corresponds to a homomorphism
\begin{equation}
  \label{eq:d_02^B}
  d^\cB_{02}:p^*_{02}\cP_{\cB/\cS,(m)}\to\cP_{\cB/\cS,(m)}(2).
\end{equation}
On the other hand there are homomorphisms
\begin{equation}
  \label{eq:d_01,12^B}
  d^\cB_{01},\ d^\cB_{12}:\cP_{\cB/\cS,(m)}\to\cP_{\cB/\cS,(m)}(2)
\end{equation}
(we will start dropping the $p^*$ again) which via \ref{eq:delta^B3}
correspond to the morphisms $x\mapsto x\ctens1$ and
$x\mapsto 1\ctens x$ for $x$ in $\cP_{\cB/\cS,(m)}$. These
homomorphisms satisfy the simplicial identities
\begin{equation}
  \label{eq:d^B-simplicial-identities}
  d^\cB_0d^\cB_{02}=d^\cB_0d^\cB_{01},\quad
  d^\cB_1d^\cB_{02}=d^\cB_1d^\cB_{12},\quad
  d^\cB_0d^\cB_{12}=d^\cB_1d^\cB_{01}
\end{equation}
and the diagrams
\begin{equation}
  \label{eq:d^B-and-d}
  \xymatrix{
    \cP_{\cX/\cS,(m)}\ar[r]^{d_{ij}}\ar[d]
    &\cP_{\cB/\cS,(m)}(2)\ar[d]\\
    \cP_{\cB/\cS,(m)}\ar[r]_{d^\cB_{ij}}
    &\cP_{\cB/\cS,(m)}(2)\\
  }
\end{equation}
commute for $(ij)=(01)$, $(02)$ and $(12)$.

\subsubsection{Modules over coefficient rings.}
\label{sec:B-modules}

Suppose $M$ is a $\cB$-module satisfying
\ref{enum:module-coeff-ring-B1}, \ref{enum:module-coeff-ring-B2} and
$\cB$ has a compatible left $\niv{\D}{m}_{\cX/\cS}$-module
structure. An $m$-PD-stratification on $M$ can be viewed as a set of
isomorphisms
\begin{equation}
  \label{eq:B-compatible-m-PD-stratification}
  \chi_n:(\cP^n_{\cX/\cS,(m)}\tens_{\O_\cX} p^*_1\cB)\tens_{p^*_1\cB} M
  \isom
  M\tens_{p^*_1\cB}(p^*_0\cB\tens_{\O_\cX}\cP^n_{\cX/\cS,(m)})
\end{equation}
and we say that the $m$-PD-stratification of $M$ is compatible with
the $\cB$-module structure if the isomorphism
\ref{eq:B-compatible-m-PD-stratification} is semilinear with respect
to the stratification
\begin{displaymath}
  \cP^n_{\cX/\cS,(m)}\tens_{\O_\cX}p_1^*\cB
  \isom
  p_0^*\cB\tens_{\O_\cX}\cP^n_{\cX/\cS,(m)}
\end{displaymath}
of $\cB$. A left $\niv{\D}{m}_{\cX/\cS}$-module structure on $M$ is
compatible with the $\cB$-module structure if this is the case for the
corresponding $m$-PD-stratification. An equivalent condition is that
the map $\cB\tens_{\O_\cX}M\to M$ defining the $\cB$-module structure
is compatible with the $m$-PD-stratifications, or in other words is
$\niv{\D}{m}_{\cX/\cS}$-linear. In local coordinates, this
condition says that
\begin{equation}
  \label{eq:Leibnitz-B-modules}
  \dpabniv{\d}{K}{m}(ax)
  =\sum_{I+J=K}\bbinom{K}{I}{m}\dpabniv{\d}{I}{m}(a)\dpabniv{\d}{J}{m}(x)
\end{equation}
for local sections $a$ of $\cB$ and $x$ of $M$. 

A $\cB$-module with a compatible left $\niv{\D}{m}_{\cX/\cS}$-module
structure evidently gives rise to a left
$\niv{\D}{m}_{\cB/\cS}$-module structure. Conversely a left
$\niv{\D}{m}_{\cB/\cS}$-module $M$ gets compatible 
$\cB$-module and $\niv{\D}{m}_{\cX/\cS}$-module structures from the
canonical inclusions of $\cB$ and $\niv{\D}{m}_{\cX/\cS}$ into
$\niv{\D}{m}_{\cB/\cS}$.

If we are given a compatible pair of isomorphisms \ref{eq:cP_B^n}, the
isomorphisms \ref{eq:B-compatible-m-PD-stratification} defining the
stratification may be rewritten
\begin{equation}
  \label{eq:B-compatible-m-PD-stratification2}
  \chi^\cB_n:\cP^n_{\cB/\cS,(m)}\tens_{p^*_1\cB}p^*_1M\isom p^*_0M
  \tens_{p^*_0\cB}\cP^n_{\cB/\cS,(m)}
\end{equation}
and \ref{eq:B-compatible-m-PD-stratification} is semilinear for the
stratification of $\cB$ if and only if the isomorphisms
\ref{eq:B-compatible-m-PD-stratification2} are
$\cP^n_{\cB/\cS}$-linear. The compatibility of \ref{eq:cP_B^n}
guarantees the compatiblility of the $\chi^\cB_n$, while the cocycle
condition can be expressed in various ways. One is to introduce the
morphisms
\begin{equation}
  \label{eq:theta^B}
  \theta^\cB_n:p^*_1M\to p^*_0M\tens_{p^*_0\cB}\cP^n_{\cB/\cS,(m)}
\end{equation}
induced by \ref{eq:B-compatible-m-PD-stratification2}, which are
$\cB$-linear for the right structure of $\cP^n_{\cB/\cS,(m)}$; they
are compatible in an obvious sense, and $\chi_n$ satisfies the cocycle
condition if and only if diagram
\begin{equation}
  \label{eq:B-cocycle-condition}
  \xymatrix{
    M\ar[r]^{\theta^\cB_{n+n'}}\ar[d]_{\theta^\cB_{n'}}
    &M\tensu\cB\cP^{n+n'}_{\cB/\cS,(m)}\ar[d]^{\delta^{n,n'}_\cB}\\
    M\tensu\cB\cP^{n'}_{\cB/\cS,(m)}
    \ar[r]^{\theta^\cB_n}
    &M\tensu\cB\cP^n_{\cB/\cS,(m)}\tensu\cB\cP^{n'}_{\cB/\cS,(m)}
  }
\end{equation}
commutes for all $n$, $n'\ge0$. On the other hand the commutativity of
\ref{eq:B-cocycle-condition} shows that the $\theta^\cB_n$ give $M$
the structure of a left $\niv{\D}{m}_{\cB/\cS}$-module, and thus
yields another way of understanding the equivalence of the
category of left $\niv{\D}{m}_{\cB/\cS}$-modules with the category of
$\cB$-modules endowed with a 
left $\niv{\D}{m}_{\cX/\cS}$-module structure compatible with the
$\cB$-algebra structure. 

When $\cB$ has a quasi-nilpotent left $\niv{\D}{m}_{\cX/\cS}$-module
structure compatible with its algebra structure, the same picture
holds for quasi-nilpotent left $\niv{\D}{m}_{\cX/\cS}$-modules endowed
with a compatible $\cB$-module structure. The $m$-HPD-stratification
\begin{displaymath}
  \chi:\cP_{\cX/\cS,(m)}\ctens p^*_1M\isom
  p^*_0M\ctens\cP_{\cX/\cS,(m)}  
\end{displaymath}
can be rewritten
\begin{equation}
  \label{eq:B-compatible-HPD-stratification}
  \chi^\cB:\cP_{\cB/\cS,(m)}\ctens_{p^*_1\cB}p^*_1M\isom
  p^*_0M\ctens_{p^*_0\cB}\cP_{\cB/\cS,(m)}
\end{equation}
as before, or as a morphism
\begin{equation}
  \label{eq:B-compatible-HPD-stratification2}
  \theta^\cB:p^*_1M\to p^*_0M\ctens_{p^*_0\cB}\cP_{\cB/\cS,(m)}  
\end{equation}
linear for the right $\cP_{\cB/\cS,(m)}$-structure. The analogue of
\ref{eq:B-cocycle-condition} is the commutative diagram
\begin{equation}
  \label{eq:B-cocycle-condition-HPD}
  \xymatrix{
    M\ar[r]^{\theta^\cB}\ar[d]_{\theta^\cB}
    &M\ctensu\cB\cP_{\cB/\cS,(m)}\ar[d]^{\delta_{\cB,(m)}}\\
    M\ctensu\cB\cP_{\cB/\cS,(m)}
    \ar[r]^{\theta^\cB_n}
    &M\ctensu\cB\cP_{\cB/\cS,(m)}\ctensu\cB\cP_{\cB/\cS,(m)}
  }
\end{equation}
where $\delta_{\cB,(m)}$ is \ref{eq:delta^B3}. But we can also
linearize \ref{eq:B-cocycle-condition-HPD} and use the isomorphism 
\ref{eq:P_B/S(2)} and the maps $d^\cB_{ij}$; then the cocycle
condition takes the usual form of an equality
\begin{equation}
  \label{eq:B-compatible-cocycle-condition}
  (d^\cB_{02})_*(\chi^\cB)
  =(d^\cB_{01})_*(\chi^\cB)\circ (d^\cB_{12})_*(\chi^\cB)
\end{equation}
of $\cP_{\cB/\cS,(m)}(2)$-modules, where $(d^\cB_{ij})_*(\chi^\cB)$
denotes the extension of scalars of $\chi^\cB$ by $d^\cB_{ij}$
(corresponding to the pullbacks by the projections for the nonexistent
formal schemes corresponding the the $\O_\cX$-algebras
$\cP_{\cB/\cS,(m)}$ and $\cP_{\cB/\cS,(m)}(2)$).

The following proposition is proven in exactly the same way in
\cite[Prop. 3.1.3]{berthelot:1990}. 

\begin{prop}\label{prop:Dm-and-Dmhat-coherent}
  Suppose $M$ is a coherent $\niv{\D}{m}_{\cB/\cS}$-module, that is
  coherent as a $\cB$-module. Then $M$ is coherent as a
  $\niv{\hD}{m}_{\cB/\cS}$-module, and the canonical homomorphism
  \begin{displaymath}
    M\to\niv{\hD}{m}_{\cB/\cS}\tens_{\niv{\D}{m}_{\cB/\cS}}M
  \end{displaymath}
  is an isomorphism. \nodemo
\end{prop}

\subsubsection{Change of coefficient ring.}
\label{sec:change-of-coefficient-ring}

Suppose $\cC$ is a second $\O_\cX$-algebra with a left
$\niv{\D}{m}_{\cX/\cS}$-module structure compatible with its
$\O_\cX$-algebra structure. If $\cB\to\cC$ is an $\O_\cX$-algebra
homomorphism linear for the $\niv{\D}{m}_{\cX/\cS}$-module structures,
there is an obvious ring homomorphism
$\niv{\D}{m}_{\cB/\cS}\to\niv{\D}{m}_{\cC/\cS}$ inducing a canonical
$\niv{\D}{m}_{\cX/\cS}$-linear and $\cC$-linear isomorphism
\begin{equation}
  \label{eq:change-of-coefficient-ring}
  \cC\tens_\cB\niv{\D}{m}_{\cB/\cS}\to\niv{\D}{m}_{\cC/\cS}.
\end{equation}
If $M$ is a left $\niv{\D}{m}_{\cB/\cS}$-module, the transitivity of
tensor shows that there is an isomorphism
\begin{equation}
  \label{eq:change-of-coefficient-ring2}
  \cC\tens_\cB M\isom\niv{\D}{m}_{\cC/\cS}\tens_{\niv{\D}{m}_{\cB/\cS}}M.
\end{equation}
In particular, if $M$ is a coherent $\niv{\D}{m}_{\cB/\cS}$-module
then $\cC\tens_\cB M$ is a coherent $\niv{\D}{m}_{\cC/\cS}$-module. 

If $\cB$ and $\cC$ satisfy conditions
\ref{enum:coeff-ring-B1}-\ref{enum:coeff-ring-B2}, the same argument
can be applied to $M/JM$ for any $m$-bilateralising ideal
$J\subset\O_\cX$, and passing to the inverse limit yields an
isomorphism
\begin{equation}
  \label{eq:change-of-coefficient-ring-hat}
  \cC\ctens_\cB\niv{\D}{m}_{\cB/\cS}\to\niv{\hD}{m}_{\cC/\cS}.
\end{equation}
and, for any coherent left $\niv{\hD}{m}_{\cB/\cS}$-module $M$, an
isomorphism
\begin{equation}
  \label{eq:change-of-coefficient-ring2-hat}
  \cC\ctens_\cB M\isom\niv{\hD}{m}_{\cC/\cS}\tens_{\niv{\hD}{m}_{\cB/\cS}}M.
\end{equation}
Thus if $M$ is a coherent $\niv{\hD}{m}_{\cB/\cS}$-module then
$\cC\tens_\cB M$ is a coherent $\niv{\hD}{m}_{\cC/\cS}$-module.

\begin{prop}\label{prop:quasi-nilpotence-and-base-change}
  Suppose $\cB$ (resp. $\cC$) is an $\O_\cX$-algebra with a compatible
  $\niv{\hD}{m'}_{\cX/\cS}$-module structure (resp. a compatible
  $\niv{\hD}{m}_{\cX/\cS}$-module structure), $m'>m$ and $\cB\to\cC$
  is a $\niv{\hD}{m}_{\cX/\cS}$-linear homomomorphism of
  $\O_\cX$-algebras. If $M$ is a coherent
  $\niv{\hD}{m'}_{\cB/\cS}$-module and $\cC$ is topologically
  quasi-nilpotent as a $\niv{\hD}{m}_{\cX/\cS}$-module,
  $\cC\ctens_\cB M$ is topologically quasi-nilpotent for its induced
  $\niv{\hD}{m}_{\cX/\cS}$-module structure.
\end{prop}
\begin{demo}
  We may work locally, so suppose $M$ is generated as a
  $\niv{\hD}{m'}_{\cB/\cS}$-module by $a_1,\ldots,a_n$. Since $m'>m$,
  the formula \ref{eq:Dm-change-of-level-explicit} shows that
  $\dpabniv{\d}{K}{m}(a_i)\to0$ for $|K|\to\infty$. Since $\cC$ is a
  quasi-nilpotent $\niv{\hD}{m}_{\cX/\cS}$-module the Leibnitz formula
  \ref{eq:Leibnitz-B-modules} shows that $\dpabniv{\d}{K}{m}(x)\to0$
  for any local section of $\cC\ctens_\cB M$, and we are done.
\end{demo}

\subsection{The isogeny category.}
\label{sec:isogeny-category}

To extend these results to coherent $\niv{\D}{m}_{\cB/\cS\bQ}$-modules
we use the following theorem of Ogus, which was proven by him when
$\cB=\O_\cX$ and $m=0$, but the argument in the general case is the
same; c.f. also \cite[Prop. 3.1.2]{berthelot:1990} for the case
$\cB=\O_\cX$ and $m\ge0$.

\begin{prop}\label{prop:Ogus-integrality-thm}
  Suppose $\cB$ is an $\O_\cX$-algebra with a compatible
  $\niv{\hD}{m}_{\cB/\cS}$-module structure satisfying conditions
  \ref{sec:coefficients}.5--7. If $M$ is a $\cB_\bQ$-coherent left
  $\niv{\hD}{m}_{\cB/\cS\bQ}$-module, then locally on $\cX$ there is a
  $\cB$-coherent left $\niv{\hD}{m}_{\cB/\cS}$-module $M^0$ such that
  $(M^0)_\bQ\simeq M$.
\end{prop}
\begin{demo}
  Let $M'$ be any coherent $\cB$-module such that $M'_\bQ=M$.  Pick
  an open affine on which $\cX/\cS$ is parallelizable and choose local
  coordinates $x_1,\ldots,x_d$; then $\cP^n_{\cX/\cS,(m)}$ is free on
  the corresponding $\dpbrniv{\xi}{K}{m}$. For $n\ge0$ the map
  \begin{displaymath}
    \theta_n:M\to M\tens_{\O_\cX}\cP^n_{\cX/\cS,(m)}
  \end{displaymath}
  arising from the $m$-PD-stratification of $M$ is
  \begin{displaymath}
    \theta_n(x)=\sum_{|K|\le n}\dpabniv{\d}{K}{m}(x)\tens\dpbrniv{\xi}{K}{m}.
  \end{displaymath}
  We will show that
  \begin{displaymath}
    M^0=\bigcap_{n\ge0}\theta_n^{-1}(M'\tens_{\O_\cX}\cP^n_{\cX/\cS,(m)})
  \end{displaymath}
  satisifies the conditions of the proposition.  From the definition
  we see that $x\in M^0$ if and only if $\dpabniv{\d}{K}{m}(x)\in M'$
  for all $K$, and thus $M^0$ is preserved by the
  $\dpabniv{\d}{K}{m}$. Since the $m$-PD-stratification of $M$ is
  semilinear with respect to the $m$-PD-stratification of $\cB$, $M^0$
  is a $\cB$-submodule of $M'$; this is most easily seen from
  condition \ref{eq:Leibnitz-B-modules}. Clearly $M^0_\bQ\simeq M$,
  and $M^0\sset M'$ since $\theta_0(M')=M'$. Thus $M^0$ is
  $\cB$-coherent, and by \ref{sec:coefficients}.5 $M^0$ is complete
  for the adic topology of $\cB$ (induced by the adic topology of
  $\O_\cX$). Since $\dpabniv{\d}{K}{m}(M^0)\sset M^0$ for all $K$ it
  follows that $M^0$ is a $\niv{\hD}{m}_{\cX/\cS}$-submodule of
  $M$. Finally, since the $\cB$-module and
  $\niv{\hD}{m}_{\cB/\cS}$-module structures of $M$ are compatible
  this will also be the case for $M^0$, and thus $M^0$ is a
  $\niv{\hD}{m}_{\cB/\cS}$-module.
\end{demo}

If $\cX$ is quasicompact, the $M$ in the proposition can be chosen
globally, but we will not need this refinement. An immediate
consequence of Ogus's theorem is a version of
\ref{prop:Dm-and-Dmhat-coherent} for
$\niv{\D}{m}_{\cB/\cS\bQ}$-modules:

\begin{prop}\label{prop:Dm-and-Dmhat-coherent-Q}
  Suppose $M$ is a left $\niv{\D}{m}_{\cB/\cS\bQ}$-module that is
  coherent as a $\cB_\bQ$-module. Then $M$ is coherent as a left
  $\niv{\D}{m}_{\cB/\cS\bQ}$-module, and the canonical homomorphism
  \begin{displaymath}
    M\to\niv{\hD}{m}_{\cB/\cS\bQ}\tens_{\niv{\D}{m}_{\cB/\cS\bQ}}M
  \end{displaymath}
  is an isomorphism. \nodemo
\end{prop}

The following is also useful when going back and forth between
$\niv{\D}{m}_{\cB/\cS}$-modules and
$\niv{\D}{m}_{\cB/\cS\bQ}$-modules:

\begin{lemma}\label{lemma:choice-of-integral-model}
  Let $A$ be a noetherian topological $\bZ_p$-algebra (not necessarily
  commutative) and $A\to B$ a ring homomorphism. Suppose $A$ has a
  ideal $J$ centralising in both $A$ and $B$, and such that $A$ and
  $B$ both have the $J$-adic topology. Let $f:M\to M'$ be a
  homomorphism of finitely generated left $A$-modules. If the kernel
  and cokernel of $f$ are $p$-torsion, $f$ induces an isomorphism
  $(B\ctens_AM)_\bQ\isom(B\ctens_AM')_\bQ$.
\end{lemma}
\begin{demo}
  It suffices to treat the cases in which $f$ is surjective or
  injective. If $f$ is injective we define $M''$ by the exactness of
  \begin{displaymath}
    0\to M\Xto{f} M'\to M''\to 0
  \end{displaymath}
  and there is an exact sequence
  \begin{displaymath}
    0\to K\to B\tens_AM\to B\tens_AM'\to B\tens_AM''\to 0
  \end{displaymath}
  where $K$ is a quotient of $\Tor^1_A(B,M'')$. Since $M''$ is
  $p$-torsion and a finitely generated $A$-module, $K$ and
  $B\tens_AM''$ are $p$-torsion finitely generated $B$-modules, and
  thus annihilated by some power of $p$. By the hypotheses on the
  topologies of $A$ and $B$, the completion
  \begin{displaymath}
    0\to\hat K\to B\ctens_AM\to B\ctens_AM'\to B\ctens_AM''\to 0
  \end{displaymath}
  of the previous exact sequence is exact. Furthermore $\hat K$ and
  $B\ctens_AM''$ are annihilated by some power of $p$, and the
  assertion follows. The case where $f$ is surjective is similar.
\end{demo}



The following lemma is evident:

\begin{lemma}\label{lemma:quasi-nilpotent-Q}
  Suppose $\cB$ is an $\O_\cX$-algebra with a compatible
  $m$-HPD-stra\-ti\-fi\-cation and $M$ is a \emph{coherent} left
  $\niv{\hD}{m}_{\cB/\cS\bQ}$-module. The following are equivalent:
  \begin{enumerate}
  \item Any finitely generated $\niv{\hD}{m}_{\cB/\cS}$-submodule of
    $M$ is topologically quasi-nilpotent.
  \item There is a topologically quasi-nilpotent left
    $\niv{\hD}{m}_{\cB/\cS}$-module $M^0$ and an isomorphism
    $M^0_\bQ\simeq M$ of $\niv{\hD}{m}_{\cB/\cS}$-modules.
  \end{enumerate}\nodemo
\end{lemma}

\begin{defn}\label{defn:topologically-quasi-nilpotent-isogeny}
  Suppose $\cB$ is as in lemma \ref{lemma:quasi-nilpotent-Q}. A
  coherent left $\niv{\hD}{m}_{\cB/\cS\bQ}$ $M$ is topologically
  quasi-nilpotent if it satisfies the equvalent conditions of
  lemma \ref{lemma:quasi-nilpotent-Q}.
\end{defn}

If $M$ is a coherent $\niv{\hD}{m}_{\cB/\cS\bQ}$-module we define the
$\cC\ctens_\cB M$ as follows: choose locally a coherent
$\niv{\hD}{m}_{\cB/\cS}$-submodule $M^0$ of $M$ such that
$M^0_\bQ\simeq M$; then
\begin{equation}
  \label{eq:rigid-analytic-base-change}
  \cC\ctens_\cB M:=(\cC\ctens_\cB M^0)_\bQ
\end{equation}
where in the right hand side we have the usual completed tensor
product. Lemma \ref{lemma:choice-of-integral-model} shows that the
left hand side of \ref{eq:rigid-analytic-base-change} is independent
of the choice of $M^0$, so the definition is justified. Note that if
$M$ is coherent as a $\cB_\bQ$-module we may choose an $M^0$ that is
coherent as a $\cB$-module. In any case the resulting
$\niv{\hD}{m}_{\cC/\cS}$-module is coherent.

From \ref{prop:quasi-nilpotence-and-base-change}, lemma
\ref{lemma:quasi-nilpotent-Q} and definition
\ref{defn:topologically-quasi-nilpotent-isogeny} we find:

\begin{cor}\label{cor:quasi-nilpotence-and-base-change}
  Let $\cB$ and $\cC$ be as in proposition
  \ref{prop:quasi-nilpotence-and-base-change}, and let $M$ be a
  coherent left $\niv{\hD}{m'}_{\cB/\cS\bQ}$-module. If $m'>m$,
  $\cC\ctens_\cB M$ is a quasi-nilpotent left
  $\niv{\hD}{m}_{\cC/\cS\bQ}$-module.
\end{cor}

In particular:

\begin{cor}\label{cor:quasi-nilpotence-change-of-level}
  Suppose $\cB$ is an $\O_\cX$ with a compatible
  $m$-HPD-stratification satisfying
  \ref{enum:coeff-ring-B1}-\ref{enum:coeff-ring-B2} and $M$ is a
  coherent left $\niv{\hD}{m+1}_{\cB/\cS\bQ}$-module. Then $M$ with
  the induced $\niv{\hD}{m}_{\cB/\cS\bQ}$-module structure is
  topologically quasi-nilpotent.  
\end{cor}

\subsection{Descent by Frobenius.}
\label{sec:Frobenius-descent}

Suppose now $\cS$ has the $m$-PD-structure $(\fa,\fb,\alpha)$ and
$p\in\fa$. As before we set $S_0=V(\fa)$, set $q=p^s$ and denote by
$F_{S_0}$ the $q$th power Frobenius of $S_0$. For any formal
$\cS$-scheme $\cX$ we again set $X_0=V(\fa\O_\cX)$ and denote by
$F_{X_0}$ the $q$th power Frobenius; then the relative Frobenius
$F_{X_0/S_0}:X_0\to\niv{X_0}{q}$ is defined, and we denote by
$W_{X_0/S_0}:\niv{X_0}{q}\to X_0$ the canonical projection.  Suppose
$F_{X_0/S_0}$ lifts to a morphism $F:\cX\to\cX'$. The argument of
\cite[Prop. 2.2.2]{berthelot:2000} can be used as is to show:

\begin{prop}\label{prop:raising-the-level}
  For any left $\niv{\D}{m}_{\cX'/\cS}$-module $M$, $F^*M$ has a
  canonical and functorial left $\niv{\D}{m+s}_{\cX/\cS}$-module
  structure restricting to the $\niv{\D}{m}_{\cX/\cS}$-module
  structure derived from base change. If $M$ is a topologically
  quasi-nilpotent $\niv{\D}{m}_{\cX'/\cS}$-module, $F^*M$ is a
  topologically quasi-nilpotent
  $\niv{\D}{m+s}_{\cX/\cS}$-module.\nodemo
\end{prop}

An important point in the proof is that $F$ is flat, a fact used in
the proofs of \cite[Lemme 2.3.2]{berthelot:2000} and \cite[Lemme
2.3.3]{berthelot:2000}. This is true in the present case, since $F$ is
a lifting of $F_{X_0/S_0}:X_0\to\niv{X_0}{q}$ which is flat by
proposition \ref{prop:flatness-of-Frobenius}.

To prove the descent theorem we need to make the further assumption
that $F$ is finite, a fact needed in the argument of
\cite[Th\'eor\`eme 2.3.6]{berthelot:2000}. This will be true if
$F_{X_0/S_0}$ is finite, but so far we only know this when
$X_0\to S_0$ is formally of finite type, again by proposition
\ref{prop:flatness-of-Frobenius}. We therefore state the Frobenius
descent theorem for $\niv{\D}{m}$-modules as follows; no changes are
needed in adapting the argument of \cite[\S2.3]{berthelot:2000}.

\begin{thm}\label{thm:Frobenius-descent-D}
  Suppose $f:\cX\to\cS$ is quasi-smooth and the relative Frobenius
  $F_{X_0/S_0}$ is finite. The functor $F^*$ induces an equivalence of
  the category of left $\niv{\D}{m}_{\cX'/\cS}$-modules with the
  category of left $\niv{\D}{m+s}_{\cX/\cS}$-modules.\nodemo
\end{thm}

The Frobenius descent theorem for $\niv{\hD}{m}$-modules is an
immediate consequence. If $J'\sset\O_{\cX'}$ is open and
$\niv{\D}{m}_{\cX'/\cS}$-bilateralising, $J=F^*J'\sset\O_\cX$ is open
and $\niv{\D}{m}_{\cX'/\cS}$-bilateralising. From the equivalence of
categories in the last paragraph we deduce, again with the assumption
that $F$ is finite and flat, that the category of left
$\niv{\D}{m}_{X'_{J'}/\cS}$-modules is equivalent to the category of
left $\niv{\D}{m}_{\cX/\cS,J}$-modules. If $J'$ is an ideal of
definition, so is $J$; this is because $F_{X_0/S_0}$ is a
homeomorphism.Furthermore, as $J'\sset\O_{\cX'}$ runs through a
cofinal set of $\niv{\D}{m}_{\cX'/\cS}$-bilateralising ideals of
definition, $F^*J$ runs through a cofinal set of
$\niv{\D}{m+s}_{\cX/\cS}$-bilateralising ideals of definition. From
this and the previous discussion we get the following:

\begin{thm}\label{thm:Frobenius-descent-Dhat}
  Suppose $\cX\to\cS$ is quasi-smooth and the relative Frobenius
  $F_{X_0/S_0}$ lifts to an $\cS$-morphism $F:\cX\to\cX'$. For any
  left $\niv{\hD}{m}_{\cX'/\cS}$-module $M$, $F^*M$ has a canonical
  left $\niv{\hD}{m+s}_{\cX/\cS}$-module structure. If $F_{X_0/S_0}$
  is finite, this induces an equivalence of the category of left
  $\niv{\hD}{m}_{\cX'/\cS}$-modules with the category of left
  $\niv{\hD}{m+s}_{\cX/\cS}$-modules.
\end{thm}

Suppose $F':\cX\to\cX'$ is a second lifting of $F_{X_0/S_0}$. If the
$m$-PD-structure $(\fa,\fb,\alpha)$ of $\cS$ is $m$-PD-nilpotent, the
isomorphism $\tau_{F,F'}$ of proposition
\ref{prop:invariance-of-base-change1} is
$\niv{\D}{m+s}_{\cX/\cS}$-linear, the argument being the same as
\cite[Prop. 2.2.5]{berthelot:1996}.

The Frobenius descent theorem extends to the case of coefficient
rings. Let $\cB$ be an $\O_{\cX'}$-algebra with a compatible left
$\niv{\hD}{m}_{\cX'/\cS}$-module structure. Then $F^*\cB$ has a left
$\niv{\hD}{m+s}_{\cX/\cS}$-module structure compatible with its
$\O_\cX$-algebra structure. Furthermore if $\cB$ satisfies the
assumptions \ref{sec:coefficients}.3--5 relative to $\cX'$, the
$\O_\cX$-algebra $F^*\cB$ satisfies these same conditions relative to
$\cX$. In these cases the Frobenius descent theorem can be stated as
follows:

\begin{thm}\label{thm:Frobenius-descent-with-coefficients}
  Let $\cB$ be an $\O_{\cX'}$-algebra with a left
  $\niv{\D}{m}_{\cX'/\cS}$-module structure compatible with its
  $\O_{\cX'}$-algebra structure. For any left
  $\niv{\D}{m}_{\cB/\cS}$-module $M$, $F^*$ has a canonical left
  $\niv{\D}{m+s}_{F^*\cB/\cS}$-module structure, and if $F_{X_0/S_0}$
  is finite this construction yields an equivalence of the category of
  left $\niv{\D}{m}_{\cB/\cS}$-modules with the category of left
  $\niv{\D}{m+s}_{F^*\cB/\cS}$-modules. If in addition $\cB$ satisfies
  the conditions \ref{sec:coefficients}.3--5, $F^*$ induces an
  equivalence of the category of left $\niv{\hD}{m}_{\cB/\cS}$-modules
  with the category of left $\niv{\hD}{m+s}_{F^*\cB/\cS}$-modules
\end{thm}
